\documentclass[11pt,namelimits,sumlimits,a4paper]{amsart}

\usepackage{amssymb}
\usepackage{amsmath}
\usepackage[mathscr]{eucal}
\usepackage{slashed}
\usepackage{amsthm}
\usepackage[latin1]{inputenc}
\usepackage[T1]{fontenc}
\usepackage[UKenglish]{babel}
\usepackage{amsfonts}
\usepackage{amscd}
\usepackage{latexsym}
\usepackage{cite}
\usepackage{comment}
\usepackage{times,psfrag}
\usepackage{mathtools}
\usepackage{bm}
\usepackage{esint}
\usepackage{layout}
\usepackage{lmodern}
\usepackage{cite}
\usepackage{comment}
\usepackage{color}
\usepackage{hyperref}

\numberwithin{equation}{section}
\newtheorem{thm}[equation]{Theorem}
\newtheorem{lemma}[equation]{Lemma}
\newtheorem{prop}[equation]{Proposition}

\theoremstyle{definition}
\newtheorem{definition}[equation]{Definition}
\newtheorem{example}[equation]{Example}

\theoremstyle{remark}
\newtheorem{remark}[equation]{Remark}
\newtheorem*{remark*}{Remark}

\newcommand{\ie}{\emph{i.e.} }

\newcommand{\cf}{\emph{cf.} }

\newcommand{\beq}{\begin{equation}}
\newcommand{\eeq}{\end{equation}}
\newcommand{\bea}{\begin{eqnarray}}
\newcommand{\eea}{\end{eqnarray}}

\newcommand{\C}{\mathbb{C}}

\newcommand{\R}{\mathbb{R}}

\newcommand{\Z}{\mathbb{Z}}
\newcommand{\N}{\mathbb{N}}

\newcommand{\T}{\mathbb{T}}

\newcommand{\Sph}{\mathbb{S}}
\newcommand{\ra}{\rightarrow}

\newcommand{\dvol}{\operatorname{dv}}

\newcommand{\ind}{\operatorname{ind}}

\newcommand{\Ric}{\operatorname{Ric}}
\newcommand{\Trace}{\operatorname{Trace}}

\newcommand{\Ad}{\textrm{Ad}}
\newcommand{\ad}{\textrm{ad}}
\newcommand{\Lie}[1]{\mathfrak{#1}}

\newcommand{\ext}{\text{ext}}

\newcommand{\Aut}{\text{Aut}}

\newcommand{\RS}{\R^{2} \times \Sph ^{1}}

\def\co{\colon\thinspace}

\begin{document}

\title{Deformation theory of periodic monopoles (with singularities)}

\author[L.~Foscolo]{Lorenzo~Foscolo}
\address{Mathematics Department, State University of New York at Stony Brook}
\email{lorenzo.foscolo@stonybrook.edu}

\date{\today}

\begin{abstract}
In \cite{Cherkis:Kapustin:1} and \cite{Cherkis:Kapustin:2} Cherkis and Kapustin introduced periodic monopoles (with singularities), \ie monopoles on $\R ^{2} \times \Sph ^{1}$ possibly singular at a finite collection of points. In this paper we show that for generic choices of parameters the moduli spaces of periodic monopoles (with singularities) are either empty or smooth hyperk\"ahler manifolds. Furthermore, we prove an index theorem and therefore compute the dimension of the moduli spaces.
\end{abstract}

\maketitle

\section{Introduction}\label{sec:Introduction}

Let $(X,g)$ be an oriented Riemannian $3$--manifold and $P \ra X$ a principal $G$--bundle, where $G$ is a compact Lie group. Consider the product $X \times \R _{s}$ endowed with the product metric, the volume form $ds \wedge \dvol _{g}$ and the pulled-back $G$--bundle $\hat{P}$. An anti-self-dual (ASD) connection (or instanton) on $\hat{P}$ is a connection $\hat{A}$ such that $\ast F_{\hat{A}}+F_{\hat{A}}=0$. If $\hat{A}$ is $\R$--invariant one can write $\hat{A}=A+\Phi \otimes ds$, for a connection $A$ on $P \ra X$ and a section $\Phi$ of the adjoint bundle $\ad (P)$. Monopoles on $X$ are pairs $(A,\Phi)$ such that $\hat{A}$ is an $\R$--invariant ASD connection on $X \times \R$. Working directly in $3$--dimensions we have the following defintion. 
\begin{definition}
\emph{(Magnetic) monopoles} are solutions $(A,\Phi )$ to the \emph{Bogomolny equation}
\begin{equation}\label{eqn:Bogomolny}
\ast F_{A}=d_{A}\Phi.
\end{equation}
Here $\ast$ is the Hodge star operator of $(X,g)$; $F_{A}$ is the curvature of the connection $A$ and $\Phi$ is called the \emph{Higgs field}. The \emph{moduli space} of monopoles on $P \ra X$ is the space of equivalence classes of solutions to \eqref{eqn:Bogomolny} with respect to the action of the gauge group $\Aut (P)$.
\end{definition}

An immediate consequence of equation (\ref{eqn:Bogomolny}) and the Bianchi identity is
\begin{equation}\label{eqn:Harmonic:Higgs:Field}
d_{A}^{\ast}d_{A}\Phi =0.
\end{equation}
In particular, when $X$ is compact smooth monopoles coincide with reducible (if $|\Phi| \neq 0$) flat connections. In order to find non-trivial solutions to (\ref{eqn:Bogomolny}) one has to consider a non-compact base manifold $X$, in the sense that either $X$ is complete or we allow for singularities of the fields $(A,\Phi)$, or a combination of the two possibilities, as in this paper.

The classical case of smooth monopoles on $\R ^{3}$ and the rich geometric properties of their moduli spaces have been investigated from many different points of view; a standard reference is Atiyah and Hitchin's book \cite{Atiyah:Hitchin}. Monopoles with and without singularities have also been studied on $3$--manifolds $X$ with different geometries: hyperbolic monopoles were introduced by Atiyah \cite{Atiyah:Hyperbolic:Monopoles}; Braam reduced the study of monopoles on an asymptotically hyperbolic manifold $X$ to that of $\Sph ^{1}$--invariant ASD connections on a conformal compactification \cite{Braam}; partial results were established by Floer \cite{Floer:Monopoles:1,Floer:Monopoles:2} for asymptotically Euclidean $X$; more recently, Kottke initiated the study of monopoles on asymptotically conical $3$--manifolds \cite{Kottke}. Monopoles with singularities were first considered by Kronheimer \cite{Kronheimer:Thesis}; the dimension of the moduli space of singular monopoles over a compact manifold $X$ was computed by Pauly \cite{Pauly}; Charbonneau and Hurtubise considered monopoles with singularities on the product of a compact Riemann surface with a circle \cite{Charbonneau:Hurtubise}.

An important feature of the moduli spaces of monopoles on $\R ^{3}$ is that they are hyperk\"aler manifolds by virtue of an infinite dimensional hyperk\"ahler quotient. In the lowest non-trivial dimension, the Atiyah--Hitchin manifold, \ie the moduli space of centred charge $2$ $SU(2)$ monopoles on $\R ^{3}$ (or its double cover) is a complete hyperk\"ahler $4$--manifold with finite $L^{2}$--norm of the curvature, a so-called \emph{gravitational instanton}, with an interesting asymptotic geometry: the volume of large geodesic balls of radius $r$ grows like $r^{3}$; the complement of a large ball is a circle bundle over $\R ^{3}/\Z_{2}$ and the metric is asymptotically adapted to this circle fibration. We say that the Atiyah--Hitchin metric is an ALF gravitational instanton.

Pursuing the idea that moduli spaces of solutions to dimensional reductions of the Yang--Mills ASD equations on $\R ^{4}$ are ``a natural place to look for gravitational instantons'' \cite{Cherkis:Talk}, in \cite{Cherkis:Kapustin:1}, \cite{Cherkis:Kapustin:3} and \cite{Cherkis:Kapustin:2} Cherkis and Kapustin introduced the study of \emph{periodic monopoles}, \ie monopoles on $\R ^{2} \times \Sph ^{1}$, possibly with isolated singularities at a finite collection of points. They argued that, when $4$--dimensional, moduli spaces of periodic monopoles (with singularities) are gravitational instantons of type ALG: the volume of large balls grows quadratically and the metric is asymptotically adapted to a fibration by $2$--dimensional tori.

This paper addresses some of the foundational questions opened by Cherkis and Kapustin's work. The main results are summarised in the following theorem. 

\begin{thm}\label{thm:Main:Theorem}
For generic choices of the parameters defining the boundary conditions, the moduli space $\mathcal{M}_{n,k}$ of charge $k$ $SO(3)$ periodic monopoles with $n$ isolated singularities is a smooth hyperk\"ahler manifold of dimension $4k-4$, provided it is not empty.
\end{thm}

Here the charge is a certain topological invariant of a monopole, \cf Definition \ref{def:Boundary:Conditions:SO(3)} for details. In \cite{Gluing} we construct solutions to \eqref{eqn:Bogomolny} on $\R ^{2} \times \Sph ^{1}$ by gluing methods, showing that $\mathcal{M}_{n,k}$ is indeed non-empty.

\subsection*{Plan of the paper}

In Section \ref{sec:Preliminaries} we introduce formal aspects of the construction of the moduli spaces $\mathcal{M}_{n,k}$ and fix some notation. In Section \ref{sec:Periodic Dirac} we define periodic Dirac monopoles, \ie solutions to \eqref{eqn:Bogomolny} on $\RS$ with structure group $U(1)$ and one isolated singularity. Following \cite{Cherkis:Kapustin:1} and \cite{Cherkis:Kapustin:2}, we then use this material to define boundary conditions for periodic monopoles with non-abelian structure group $G=SO(3)$.

Sections \ref{sec:Singularity} and \ref{sec:Big:End} deal with the local analysis in a neighbourhood of the singularities and on the big end of $\RS$: we introduce weighted Sobolev spaces, prove embedding and multiplication results and study the mapping properties of the relevant operators. These  analytic results are applied in Section \ref{sec:Construction:Moduli:Space} to prove that the moduli spaces $\mathcal{M}_{n,k}$ are smooth hyperk\"ahler manifolds (when non-empty) provided there are no reducible solutions of the Bogomolny equation satisfying the given boundary conditions.

The final section contains the proof of the dimension formula, \ie the computation of the index of a certain Dirac-type operator. No index theorem available in the literature applies to the situation at hand and we give a geometric proof of the index formula based on the excision principle.

\subsection*{Aknowledgments}
The results of this paper are part of the author's Ph.D. thesis at Imperial College London. He wishes to thank his supervisor Mark Haskins for his continuous support. Olivier Biquard guided early stages of this project; we thank him for suggesting us this problem. The author is grateful to Simon Donaldson and Michael Singer for their careful comments on an early version of this paper. The paper was completed as the author was a Simons Instructor at SUNY Stony Brook.

\section{Preliminaries}\label{sec:Preliminaries}

In this section, whose purpose is mainly to fix the notation, we recall formal aspects of the deformation theory of monopoles. In particular, we introduce the relevant elliptic operators and state Weitzenb\"ock formulas that will be used throughout the paper.

Let $X$ be a non-compact oriented $3$--manifold and $P \ra X$ a principal $G$--bundle. Denote by $\mathcal{C}$ the infinite dimensional space of smooth pairs $c=(A,\Phi)$, where $A$ is a connection on $P \ra X$ and $\Phi \in \Omega ^{0}(X;\ad P)$ a Higgs field. Since $X$ is not compact, elements $c \in \mathcal{C}$ have to satisfy appropriate boundary conditions, which we suppose to be included in the definition of $\mathcal{C}$. The space $\mathcal{C}$ is an affine space. The underlying vector space is the space of section $\Omega (X;\ad P)=\Omega ^{1}(X;\ad P)\oplus\Omega ^{0}(X;\ad P)$ satisfying appropriate decay conditions. Let $\mathcal{G}$ be the group of bounded smooth sections of $\Aut(P)$ which preserve the chosen boundary conditions. Here $g \in Aut(P)$ acts on a pair $c=(A,\Phi) \in \mathcal{C}$ by $c \mapsto c + (d_{1}g)g^{-1}$, where
\begin{equation}\label{eqn:Linearisation:Gauge:Action}
d_{1}g = -\left( d_{A}g, [\Phi,g] \right) \in \Omega (X; \ad P).
\end{equation}

Consider the gauge-equivariant map $\Psi\co \mathcal{C} \ra \Omega ^{1}(X;\ad P)$ defined by $(A,\Phi )\mapsto \ast F_{A}-d_{A}\Phi$. By fixing a base point $c=(A,\Phi )\in\mathcal{C}$ we write $\Psi (A+a,\Phi +\psi )=\Psi (c)+d_{2}(a,\psi )+(a,\psi )\cdot (a, \psi )$ for all $(a,\psi )\in\Omega (X; \ad P)$. The linearisation $d_{2}$ of $\Psi$ at $c$ and the quadratic term are defined by:
\begin{equation}\label{eqn:Linearisation:Bogomolny}
d_{2}(a,\psi )=\ast d_{A}a-d_{A}\psi +[\Phi ,a]
\end{equation}
\begin{equation}\label{eqn:Quadratic:Term:Bogomolny}
(a,\psi )\cdot (a, \psi )=\ast [a,a]-[a,\psi ]
\end{equation}

The linearisation at $c$ of the action of $\mathcal{G}$ on $\mathcal{C}$ is the operator $d_{1}\co \Omega ^{0}(X;\ad\,P) \ra \Omega (X; \ad P)$ defined as in \eqref{eqn:Linearisation:Gauge:Action}. Couple $d_{2}$ with $d_{1}^{\ast}$ to obtain an elliptic operator
\begin{equation}\label{eqn:Dirac:Operator}
D=d_{2}\oplus d_{1}^{\ast }\co \Omega (X; \ad P)\longrightarrow \Omega (X; \ad P).
\end{equation}

The moduli space $\mathcal{M}$ of monopoles in $\mathcal{C}$ is $\mathcal{M}=\Psi ^{-1}(0)/\mathcal{G}$. Suppose that $c=(A,\varphi )$ is a solution to the Bogomolny equation and consider the elliptic complex
\begin{equation}\label{eqn:Deformation:Complex}
\Omega ^{0}(X;\ad\, P) \stackrel{d_{1}}{\longrightarrow} \Omega (X; \ad P)\stackrel{d_{2}}{\longrightarrow}\Omega ^{1}(X;\ad\, P)
\end{equation}
(this is a complex precisely when $\Psi (A,\Phi)=0$). Standard theory \cite[Chapter 4]{Donaldson:Kronheimer} implies that $\mathcal{M}$ is a smooth manifold if---after choosing Sobolev completions of the spaces of $\ad(P)$--valued forms so that $\Psi$ and the action of gauge transformations $\mathcal{G} \times \mathcal{C} \ra \mathcal{C}$ extend to smooth maps of Banach spaces and \eqref{eqn:Deformation:Complex} is a Fredholm complex---the cohomology groups of \eqref{eqn:Deformation:Complex} in degree $0$ and $2$ vanish. Then the tangent space $T_{[c]}\mathcal{M}$ at the point $[c]$ is identified with $\ker D_{c}$, \ie the cohomology of \eqref{eqn:Deformation:Complex} in degree $1$.

We can interpret $D$ as a twisted Dirac operator on $\Omega (X; \ad P)$. The Clifford multiplication of a $1$--form $\alpha$ and a $k$--form $\beta$ on $X$ is
\begin{equation}\label{eqn:Clifford:Multiplication}
\gamma (\alpha)\beta = \alpha \wedge \beta - \alpha ^{\sharp} \lrcorner\, \beta
\end{equation}
Define a twisted Dirac operator $\slashed{D}_{A}$ on $\Omega (X; \ad P)$ by
\begin{equation}\label{eqn:Dirac:Operator:1}
\Omega ^{1} \oplus \Omega ^{0} \xrightarrow{(\text{id},\ast)} \Omega ^{1} \oplus \Omega ^{3} \xrightarrow{\gamma\, \circ \nabla _{A}} \Omega ^{2} \oplus \Omega ^{0} \xrightarrow{(\ast, \text{id})} \Omega ^{1} \oplus \Omega ^{0}.
\end{equation}
The operator $D$ of \eqref{eqn:Dirac:Operator} is $D = \tau\slashed{D}_{A} + [\Phi ,\cdot\, ]$, where $\tau$ is a sign operator with $\tau = 1$ on $1$--forms and $\tau =-1$ on $0$--forms. From this point of view, the product \eqref{eqn:Quadratic:Term:Bogomolny} is the multiplication on $\Omega (X; \ad P)$ obtained combining Clifford multiplication of forms and the Lie bracket on $\ad(P)$. The formal $L^{2}$--adjoint of $D$ is $D^{\ast} = D - 2[\Phi, \cdot\, ]$ and we have the following Weitzenb\"ock formulas.

\begin{lemma}[see for example {\cite[Lemma 18]{Floer:Monopoles:2}}]\label{lem:Weitzenbock}
\begin{alignat*}{2}
DD^{\ast}=\nabla _{A}^{\ast}\nabla _{A}-\text{ad}(\Phi)^{2}+\Psi + \Ric
\qquad \text{ and } \qquad &
D^{\ast}D=DD^{\ast}+2d_{A}\Phi,
\end{alignat*}
where $\Psi = \ast F_{A} - d_{A}\Phi$.
\end{lemma}

As a final remark in this general setting, observe that if one fixes boundary conditions so that infinitesimal deformations are $L^{2}$--integrable, the $L^{2}$--product restricted to $\ker D$ defines a Riemannian metric on the moduli space $\mathcal{M}$. As in the Euclidean case, if $X=\RS$ this \emph{$L^{2}$--metric} is hyperk\"ahler by virtue of an infinite dimensional hyperk\"ahler quotient \cite{Atiyah:Hitchin}.

\section{Periodic Dirac monopole}\label{sec:Periodic Dirac}

When the structure group $G=U(1)$, the Bogomolny equation \eqref{eqn:Bogomolny} reduces to a linear equation. By \eqref{eqn:Harmonic:Higgs:Field} the Higgs field $\Phi$ is a harmonic function such that $\frac{\ast d\Phi}{2\pi i}$ represents the first Chern class of a line bundle. Global solutions are necessarily trivial; on $\R ^{3}$ they are given by pairs $(A,\Phi)=(0,v)$ while on $\R ^{2} \times \Sph ^{1}$ by $(A,\Phi)=(d+ib\, dt,v)$, where $v \in \R$ and $b \in \R /\Z$. We call such pairs flat (or vacuum) abelian monopoles. Non-trivial abelian solutions are obtained if one allows an isolated singularity.

\begin{definition}\label{def:Euclidean:Dirac:Monopole}
Fix a point $q \in \R ^{3}$ and let $H_{q}$ denote the radial extension of the Hopf line bundle to $\R ^{3} \setminus \{ q \}$. Fix $k \in \Z$ and $v \in \R$. The \emph{Euclidean Dirac monopole} of charge $k$ and mass $v$ with singularity at $q$ is the abelian monopole $(A, \Phi)$ on $H_{q}^{k}$, where
\[
\Phi =i\left( v-\frac{k}{2|x-q|} \right),
\]
$x \in \R ^{3}$, and $A$ is the $SO(3)$--invariant connection on $H_{q}^{k}$ with curvature $\ast d\Phi$.
\end{definition}

Periodic Dirac monopoles are defined in a similar way. Fix coordinates $(z,t)\in \C \times \R /2\pi \Z$ and a point $q=(z_{0},t_{0})\in \R ^{2} \times \Sph ^{1}$. Line bundles of a fixed degree on $(\RS) \setminus \{ q \}$ differ by tensoring by flat line bundles. We can distinguish connections with the same curvature by comparing their holonomy around loops $\gamma _{z}:=\{ z \} \times \Sph ^{1}_{t}$ for $z \neq z_{0}$. Set $\theta _{q}=\text{arg}(z-z_{0})$ and fix an origin in the circle parametrised by $\theta _{q}$. It follows from Remark \ref{rmk:Holonomy} below that the holonomy around $\gamma _{z}$ of a connection on a degree $k$ line bundle over $(\RS) \setminus \{ q \}$ is of the form $e^{-ik\theta _{q}}e^{-2\pi ib}$ for some $b \in \R /\Z$. Denote by $L_{q}$ the degree $1$ line bundle on $(\RS) \setminus \{ q \}$ with connection $A_{q}$ whose holonomy around $\gamma _{z}$ is $e^{-i\theta _{q}}$. Any line bundle of degree $1$ is of the form $L_{q} \otimes L_{b}$ for some flat line bundle $L_{b}$.

\begin{definition}\label{def:Periodic:Dirac:Monopole}
Fix a point $q \in \RS$. The \emph{periodic Dirac monopole} of charge $k \in \Z$, with singularity at $q$ and twisted by the flat line bundle $L_{v,b}$ for some $v \in \R$ and $b \in \R /\Z$ is the pair $(A, \Phi)$ on $L_{q}^{k} \otimes L_{v,b}$, where
\[
-i\Phi =v+kG_{q}
\]
and up to gauge transformations the connection $A= kA_{q}+ib\, dt$. Here $G_{q}$ defined in \eqref{eqn:Periodic:Dirac:Higgs:Field} below is a Green's function of $\RS$ with singularity at $q$.
\end{definition}

In the rest of the section we derive asymptotic expansions for the Green's function $G_{q}$ and the connection $A_{q}$, both at infinity and close to the singularity.

\subsection{The Green's function of $\R ^{2} \times \Sph ^{1}$}

By taking coordinates centred at $q \in \RS$, we can assume that the singularity is located at $q=0$. We use polar coordinates $z=re^{i\theta} \in \C$. Consider the series
\begin{equation}\label{eqn:Periodic:Dirac:Higgs:Field}
G(z,t)=-\frac{1}{2}\sum _{m\in\mathbb{Z}}{\left[ \frac{1}{\sqrt{r^{2}+(t-2m\pi )^{2}}}-a_{|m|}\right] },
\end{equation}
where
\begin{alignat*}{2}
a_{|m|}=\frac{1}{2|m|\pi } \text{ if } m\neq 0
 \qquad &
a_{0}=2\frac{\log{4\pi }-\gamma }{2\pi }
\end{alignat*}
($\gamma $ is the Euler--Mascheroni constant, $\gamma =\lim _{n\rightarrow\infty}\sum _{k=1}^{n}k^{-1}-\log{n}$).

\begin{lemma}\label{lem:Asymptotics:Periodic:Dirac:Higgs:Field}
The series (\ref{eqn:Periodic:Dirac:Higgs:Field}) converges uniformly on compact sets of $(\RS)\setminus \{0\}$ to a Green's function of $\RS$ with singularity at $0$.
\begin{itemize}
\item[(i)] Whenever $z\neq 0$, $G$ can be expressed as
\[
G(z,t)=\frac{1}{2\pi }\log {r}-\frac{1}{2\pi }\sum _{m\in\mathbb{Z}^{\ast}}{K_{0}(|m|r)e^{imt}},
\]
where $K_{0}$ is the second modified Bessel function.
\item[(ii)] There exists a constant $C_{1}>0$ such that
\[
\left| \nabla ^{k}\left( G(z,t)-\frac{1}{2\pi }\log {r}\right)\right| \leq C_{1}e^{-r}
\]
for all $r\geq 2$ and $k=0,1,2$.
\item[(iii)] There exists a constant $C_{2}>0$ such that 
\[
\left| \nabla ^{k}\left( G(z,t)-\frac{a_{0}}{2}+\frac{1}{2\rho }\right)\right| \leq C_{2}\rho ^{2-k}
\]
for all $(z,t)$ with $\rho=\sqrt{r^{2}+t^{2}}< \frac{\pi}{2}$ and $k=0,1,2$.
\end{itemize}
\proof
The convergence of \eqref{eqn:Periodic:Dirac:Higgs:Field}, the expansion in $(i)$ and the estimate in (ii) are proved in \cite[Lemma 3.1 (a),(b)]{Gross:Wilson}. $(iii)$ follows from the classical multipole expansion.
\qed
\end{lemma}

\subsection{The connection.}

Fix a constant $v \in \mathbb{R}$ and consider the Higgs field $\Phi =iv+iG$. The $2$--form $i\ast dG$ represents the curvature of a line bundle $L=L_{q}$ over $(\RS)\setminus\{ q\}$. A connection $A=A_{q}$ on $L$ is uniquely determined up to the addition of a closed $1$--form. The action of gauge transformations is the addition of exact forms, so the gauge equivalence class of $A$ is uniquely determined up to the addition of an imaginary multiple of $dt$, corresponding to tensoring $L$ by a flat line bundle.

\begin{remark}\label{rmk:Holonomy}
In order to calculate the holonomy of $A$ around a loop $\gamma _{z}=\{ z \} \times \Sph ^{1}_{t}$ one can use Lemma \ref{lem:Asymptotics:Periodic:Dirac:Higgs:Field}.(i) to show that $d\left( \int_{\gamma _{z}}{ A }\right) = -\int_{\gamma _{z}}{ F_{A} } = i \int_{\gamma _{z}}{ r (\partial _{r}G) dt } = i\, d\theta_{q}$. Here, as before, $\theta _{q}=\text{arg}(z-z_{0})$ if $q=(z_{0},t_{0}) \in X$.
\end{remark}

In a neighbourhood of the singularity $L$ is isomorphic to the Hopf line bundle extended radially from a small sphere $\Sph ^{2}$ enclosing the origin. At infinity $L$ is isomorphic to the radial extension of a line bundle of degree $1$ over the torus $\mathbb{T}^{2}_{\infty}$. Representatives for the connection in these asymptotic models are given by:
\begin{itemize}
\item[$\bullet$] Introduce spherical coordinates $(z,t)=(\rho\sin{\phi}\, e^{i\theta},\rho\cos{\phi})$ on a $3$-ball $B_{\sigma}$ around the singularity. The unique connection $A^{0}$ on $H$ with harmonic curvature $\frac{i}{2}\dvol_{\Sph ^{2}}$ is defined by $\frac{i}{2}(\pm 1 - \cos{\phi})d\theta$ in the standard cover $U_{\pm}=\Sph ^{2}\setminus (0,0,\pm 1)$ of $\Sph ^{2}$.

\item[$\bullet$] Consider the connection $A^{\infty}=-i\frac{t}{2\pi}d\theta$ on the trivial line bundle $\underline{\mathbb{C}}$ over $\Sph ^{1}_{\theta} \times \mathbb{R}_{t}$. If $(e^{i\theta },t,\xi ) \in \underline{\mathbb{C}}$, the map $\tau (e^{i\theta },t,\xi )=(e^{i\theta },t+2\pi,e^{i\theta}\xi )$ satisfies $\tau ^{\ast }A^{\infty}=A^{\infty}$. Define a line bundle with connection over $\mathbb{T}^{2}_{\theta ,t}$ as the quotient $(\underline{\mathbb{C}},A^{\infty })/\tau$.
\end{itemize}

Any connection $A$ on $L$ with $F_{A}=\ast d\Phi$ is asymptotically gauge equivalent to $A^{0}$ as $\rho \ra 0$. As $r\rightarrow\infty$, up to gauge transformations, $A$ is asymptotic to $A^{\infty}+i\alpha\, d\theta +ib\, dt$ for some $\alpha ,b\in\mathbb{R}/\mathbb{Z}$. The monodromy of this limiting connection is $e^{-i\theta-2\pi ib}$ around the circle $\{\theta\}\times \Sph ^{1}_{t}$ and $e^{it-2\pi i\alpha}$ around the circle $\Sph ^{1}_{\theta}\times\{ t\}$. While $b$ can be chosen arbitrarily, $\alpha$ is fixed by the Bogomolny equation \eqref{eqn:Bogomolny}. Indeed, \eqref{eqn:Periodic:Dirac:Higgs:Field} implies that $\partial _{t}G(z,t)=0$ if $t\in \pi\mathbb{Z}$ and therefore the connection $A$ restricted to the plane $\{ t=\pi\}$ is flat. On the other hand, as we approach infinity the limiting holonomy of $A$ on large circles $\{ r=\text{const},t=\pi\}$ converges to $e^{i(\pi -2\pi\alpha )}$. Thus $\alpha =\frac{1}{2}$ modulo $\Z$.

\begin{lemma}\label{lem:Asymptotics:Periodic:Dirac:Connection}
Fix parameters $(v,b) \in \R \times \R /\Z$. Let $(A,\Phi )$ be a solution to (\ref{eqn:Bogomolny}) with $\Phi =i \left( v+G \right)$ and such that the holonomy of $A$ around circles $\{ re^{i\theta} \} \times \Sph ^{1}_{t}$, $r \neq 0$, is $e^{-i\theta-2\pi ib}$.
\begin{itemize}
\item[(i)] In the region where $r\geq 2$ the connection $A$ is gauge equivalent to
\[
A^{\infty}+\frac{i}{2}\, d\theta +ib\, dt+a
\]
for a $1$--form $a$ such that $d^{\ast}a=0=\partial _{r}  \lrcorner\, a$ and $|a|+|\nabla a|=O(e^{-r})$.
\item[(ii)] In a ball of radius $\frac{\pi}{2}$ centred at the singular point $z=0=t$, $A$ is gauge equivalent to $A^{0}+a'$ where $|a'|+\rho |\nabla a'|=O(\rho ^{2})$ and $d^{\ast}a'=0=\partial _{\rho}\lrcorner \, a'$.
\end{itemize}
\proof
In order to prove (i), write $\Phi =i \left( v+\frac{1}{2\pi }\log{r} \right) +\psi$ and solve (\ref{eqn:Bogomolny}) in a radial gauge. Write $A=A^{\infty}+a$, where $a=a_{\theta }d\theta +a_{t}dt$ solves $da = \ast d\psi$:
\[
\left\{\begin{array}{l}
\partial _{r}a_{\theta}=r\partial _{t}\psi\\
\partial _{r}a_{t}=-\frac{1}{r}\partial _{\theta}\psi=0\\
\partial _{\theta }a_{t}-\partial _{t}a_{\theta}=r\partial _{r}\psi \\
\end{array}\right .
\]
Since $|\psi|+|\nabla \psi|=O(e^{-r})$, we can solve the system integrating along rays. Up to exponentially decaying terms, $a$ has a flat limit $a^{\infty}=a^{\infty}_{\theta }d\theta +a^{\infty}_{t}dt$ over the torus at infinity. By holonomy considerations as above, up to gauge transformations $a^{\infty}_{\theta}=\frac{i}{2}$ and $a^{\infty}_{t}=ib$. Then set $a_{\theta}-a^{\infty}_{\theta}=-\int _{r}^{\infty}{r\partial _{t}\psi}$ and $a_{t}=a^{\infty}_{t}$. Using these expressions one can check that $a$ is a solution to the system above because $\psi$ is harmonic; moreover, $d^{\ast}a=0$ because $\psi$ is independent of $\theta$. Finally, the decay of $\psi$ and its gradient imply the desired estimates. (ii) is proved similarly using Lemma \ref{lem:Asymptotics:Periodic:Dirac:Higgs:Field}.(iii).
\endproof
\end{lemma}

\subsection{The action of translations, rotations and scaling.}
Given an arbitrary point $q=(z_{0},t_{0})$ in $\RS$ the same formulas describe the asymptotic behaviour of the periodic Dirac monopole $(A_{q},\Phi _{q})$ with singularity at $q$ in coordinates centred at $q$. It will be useful to express the behaviour of $(A_{q},\Phi _{q})$ at large distances from $q$ in a fixed coordinate system.

\begin{lemma}\label{lem:Asymptotics:Periodic:Dirac:Translations}
For $r\geq 2|z_{0}|$ we have
\[
\begin{gathered}
\frac{1}{i}\Phi _{q}(z,t) = v+\frac{1}{2\pi}\log{r}-\frac{1}{2\pi}\, \text{Re}\left( \frac{z_{0}}{z}\right) +O(r^{-2}) \\
A_{q}(z,t) = A^{\infty } + ib\, dt + i\, \frac{t_{0}+\pi}{2\pi}\, d\theta -\frac{i}{2\pi}\, \text{Im}\left( \frac{z_{0}}{z}\right) dt+O(r^{-2}).
\end{gathered}
\]
\proof
Write $z=re^{i\theta}$ and $z_{0}=r_{0}e^{i\theta _{0}}$ and expand the logarithm for $r>r_{0}$
\[
\log{|z-z_{0}|}=\log{r}-\sum _{n=1}^{\infty} { \frac{(-1)^{n}}{n} \left( \frac{r_{0}}{r} \right) ^{n}\cos{\left[ n(\theta -\theta _{0}) \right]} } = \log{r}-\text{Re}\left( \frac{z_{0}}{z}\right) +O\left( \frac{r_{0}^{2}}{r^{2}}\right) .
\]
Together with Lemma \ref{lem:Asymptotics:Periodic:Dirac:Higgs:Field}.(ii), this proves the asymptotic expansion for the Higgs field. In order to derive an asymptotic expansion for the connection $A_{q}$, solve the abelian Bogomolny equation (\ref{eqn:Bogomolny}) using this asymptotic expansion for $\Phi$ as in the proof of Lemma \ref{lem:Asymptotics:Periodic:Dirac:Connection}.(i).
\qed
\end{lemma}

The choice of the parameters $(v,b) \in \R \times \R /\Z$ is related to rotations and dilations. By a rotation in the $z$--plane, we can always assume that $b=0$. On the other hand, given any $\lambda>0$ consider the homothety
\[
h_{\lambda }:\R ^{2} \times \R/2\pi\Z\longrightarrow \R ^{2} \times \R/2\pi\lambda\Z
\]
of ratio $\lambda$. We saw that the Bogomolny equation is the dimensional reduction of the ASD equation, which is conformally invariant. Then, forcing the Higgs field to scale as a $1$--form, $(h_{\lambda }^{\ast}A ,\lambda\, h_{\lambda }^{\ast }\Phi )$ is a monopole on $\R ^{2} \times \R/2\pi\Z$ if and only if $(A,\Phi )$ solves the Bogomolny equation on $\R ^{2} \times \R/2\pi\lambda\Z$. Now, given a periodic Dirac monopole $(A_{q},\Phi_{q})$ with mass $v$, set $\lambda =v+\frac{a_{0}}{2}$. Then as $v\rightarrow \infty$
\[
\lambda^{-1}h_{\lambda ^{-1}}^{\ast }\Phi\longrightarrow i\left( 1-\frac{1}{2\sqrt{r^{2}+t^{2}}}\right),
\]
\ie the limit $v\rightarrow \infty$ corresponds to the limit $\mathbb{R}^{2}\times \Sph ^{1}\rightarrow\mathbb{R}^{3}$ and in this limit a periodic Dirac monopole converges to an Euclidean Dirac monopole.

\section{Boundary conditions}\label{sec:Boundary:Conditions}

Having described the abelian periodic solutions to the Bogomolny equation, we proceed to state and discuss the boundary conditions for periodic monopoles (with singularities) introduced by Cherkis and Kapustin in \cite{Cherkis:Kapustin:1} and  \cite{Cherkis:Kapustin:2}. Periodic monopoles will be required to approach periodic Dirac monopoles of appropriate charges both at infinity and at the singularities. This is analogous to the case of $SU(2)$ monopoles on $\mathbb{R}^{3}$ without singularities, in which case it is well-known (\cf for example \cite[Chapter IV, Part II]{Jaffe:Taubes}) that every monopole with finite energy is asymptotic to an Euclidean Dirac monopole. Before giving precise definitions, we need to address the issue of which structure group to consider.

\subsection{The structure group: $SO(3)$ vs. $SU(2)$}\label{sec:Structure:Group}

Limiting ourselves to compact Lie groups of rank $2$, the simplest choice would be to take $G=SU(2)$. However, in order to introduce singularities of the fields while hoping to obtain smooth moduli spaces, it is necessary to pick $SO(3)$ as structure group. Indeed, Kronheimer \cite{Kronheimer:Thesis} showed that the moduli space of framed monopoles of charge $1$ on $\R ^{3}$ with one singularity at a point $p$ and structure group $G=SU(2)$ has a singularity of the form $\C ^{2}/\Z_{2}$. In \cite{Cherkis:Kapustin:2} Cherkis and Kapustin define periodic $U(2)$ and $SO(3)$--monopoles with singularities. We briefly discuss the relation between the two choices of structure group, following Braam--Donaldson \cite[\S 1.1-1.2, Part II]{Braam:Donaldson} and Donaldson \cite[\S 5.6]{Donaldson:Floer:Homology}.

Given a collection $S$ of $n$ distinct points $p_{1},\dots,p_{n} \in \RS$, let $V \ra (\RS) \setminus S$ be an $SO(3)$--bundle. By a result of Whitney \cite[\S III.7]{Whitney}, isomorphism classes of $SO(3)$--bundles over a CW--complex of dimension at most 3 are completely classified by the second Stiefel--Whitney class $w_{2}$. The second homology of $(\RS) \setminus S$ is generated by the classes of $2$--spheres $\mathbb{S}^{2}_{p_{i}}$ each enclosing the point $p_{i} \in S$. We fix the isomorphism class of $V$ by requiring that $w_{2}(V)\cdot [\mathbb{S}^{2}_{p_{i}}]=1$ for all $i=1,\ldots ,n$. $V$ does not lift to an $SU(2)$--bundle whenever $n>0$.

However, $V$ does always lift to a $U(2)$--bundle $E \ra (\RS) \setminus S$ with $c_{1}(E) \equiv w_{2}(V)$ (mod $2$). The adjoint bundle $\Lie{g}_{E}$ splits into a direct sum $\Lie{g}_{E} = \underline{\R} \oplus \Lie{g}^{(0)}_{E}$ of a trivial real line bundle, the trace part, and the trace-less part $\Lie{g}^{(0)}_{E} \simeq V$, a $PU(2) \simeq SO(3)$ bundle. A pair $(A,\Phi)$ on $E$ satisfying the Bogomolny equation induces an abelian monopole $(A_{tr},\Phi _{tr})$ on $\det(E)$ and an $SO(3)$--monopole $(A^{(0)},\Phi ^{(0)})$ on $V$. The moduli space of $U(2)$ monopoles on $E$ with fixed determinant is a double cover of the moduli space of $SO(3)$ monopoles on $V$, with $H^{1}((\RS) \setminus S;\Z _{2}) \simeq \Z_{2}$ as the group of deck transformations. Very concretely, the $\Z_{2}$--action is given by tensoring $E$ with the flat line bundle $L_{\frac{1}{2}}$ with holonomy $-\text{id}$ around circles $\gamma _{z} = \{ z \} \times \Sph ^{1}_{t}$.

We conclude that, up to a finite cover, it makes no difference to consider $U(2)$ monopoles with fixed central part and $SO(3)$ monopoles. Moreover, fixing boundary conditions resolves this ambiguity. We will work with structure group $G=SO(3)$, referring the reader to \cite{Cherkis:Kapustin:2} on how to adapt the definitions to the case $G=U(2)$.

\subsection{Boundary conditions for $SO(3)$--monopoles.}\label{sec:Boundary:Conditions:SO(3)}

We begin with some preliminary notational remarks. With the normalisation $|A|^{2}=-2\Trace{(A^{2})}$ of the norm on $\Lie{su}(2)$, the isomorphism $\Lie{so}(3) \simeq \Lie{su}(2)$ via the adjoint representation is an isometry. Observe that if $V \ra (\RS) \setminus S$ is a rank $3$ real oriented Riemannian vector bundle and $P$ is the principal $SO(3)$--bundle of orthonormal frames of $V$, then $V \simeq \ad\,P$. Finally, a reducible $SO(3)$--bundle $V$ is an oriented Riemannian rank $3$ vector bundle with a decomposition $V \simeq \underline{\R} \oplus M$ for an $SO(2)$--bundle $M$. We denote by $\hat{\sigma}$ the trivialising unit-norm section of the first factor. We will use the isomorphism $V \simeq \ad\,P$ to identify $\hat{\sigma}$ with $[\sigma _{3}, \,\cdot\, ]$, where $\sigma _{3} = \frac{1}{2}\text{diag}(i,-i)$, in a local trivialisation $\ad\, P \simeq U \times \Lie{su}_{2}$ over an open set $U$. In this sense we will talk of diagonal and off-diagonal sections of $V$ to denote the sections of the two factors in the decomposition $V \simeq \underline{\R} \oplus M$.

Fix a collection $S$ of $n$ distinct points $p_{1}, \ldots, p_{n} \in \RS$ and an $SO(3)$--bundle $V$ on $(\RS) \setminus S$ with the topology described above. We also fix an origin and a frame in $\RS$ and use coordinates $(z,t) \in \C \times \R /2\pi\Z$ with $z=x+iy=re^{i\theta}$. In \cite[\S 1.4]{Cherkis:Kapustin:1} and \cite[\S 2]{Cherkis:Kapustin:2} Cherkis and Kapustin consider the following boundary conditions for periodic monopoles (with singularities).

\begin{definition}\label{def:Boundary:Conditions:SO(3)}
Given a non-negative integer $k_{\infty} \in \Z _{\geq 0}$, parameters $(v,b) \in \R \times \R / \Z$ and a point $q=(\mu ,\alpha) \in \RS$, let $\mathcal{C}=\mathcal{C}(p_{1},\ldots ,p_{n},k_{\infty},v,b,q)$ be the space of smooth pairs $c=(A,\Phi )$ of a connection $A$ on $V$ and a section $\Phi$ of $V$ satisfying the following boundary conditions.

\begin{enumerate}
\item For each $p_{i} \in S$ there exists a ball $B_{\sigma}(p_{i})$ and a gauge $V|_{B_{\sigma}(p_{i}) \setminus \{ p_{i} \} } \simeq \underline{\R} \oplus H_{p_{i}}$ such that $(A,\Phi)$ can be written
\begin{alignat*}{2}
\Phi = -\frac{1}{2\rho_{i} }\, \hat{\sigma} +\psi\qquad & A = A^{0}\, \hat{\sigma} +a
\end{alignat*}
with $\xi= (a,\psi)=O(\rho_{i}^{-1+\tau})$ and  $|\nabla _{A}\xi| + |[\Phi,\xi]|=O(\rho_{i}^{-2+\tau})$ for some rate $\tau >0$. Here $\rho_{i}$ is the distance from $p_{i}$ and $A^{0}$ is the $SO(3)$--invariant connection on $H_{p_{i}}$.

\item There exists $R>0$ and a gauge $V \simeq \underline{\R} \oplus \left( L^{k_{\infty}}_{q} \otimes L_{v,b} \right)$ over $\left( \R^{2} \setminus B_{R} \right) \times \Sph ^{1}$ such that $(A,\Phi)$ can be written
\begin{gather*}
\Phi =\left[v + \frac{k_{\infty}}{2\pi }\log{r}-\frac{k_{\infty}}{2\pi}\text{Re}\left( \frac{\mu}{z}\right)\right] \,\hat{\sigma} +\psi\\
\\
A = \left[b\, dt +k_{\infty}A^{\infty }+\frac{k_{\infty}}{2\pi }(\alpha + \pi)d\theta -\frac{k_{\infty}}{2\pi}\text{Im}\left( \frac{\mu}{z}\right) dt\right] \, \hat{\sigma} + a
\end{gather*}
with $\xi= (a,\psi)=O(r^{-1-\tau})$ and  $|\nabla _{A}\xi| + |[\Phi,\xi]|=O(r^{-2-\tau})$ for some $\tau >0$. Here $A^{\infty}$ is the connection on $L_{q}$ of Lemma \ref{lem:Asymptotics:Periodic:Dirac:Connection}.
\end{enumerate}
\end{definition}

We refer to Lemmas \ref{lem:Decay:Singularity} and \ref{lem:Decay:Construction:Moduli:Space} for some discussion of the optimal rate of convergence of a monopole $(A,\Phi) \in \mathcal{C}$ to the asymptotic models. Here we collect some comments on Definition \ref{def:Boundary:Conditions:SO(3)}.

There is a topological constraint on the choice of the charge at infinity $k_{\infty}$. Indeed, since $[\T _{\infty}]$ is homologous to the sum $[\Sph ^{2}_{p_{1}}]+ \ldots +[\Sph ^{2}_{p_{n}}]$ and $k_{\infty}$ (mod $2$) is the value of the second Stiefel--Whitney class $w_{2}(V)$ on $[\T _{\infty}]$, we must have $k_{\infty} \equiv n$ modulo $2$. The (non-abelian) \emph{charge} of an $SO(3)$--monopole $(A,\Phi) \in \mathcal{C}$ is the non-negative integer $k$ defined by $2k=k_{\infty}+n$. In particular, for each charge $k$ the number of singularities cannot exceed $2k$. In the extremal case $k_{\infty}=0$ we require that $v>0$, so that $\Phi$ still defines a reduction $V \simeq \underline{\R} \oplus M$ of the structure group to $SO(2)$ both at infinity and close to the singularities.

The parameter $q$ in Definition \ref{def:Boundary:Conditions:SO(3)} is referred to as the \emph{centre} of the monopole. It is necessary to fix $q$ in order to have $L^{2}$--integrable infinitesimal deformations. Thus, differently from the Euclidean case, only moduli spaces of centred periodic monopoles carry a Riemannian metric induced by the $L^{2}$--norm of infinitesimal deformations. Notice that the boundary conditions of Definition \ref{def:Boundary:Conditions:SO(3)} depend on the choice of an origin and a frame in $\RS$.

Finally, Definition \ref{def:Boundary:Conditions:SO(3)} implies that non-trivial periodic monopoles have infinite energy
\begin{equation}\label{eqn:Yang:Mills:Higgs}
\mathcal{A}(A,\Phi)=\frac{1}{2}\int _{X}{|F_{A}|^{2}+|d_{A}\Phi |^{2}}.
\end{equation}

\section{Monopoles with a Dirac-type singularity}\label{sec:Singularity}

This and the next section, of a technical nature, are aimed to introduce the analytical tools needed to work with Definition \ref{def:Boundary:Conditions:SO(3)}. We begin in this section by studying monopoles on a punctured ball with a Dirac type singularity at the origin. We review the approach of Kronheimer \cite{Kronheimer:Thesis}, who showed that the Hopf fibration induces a bijection between monopoles on $\R ^{3}$ with Dirac type singularities and $\Sph ^{1}$--invariant instantons on $\R ^{4}$. This discussion will serve as motivation for the singular behaviour imposed in Definition \ref{def:Boundary:Conditions:SO(3)}. Moreover, in a number of points throughout the paper we will deduce decay properties of monopoles with Dirac type singularities from the $4$--dimensional theory. Next, we will introduce weighted Sobolev spaces and check that the necessary embedding and multiplication properties hold. Finally, we will study the mapping properties of the Laplacian $DD^{\ast}$, where $D$ is the Dirac operator of \eqref{eqn:Dirac:Operator}, in these weighted spaces. 

\subsection{Hopf lift of a monopole with a Dirac-type singularity}

Let $B^{3}=B_{\sigma}(0)$ be a ball in $\R ^{3}$. Fix complex coordinates $(z_{1},z_{2})$ on $\C^{2}\simeq \R ^{4}$ and consider the Hopf projection $\pi\co B^{4} \ra B^{3}$, $(z_{1},z_{2}) \mapsto (|z_{1}|^{2}-|z_{2}|^{2},2z_{1}z_{2}) \in \R \oplus \C$, which exhibits $B ^{4} \setminus \{ 0 \}$ as a circle bundle over $B^{3} \setminus \{ 0\}$ with fibre-wise circle action $e^{is}\cdot (z_{1},z_{2}) = (e^{is}z_{1},e^{-is}z_{2})$. Here $B^{4}=B_{\sqrt{\sigma}}(0) \subset \R ^{4}$. The Euclidean metric on $B^{4}\setminus \{ 0 \}$ can be expressed in Gibbons--Hawking coordinates \cite{Gibbons:Hawking} as
\begin{equation}\label{eqn:Gibbons:Hawking}
g_{\R ^{4}} = h\pi ^{\ast}g_{\R ^{3}} + h^{-1}\theta _{0}^{2},
\end{equation}
where $h$ is the harmonic function $h=\frac{1}{2\rho}$, $\rho$ is the distance from the origin in $\R ^{3}$ and $\theta _{0}$ is a connection of $\pi$ with $\ast dh = d\theta _{0}$.

Let $V \ra B^{3} \setminus \{ 0 \}$ be an $SO(3)$--bundle and $(A,\Phi)$ a connection and Higgs field on $V$. Define a connection $\hat{A}$ on $\pi ^{\ast}V \ra B^{4} \setminus \{ 0 \}$ by
\begin{equation}\label{eqn:Lift:Connection}
\hat{A}=\pi^{\ast}A-\pi ^{\ast}\left(h^{-1}\Phi \right) \otimes \theta _{0}.
\end{equation}
Then $\hat{A}$ is an $\Sph ^{1}$--invariant ASD connection on $B^{4} \setminus \{ 0 \}$. The following lemma is proved by Kronheimer as an application of Uhlenbeck's Removable Singularities Theorem \cite[Theorem 4.1]{Uhlenbeck:Removable:Singularities}.

\begin{lemma}[Lemma 3.5 of \cite{Kronheimer:Thesis}]\label{lem:Dirac:Monopole:Lift:ASD:Connection}
A smooth pair $(A,\Phi)$ is a monopole on $B^{3} \setminus \{ 0 \}$ such that
\begin{itemize}
\item[(i)] $h^{-1} |\Phi | \ra k \in \N$ as $\rho \ra 0$, and
\item[(ii)] $\int_{B^{3}}{ |d_{A}(h^{-1}\Phi)|^{2} h\dvol _{\R ^{3}} } < \infty$
\end{itemize}
if and only if $\hat{A}$ defined by (\ref{eqn:Lift:Connection}) is gauge equivalent to a smooth $\Sph ^{1}$--invariant ASD connection on $B^{4}$ and the $\Sph ^{1}$--action on the fibre over the origin of the extension of $\pi ^{\ast}V$ has weight $k$.
\end{lemma}

\begin{example}[Euclidean Dirac monopole]\label{rmk:Decay:Singularity}
Consider the model case of an Euclidean Dirac monopole $k(A_{0},\Phi _{0})$ of charge $k$ and vanishing mass on the reducible $SO(3)$--bundle $\underline{\R} \oplus H^{k}$ and the corresponding ASD connection $\hat{A}$. Writing $z_{i}=|z_{i}|e^{i\theta _{i}}$, a simple computation shows that the singular gauge transformation
\begin{equation}\label{eqn:Hopf:Lift:Gauge}
g = \left\{
\begin{array}{cc}
e^{k\theta_{1}\sigma_{3}}  & \text{if } z_{1}\neq 0 \\
e^{-k\theta_{2}\sigma_{3}}  & \text{if } z_{2}\neq 0 \\
\end{array}
\right . ,
\end{equation}
is an isomorphism $\pi ^{\ast}V \simeq (B^{4} \setminus \{ 0\} ) \times \Lie{su}(2)$ such that $g(\hat{A})=g\hat{A}g^{-1}-(dg)g^{-1}$ is the trivial connection. In this gauge the $\Sph^{1}$--action is given by
\begin{equation}\label{eqn:Hopf:Lift:Circle:Action}
e^{is} \cdot (z_{1},z_{2},X)=\left( e^{is}z_{1},e^{-is}z_{2},\Ad\left( e^{ks\sigma_{3}}\right) X \right)
\end{equation}
for $(z_{1},z_{2}) \in \C ^{2}$ and $X \in \Lie{su}(2)$.
\end{example}

In the general case of a monopole $(A,\Phi)$ with a Dirac type singularity of charge $k$ we deduce the decay of $(A,\Phi)$ to the model $k(A_{0},\Phi _{0})$ from Lemma \ref{lem:Dirac:Monopole:Lift:ASD:Connection}.

\begin{lemma}[\cf {\cite[Appendix A]{Cherkis:Kapustin:2}}]\label{lem:Decay:Singularity}
Let $(A,\Phi)$ be a monopole on $B^{3} \setminus \{ 0 \}$ such that
\begin{itemize}
\item[(i)] $\int_{B^{3}}{ |d_{A}(h^{-1}\Phi)|^{2}h } < \infty$;
\item[(ii)] $h^{-1}|\Phi| \ra k$ as $\rho \ra 0$.
\end{itemize}
Then there exists a gauge such that
\[
(A,\Phi) = k(A_{0},\Phi _{0}) + O(1).
\]
\proof
Let $\hat{A}$ be the corresponding smooth $\Sph ^{1}$--invariant connection on $B^{4}$. By parallel transport from $0 \in B^{4}$ we can define a trivialisation of $B^{4} \times \Lie{su}(2)$ such that
\begin{itemize}
\item[(a)] $|\hat{A}| \leq C|z|$, where $C$ depends on $\|F_{\hat{A}}\| _{L^{\infty}}$ and $|z|$ is the Euclidean distance from the origin in $\R ^{4}$;
\item[(b)] the $\Sph ^{1}$--action on $B^{4} \times \Lie{su}(2)$ takes the standard form \eqref{eqn:Hopf:Lift:Circle:Action}.
\end{itemize}
Consider the action of the singular gauge transformation \eqref{eqn:Hopf:Lift:Gauge}:
\[
g(\hat{A}) = g\hat{A}g^{-1} - (dg) g^{-1} = \pi ^{\ast}A - \pi ^{\ast}(h^{-1}\Phi) \otimes \theta _{0}
\]
and $- (dg) g^{-1} = k\pi ^{\ast}A_{0} - k\pi ^{\ast}(h^{-1}\Phi_{0}) \otimes \theta _{0}$ by Example \ref{rmk:Decay:Singularity}. Thus we have found a gauge such that $(A,\Phi) = k(A_{0},\Phi_{0}) + (a,\psi)$ with $\pi ^{\ast}a - \pi ^{\ast}(h^{-1}\psi) \otimes \theta _{0} = g\hat{A}g^{-1}$. Computing norms using the expression \eqref{eqn:Gibbons:Hawking} for the Euclidean metric, we find
\[
h^{-1}(|a|^{2}+|\psi|^{2}) = |g\hat{A}g^{-1}|^{2}=|\hat{A}|^{2} \leq C|z|^{2} = 2Ch^{-1}. \qedhere
\]
\end{lemma}

Finally, we observe that via the Hopf map the deformation theory of monopoles with a Dirac type singularity on $B^{3} \setminus \{ 0 \}$ corresponds to the one of $\Sph ^{1}$--invariant instantons on $B^{4}$. More precisely, the deformation theory of instantons is governed by the Dirac-type operator
\begin{equation}\label{eqn:Dirac:Operator:ASD}
\hat{D}:=2d_{\hat{A}}^{+} \oplus d_{\hat{A}}^{\ast}:\Omega ^{1}(B^{4};\pi ^{\ast}V) \ra \Omega ^{+}(B^{4};\pi ^{\ast}V) \oplus \Omega ^{0}(B^{4};\pi ^{\ast}V),
\end{equation}
where $\Omega ^{+}$ denotes the space of self-dual forms. Use the Hopf map to define lifts of forms as follows:
\begin{itemize}
\item[$(i)$] If $u \in \Omega ^{0}(B^{\ast};V)$ and $\alpha \in \Omega ^{1}(B^{\ast};V)$ set $\hat{u}=\pi^{\ast}u$ and $\hat{\alpha}=\pi^{\ast}(\ast h\alpha)+\pi^{\ast}\alpha\wedge \theta_{0}$. Observe that $|u|=|\hat{u}|$ and $|\hat{\alpha}|=|\alpha|$.
\item[$(ii)$] If $\xi =(a,\psi) \in \Omega (B^{\ast};V)$ define a $1$--form $\hat{\xi}$ with values in $\pi ^{\ast}V$ by:
\begin{equation}\label{eqn:Lift:Forms}
\hat{\xi}=\pi^{\ast}a-\pi^{\ast}(h^{-1}\psi)\otimes \theta_{0}
\end{equation}
We have already observed that $|\hat{\xi}|^{2}=h^{-1}\left( |a|^{2}+|\psi|^{2}\right)$.
\end{itemize}
Under these identifications the Dirac operator $\hat{D}$ and its adjoint $\hat{D}^{\ast}$ correspond to $h^{-1}D$ and $D^{\ast}$, respectively.

\subsection{Function spaces for gauge theory}

It is therefore possible to study the deformation theory of monopoles with a Dirac type singularity by studying the deformation theory of $\Sph ^{1}$--invariant instantons. This is the approach adopted by Pauly \cite{Pauly} to study singular monopoles on compact $3$--manifolds. On the other hand, it also makes sense to work directly in $3$--dimensions using weighted Sobolev spaces and a Dirac monopole as a background for the analysis. Some advantages of the latter approach are that one can work with stronger norms in terms of decay at the puncture and with $L^{2}$--spaces, because $W^{2,2} \hookrightarrow C^{0}$ in $3$ dimensions.

The theory of weighted Sobolev spaces is by now a fairly standard tool in many geometric problems. Classical references are Lockhart--McOwen's paper \cite{Lockhart:McOwen} and Melrose's book \cite{Melrose}. Our analysis is modelled on the work of Biquard \cite{Biquard:Paraboliques, Biquard:Prolongement} on singular connections on punctured Riemann surfaces and the work of Kronheimer--Mrowka \cite{Kronheimer:Mrowka} and R\aa de \cite{Rade:Local:Theory:1, Rade:Local:Theory:2, Rade:Global:Theory} on ASD connections with codimension $2$ singularities.

The exposition is standard except for a minor technical difficulty. The choice of weight function is dictated by two requirements: on one side, we want certain multiplicative properties to hold; on the other, we have to show that the Dirichlet problem for $DD^{\ast}$ can be solved for every appropriate boundary data. At first sight it seems that no choice of weighted spaces can satisfy both conditions. However, one can exploit the fact that we work on a reducible $SO(3)$--bundle $V = \underline{\R} \oplus H$ to resolve this issue. First, one defines weighted spaces so that the necessary multiplicative properties hold. The lack of surjectivity of the operator $DD^{\ast}$ acting between these spaces is easy to understand: it is necessary to enlarge the domain by adding constant diagonal sections. After this modification, it is crucial that the product on sections of $V$ is induced by the Lie bracket on $\Lie{su}_{2}$ to guarantee that the multiplicative properties are not destroyed.

\begin{definition}\label{def:Weighted:Spaces:Singularities:1}
Let $B^{\ast}$ be the punctured unit ball in $\RS$ and $V \ra B^{\ast}$ a Riemannian vector bundle endowed with a metric connection $A$. Given $\delta \in \R$ define the space $W^{m,p}_{\rho,\delta}$ as the closure of the space of sections $u \in C^{\infty}(B^{\ast};V)$ vanishing in a neighbourhood of the origin with respect to the norm:
\[
\| u \| _{W^{m,p}_{\rho,\delta}}^{p} = \sum_{j=0}^{m} \int{ \left| \rho ^{-\delta -\frac{3}{p} + j}\nabla _{A}^{j}u \right|^{p} \dvol _{\R ^{3}} }
\]
We will use the notation $L^{p}_{\rho,\delta}$ for $W^{0,p}_{\rho,\delta}$.
\end{definition}

\begin{remark}\label{rmk:Weighted:Sobolev:Spaces:Puncture:Cylinder}
\begin{itemize}
\item[(i)] $\rho ^{\beta} \in L^{p}_{\rho,\delta}$ if and only if $\beta > \delta$.
\item[(ii)] Pass to the conformal cylinder $(0,+\infty) \times \Sph ^{2}$ with metric
\[
g_{cyl} = d\tau ^{2} + g_{\Sph ^{2}} = \frac{d\rho^{2}}{\rho ^{2}} + g_{\Sph ^{2}},
\]
where we set $\tau = -\log{\rho}$. Then $u \in W^{m,p}_{\rho,\delta}$ if and only if $e^{\delta\tau}u \in W^{m,p}_{cyl}$, where the last symbol denotes the standard Sobolev space defined with respect to the cylindrical metric.
\end{itemize}
\end{remark}
The latter observation and the lemmas below are useful tools to work with these weighted spaces.

\begin{lemma}[{\cf \cite[Lemma 3.1]{Kronheimer:Mrowka}}]
If $u \in W^{m,p}_{loc}(B^{\ast})$ and $\| u \| _{W^{k,p}_{\rho,\delta}} < \infty$ then $u \in W^{k,p}_{\rho,\delta}$.
\end{lemma}

\begin{lemma}[{\cf \cite[Theorem 1.2]{Biquard:Paraboliques}}]\label{lem:Control:Radial:Weighted:Spaces}
For all $\delta \neq 0$ and $u \in C^{\infty}_{0}(B^{\ast})$ 
\[
\| u \| _{W^{1,p}_{\rho,\delta}} \leq \frac{1}{|\delta|} \| \nabla _{A}u \| _{L^{p}_{\rho,\delta -1}}.
\]
If $\delta >0$ it is not necessary to require that $u \equiv 0$ on $\partial B$.
\end{lemma}

We will now define spaces for gauge theory on the punctured ball modelled on the spaces $W^{m,2}_{\rho,\delta}$. Let $V$ be the reducible $SO(3)$--bundle $V=\underline{\R} \oplus H^{k} \ra B^{\ast}$ endowed with a pair $c=k(A_{0},\Phi _{0}) \, \hat{\sigma}$ induced by an Euclidean Dirac monopole of charge $k$, mass $0$ and singularity at the origin. For a $V$--valued form $u$ we will write $u=u_{D} + u_{T}$ in the decomposition into diagonal and off-diagonal part. We use covariant weighted $W^{m,2}_{\rho,\delta}$--norms for sections of $V$. Norms of $V$--valued differential forms are defined similarly by taking the $W^{m,2}_{\rho,\delta}$ norm of each component of the form.

\begin{definition}\label{def:Weighted:Spaces:Singularities}
Let $c=k(A_{0},\Phi_{0}) \, \hat{\sigma}$ be a Dirac monopole on $V=\underline{\R} \oplus H^{k} \ra B^{\ast}$ and fix $\delta >0$.
\begin{itemize}
\item[(i)] Define the gauge group $\mathcal{G}^{0}_{\delta}$ as the set of automorphisms $g$ of $V$ such that $(d_{1}g)g^{-1} \in L^{2}_{\rho,\delta -1}$ and $\nabla ^{2}_{A}g \in L^{2}_{\rho, \delta -2}$.
\item[(ii)] Define $\mathcal{C}^{0}_{\delta}$ as the space of configurations $c + (a,\psi)$ on $V$ with $(a,\psi) \in W^{1,2}_{\rho,\delta -1}$.
\item[(iii)] Define a space $\widetilde{W}^{2,2}_{\rho,\delta}$ of infinitesimal gauge transformations as
\[
\widetilde{W}^{2,2}_{\rho,\delta} = \left\{ (u_{D},u_{T}) \in L^{2}_{\rho,-\delta} \oplus L^{2}_{\rho,\delta} \, | \, \nabla _{A}u \in L^{2}_{\rho,\delta -1}, \nabla _{A}^{2}u \in L^{2}_{\rho,\delta -2} \right\}.
\]
\end{itemize}
\end{definition}

The fact that $\mathcal{G}^{0}_{\delta}$ is a group, at the moment unjustified, is Proposition \ref{prop:Weighted:Spaces:Singularities}.(a) below.

\begin{remark}
Since $\Phi$ acts by $-i\frac{k}{2\rho}$ on the off-diagonal component $u_{T}$ and trivially on the diagonal $u_{D}$, $\widetilde{W}^{2,2}_{\rho,\delta}$ can be defined using the equivalent norm
\[
\| u \| _{\widetilde{W}^{2,2}_{\rho,\delta}} \sim \| u \| _{L^{2}_{\rho,-\delta}} + \| \nabla _{A}u \| _{L^{2}_{\rho,\delta-1}} + \| [\Phi,u] \| _{L^{2}_{\rho,\delta-1}} + \|\nabla _{A}^{2}u \| _{L^{2}_{\rho,\delta-2}}.
\]
Similarly, the $\widetilde{W}^{2,2}_{\rho, \delta}$--norm of a $V$--valued form $u \in \Omega (B^{\ast};V)$ is defined by
\[
\| u \| ^{2}_{ \widetilde{W}^{2,2}_{\rho, \delta} } = \| u \| ^{2}_{L^{2}_{\rho, -\delta}}+ \| \nabla _{A}u \| ^{2}_{L^{2}_{\rho,\delta-1}} + |[\Phi,u] \| ^{2}_{L^{2}_{\rho,\delta-1}} + \| \nabla _{A}(D^{\ast}u)| ^{2}_{L^{2}_{\rho,\delta-2}} + |[\Phi,D^{\ast}u] \| ^{2}_{L^{2}_{\rho,\delta-2}}.
\]
If $u \in \Omega ^{0}(B^{\ast};V)$, $D^{\ast}(0,u)=-(d_{A}u,[\Phi,u])$ and this coincides with Definition \ref{def:Weighted:Spaces:Singularities}.(iii).
\end{remark}

The following lemma helps to understand the definition of the space $\widetilde{W}^{2,2}_{\rho,\delta}$.

\begin{lemma}\label{lem:Embeddings:Products:Singularity:1}
Fix $\delta >0$. There are continuous embeddings $\widetilde{W}^{2,2}_{\rho,\delta} \hookrightarrow C^{0}$ and $W^{2,2}_{\rho,\delta} \hookrightarrow \rho ^{\delta }C^{0}$. Moreover, $\| u - u(0) \| _{W^{2,2}_{\rho,\delta}} \leq C \| u \| _{\widetilde{W}^{2,2}_{\rho,\delta}}$ for all $u \in \widetilde{W}^{2,2}_{\rho,\delta}$.
\proof
The first claim is proved in three steps:
\begin{enumerate}
\item By the Sobolev embedding in $3$ dimensions and the assumption $\delta >0$, if $u\in \widetilde{W}^{2,2}_{\rho,\delta}$ then $\rho ^{-\delta +\frac{1}{2}}\nabla_{A} u \in L^{p}$ for all $2\leq p \leq 6$.
\item If $\delta \geq \frac{1}{2}$ conclude immediately that $\nabla _{A}u \in L^{p}$ for all $2 \leq p \leq 6$. Otherwise, by H\"older's inequality $\nabla _{A}u \in L^{p}$ for all $3<p<\frac{3}{1-\delta}$.
\item Kato's inequality and Morrey's estimate \cite[Theorem 4, \S 5.6.2]{Evans} imply that $u \in C^{0,\alpha}(B)$ for all $\alpha \in (0,\delta)$.
\end{enumerate}
The second statement follows from Remark \ref{rmk:Weighted:Sobolev:Spaces:Puncture:Cylinder}.(ii) and the Sobolev embedding with respect to the cylindrical metric, while the last claim is Lemma \ref{lem:Control:Radial:Weighted:Spaces}.
\qed
\end{lemma}
Thus we have an extension
\[
0 \ra W^{2,2}_{\rho, \delta} \ra \widetilde{W}^{2,2}_{\rho, \delta} \ra \R \hat{\sigma} \ra 0,
\]
where $\hat{\sigma}$ is a unit-norm section of the trivial factor in the decomposition $V = \underline{\R} \oplus H^{k}$.

\begin{remark}
By the definition of $d_{1}$, $g \in \mathcal{G}^{0}_{\delta}$ satisfies $\nabla_{A}g \in L^{2}_{\rho, \delta -1}, (g\Phi g^{-1} - \Phi) \in L^{2}_{\rho, \delta -1}$ and $\nabla ^{2}_{A}g \in L^{2}_{\rho, \delta -2}$. By Lemma \ref{lem:Embeddings:Products:Singularity:1} $g$ is continuous and has a well-defined limit over $0 \in B$; the condition $(g\Phi g^{-1} - \Phi) \in L^{2}_{\rho, \delta -1}$ forces this limiting value to lie in the stabiliser of $\Phi$.
\end{remark}

\begin{lemma}\label{lem:Embeddings:Products:Singularity}
Assume that all weighted spaces below are spaces of $V$--valued forms and the product on $V \simeq \ad (P_{V})$ is induced by the Lie bracket of $\Lie{su}_{2}$. If $\delta > 0$ the following are continuous maps:
\begin{enumerate}
\item $W^{1,2}_{\rho,\delta-1} \hookrightarrow L^{6}_{\rho,\delta-1}$
\item $\widetilde{W}^{2,2}_{\rho,\delta } \hookrightarrow C^{0}(B)$
\item $W^{1,2}_{\rho,\delta -1} \times W^{1,2}_{\rho,\delta -1} \ra L^{2}_{\rho,\delta -2}$
\item $\widetilde{W}^{2,2}_{\rho,\delta} \times L^{2}_{\rho,\delta -2} \ra L^{2}_{\rho,\delta -2}$
\item $\widetilde{W}^{2,2}_{\rho,\delta} \times \widetilde{W}^{2,2}_{\rho,\delta} \ra \widetilde{W}^{2,2}_{\rho,\delta}$
\item $\widetilde{W}^{2,2}_{\rho,\delta} \times W^{1,2}_{\rho,\delta -1} \ra W^{1,2}_{\rho,\delta -1}$
\end{enumerate}

In cases (3) and (6) the maps $W^{1,2}_{\rho,\delta -1} \ra L^{2}_{\rho,\delta -2}$ and $\widetilde{W}^{2,2}_{\rho,\delta} \ra W^{1,2}_{\rho,\delta -1}$ obtained by fixing the second factor are compact.
\proof
The embeddings (1) and (2) follow from the Sobolev embedding theorem with respect the cylindrical metric and Lemma \ref{lem:Embeddings:Products:Singularity:1}, respectively. The continuity of the products in (3)--(6) then follows easily using the embeddings (1)-(2), H\"older's inequality and the assumption $\delta >0$, as we now briefly explain. 

In order to prove (3) observe that by H\"older's inequality
\[
\| \xi \cdot \eta \| _{L^{2}_{\rho, \delta -2}} = \| \rho ^{-\delta + \frac{1}{2} }(\xi \cdot \eta) \| _{L^{2}} \leq \| \rho ^{-\delta +\frac{1}{2} }\xi \| _{L^{6}} \| \eta \| _{L^{3}} = \| \xi \| _{L^{6}_{\rho , \delta -1}} \| \eta \| _{L^{3}}
\]
and similarly
\[
\| \eta \| _{L^{3}} \leq \text{diam}(B)^{\delta} \| \eta \| ^{\frac{1}{2}}_{L^{6}_{\rho , \delta -1}} \| \eta \| ^{\frac{1}{2}}_{L^{2}_{\rho , \delta -1}}.
\]
The continuity of the product $W^{1,2}_{\rho,\delta -1} \times W^{1,2}_{\rho,\delta -1} \ra L^{2}_{\rho,\delta -2}$ now follows from (1). The compactness of the induced map $W^{1,2}_{\rho,\delta -1} \ra L^{2}_{\rho,\delta -2}$ is deduced by writing
\[
\| \xi \cdot (\eta_{i}-\eta _{i'}) \| _{L^{2}_{\rho, \delta -2}} \leq \| \xi \| _{L^{6}_{\rho , \delta -1}(B_{\sigma})} \| \eta_{i}-\eta _{i'} \| _{L^{3}(B_{\sigma })} + \| \xi \| _{L^{6}_{\rho , \delta -1}(B\setminus B_{\sigma})} \| \eta_{i}-\eta _{i'} \| _{L^{3}(B \setminus B_{\sigma })}
\] 
and using the fact that $\| \xi \| _{L^{6}_{\rho , \delta -1}(B_{\sigma})} \ra 0$ as $\sigma \ra 0$ together with the compactness of the embedding $W^{1,2} \hookrightarrow L^{3}$.

In view of the embedding in (2), the continuity of the map in (4) is immediate. For the statement in (5), observe that in the decomposition $u=u_{D} + u_{T}$ the product takes the form:
\[
(u_{D}+u_{T}) \cdot (v_{D} + v_{T}) = (u_{T} \cdot v_{T}) + ( u_{D} \cdot v_{T} + u_{T} \cdot v_{D})
\]
Therefore $\| (u \cdot v)_{D} \| _{L^{2}_{\rho, -\delta}} \leq \| u \| _{L^{\infty}} \| v \| _{L^{2}_{\rho, -\delta}}$ and
\[
\| (u \cdot v)_{T} \| _{L^{2}_{\rho, \delta}} \leq \sqrt{2} \| u_{D} \| _{L^{\infty}} \| v_{T} \| _{L^{2}_{\rho, \delta}} + \sqrt{2} \| v_{D} \| _{L^{\infty}} \| u_{T} \| _{L^{2}_{\rho, \delta}}.
\]
The rest of the proof of (5) and (6) follows easily making use of (3).
\endproof
\end{lemma}

\begin{prop}\label{prop:Weighted:Spaces:Singularities} For all $\delta > 0$
\begin{itemize}
\item[(a)] $\mathcal{G}^{0}_{\delta}$ is a Banach Lie group which acts smoothly on $\mathcal{C}^{0}_{\delta}$.
\item[(b)] The map $\Psi\co W^{1,2}_{\rho,\delta-1} \longrightarrow L^{2}_{\delta -2}$ defined by $\Psi (\xi) = \ast F_{A} - d_{A}\Phi + d_{2}\xi +\xi \cdot \xi$ is smooth.
\end{itemize}
\end{prop}

\subsection{Elliptic theory}

The proposition above shows that the spaces $\mathcal{C}^{0}_{\delta}$ and $\mathcal{G}^{0}_{\delta}$ are well-suited to study gauge theory. The next task is to find a range of values for $\delta >0$ such that the Laplacian $DD^{\ast}$ (coupled to Dirichlet boundary conditions) is an isomorphism $DD^{\ast}\co \widetilde{W}^{2,2}_{\rho, \delta} \ra L^{2}_{\rho, \delta -2}$.

We continue to work with the reducible pair $(A,\Phi)=k(A_{0},\Phi _{0})\, \hat{\sigma}$ given by an Euclidean Dirac monopole of charge $k$, zero mass and singularity at the origin. By changing variables to $\tau = -\log{\rho}$ the punctured ball $B^{\ast} = B_{\sigma} \setminus \{ 0 \}$ becomes the half cylinder $Q=(T,+\infty) \times \Sph ^{2}$, where $T=-\log{\sigma}$. The operator $\rho ^{2}DD^{\ast}$ has the form
\begin{equation}\label{eqn:Translational:Invariant:Operator}
\rho ^{2}DD^{\ast}u=-\ddot{u} + \dot{u} + Lu =: \mathcal{L}u
\end{equation}
where the dots denote derivatives with respect to $\tau$. $L$ is the positive self-adjoint operator on $\Sph ^{2}$ $L=\left( \triangle_{\Sph ^{2}}, \nabla _{A}^{\ast}\nabla _{A} + \frac{k^{2}}{4} \right)$ in the decomposition $V = \underline{\R} \oplus H^{k}$. Here $\nabla _{A}^{\ast}\nabla _{A}$ is the Laplacian of the connection $A=kA_{0}$ on $H^{k} \ra \Sph ^{2}$.

$\mathcal{L}$ is a translation-invariant operator on the cylinder $Q$. In view of Remark \ref{rmk:Weighted:Sobolev:Spaces:Puncture:Cylinder}.(ii), we want to study its mapping properties between weighted Sobolev spaces $\mathcal{L}\co e^{-\delta\tau}W^{2,2}_{cyl} \ra e^{-\delta\tau}L^{2}_{cyl}$. Lockhart--McOwen's theory \cite{Lockhart:McOwen} deals precisely with this kind of elliptic operators and their perturbations on cylinders and asymptotically cylindrical manifolds. Since we will study a boundary value problem, we introduce the appropriate spaces for the boundary data:
\begin{definition}\label{def:Weighted:Spaces:Singularity:Boundary}
Let $\partial\widetilde{W}^{2,2}_{\rho,\delta}$ be the closure of $C^{\infty}(\partial B; V|_{\partial B})$ with respect to the norm
\[
\| \varphi \| _{\partial\widetilde{W}^{2,2}_{\rho,\delta}} = \inf {\| \tilde{\varphi}\| _{\widetilde{W}^{2,2}_{\rho,\delta}}},
\]
where the infimum is taken over all $\tilde{\varphi} \in C^{\infty}(B^{\ast};V)$ such that $\tilde{\varphi}|_{\partial B}\equiv \varphi$.
\end{definition}

We associate to the operator $\mathcal{L}$ of \eqref{eqn:Translational:Invariant:Operator} a discrete set of weights, called \emph{exceptional}, as follows. Since $L$ is a self-adjoint positive operator its eigenvalues form a discrete sequence $0 \leq \lambda_{1} \leq \lambda _{2} \leq \ldots$ Moreover, we can select an orthonormal basis of $L^{2}(\Sph ^{2}; \underline{\R} \oplus H^{k})$ given by eigensections $\phi _{j}$ of $L$. Every solution to $\mathcal{L}u=0$ can be written
\[
u=\sum_{j=1}^{\infty}{ \left( A^{+}_{j}e^{-\gamma _{j}^{+}\tau } + A^{-}_{j}e^{ -\gamma _{j}^{-} \tau }\right) \phi_{j} }
\]
where $\gamma _{j}^{\pm}$ are the two solutions of $\gamma ^{2} + \gamma -\lambda _{j}$, \ie $\gamma_{j}^{\pm}= -\frac{1}{2} \pm \sqrt{\frac{1}{4}+\lambda _{j}}$. Define the set of exceptional weights of the operator $\mathcal{L}$ to be the collection $\mathcal{D}(\mathcal{L})$ of all $\gamma _{j}^{\pm}$, $j \geq 0$. The relevance of the exceptional weights is that the operator
\[
e^{-\delta \tau}W^{2,2}_{cyl} \longrightarrow e^{-\delta \tau}L^{2}_{cyl} \oplus \partial W^{2,2}_{cyl},
\]
defined by $u \longmapsto \mathcal{L}u \oplus u|_{\partial Q}$ is Fredholm for all $\delta \notin \mathcal{D}(\mathcal{L})$, \cf \cite[Theorem 6.3]{Lockhart:McOwen} and \cite[Theorems 5.60 and 6.5]{Melrose}. Here $\partial W^{2,2}_{cyl}$ is defined similarly to Definition \ref{def:Weighted:Spaces:Singularity:Boundary}.

\begin{lemma}\label{lem:Spectrum:Laplacian:Sphere}
The exceptional weights $\gamma ^{\pm}_{j} \in \mathcal{D}(\mathcal{L})$ are:
\[
\begin{dcases}
\gamma _{j}^{+}=j+\frac{|m|}{2}\\
\gamma _{j}^{-}=-j-1-\frac{|m|}{2}
\end{dcases}
\]
for $j=0,1,2,3,\ldots$ each with multiplicity $2j+|m|+1$. Here we take $m=0$ for the operator restricted to the diagonal component and $m=k$ when we restrict $\mathcal{L}$ to forms with values in $H^{k}$.
\proof
The eigenvalues of the Laplacian $\nabla _{A}^{\ast}\nabla _{A}$ of the $SO(3)$--invariant connection $mA_{0}$ on $H^{m}$ have been calculated by Kuwabara \cite[Theorem 5.1]{Kuwabara}. They are
\[
\frac{l(l+2)-m^{2}}{4},\qquad l=|m|+2j, \text{ for }j=0,1,2,\ldots 
\]
each with multiplicity $l+1$. Hence the eigenvalues of $L$ are $\frac{l(l+2)}{4}$, where we take $m=0$ on the diagonal component and $m=k$ on the off-diagonal part. The Lemma follows.
\qed
\end{lemma}
In particular, $0$ is an exceptional weight with multiplicity $1$ (the constant functions) for the operator $\mathcal{L}$ restricted to the diagonal part, while none of the weights in the interval $(-1-\frac{|k|}{2},\frac{|k|}{2})$ is exceptional for the operator restricted to the off-diagonal part.

\begin{prop}\label{prop:Dirichlet:Problem:Weighted:Spaces:Singularity}
Fix $0 < \delta <\min{\{ 1,\frac{|k|}{2}\}}$. The Dirichlet problem
\[
\begin{dcases}
\nabla _{A}^{\ast}\nabla _{A}u - \ad^{2}(\Phi)u = f\\
u|_{\partial B} = \varphi
\end{dcases}
\]
has a unique solution $u \in \widetilde{W}^{2,2}_{\rho,\delta}$ for all $f \in L^2_{\rho, \delta -2}$ and $\varphi \in \partial \widetilde{W}^{2,2}_{\rho,\delta}$. Moreover there exists a constant $C$ independent of $u,f,\varphi$ such that
\[
\| u \| _{\widetilde{W}^{2,2}_{\rho,\delta}} \leq C \left( \| f \| _{L^{2}_{\rho, \delta}} + \| \varphi \| _{\partial \widetilde{W}^{2,2}_{\rho,\delta}} \right).
\]
\end{prop}
The proposition is proved easily by separation of variables. Notice that introducing $\widetilde{W}^{2,2}_{\rho, \delta}$ , which is an extension of $W^{2,2}_{\rho, \delta}$ by constant diagonal sections, is necessary to be able to solve the Dirichlet problem for arbitrary boundary data.

\section{Analysis on the big end of $\R ^{2} \times \Sph ^{1}$}\label{sec:Big:End}

We move on to discuss the framework to tackle the analysis on the big end of $\RS$. The local model is provided in this case by a periodic Dirac monopole, or better its asymptotic form analysed in Lemmas \ref{lem:Asymptotics:Periodic:Dirac:Higgs:Field} and \ref{lem:Asymptotics:Periodic:Dirac:Connection}: we work on the $SO(3)$--bundle $V=\underline{\R} \oplus (L_{v,b} \otimes L^{k_{\infty}}_{q})$ endowed with the reducible pair $(A_{\infty},\Phi_{\infty})$ induced by a periodic Dirac monopole of centre $q$, charge $k_{\infty}$ and vacuum asymptotic parameters $v,b$. We will drop the subscript $\,_{\infty}$ for most of the section.

Fix $R>0$ so that for $r\geq R$ we can write $|\Phi| = v+\frac{ k_{\infty} }{2\pi}\log{r} + O(r^{-1})$. Hence we can find a constant $c=c(R,v,q)>0$ such that
\begin{alignat}{2}\label{eqn:Growth:Higgs:Field:Infinity}
|\Phi| \geq c \qquad & |d_{A}\Phi| \leq \frac{c}{r}
\end{alignat}
if $r \geq R$ (recall that we assume $v>0$ if $k_{\infty}=0$). Let $U_{R}$ be the open exterior domain $\R ^{2} \setminus \overline{B}_{R}$; we will drop the subscript $\,_{R}$ when it is not essential in the discussion. If $u$ is a section of $V$ we write $u=u_{D} + u_{T}$ in the decomposition into diagonal and off-diagonal part. Then in the region $U \times \Sph ^{1}$
\begin{equation}\label{eqn:Control:OffDiagonal:Infinity}
|[\Phi , u]|^{2} \geq c |u_{T}|^{2}.
\end{equation}
By Fourier analysis with respect to the circle variable $t$ we can further decompose $u_{D} = \Pi _{0}u_{D} + \Pi _{\perp}u_{D}$ into $\Sph ^{1}$--invariant and oscillatory part. On each circle $\{ z \} \times \Sph _{t}^{1}$ the following Poincar\'e inequality holds
\begin{equation}\label{eqn:Control:Oscillatory:Infinity}
\int_{\Sph ^{1}}{ |\nabla ( \Pi _{\perp}u_{D} )|^{2} } \geq \int_{\Sph ^{1}}{ |\Pi _{\perp}u_{D}|^{2} }.
\end{equation}

The inequalities \eqref{eqn:Control:OffDiagonal:Infinity} and \eqref{eqn:Control:Oscillatory:Infinity} suggest that, via the Weitzenb\"ock formula Lemma \ref{lem:Weitzenbock}, we have extremely good control of the off-diagonal and oscillatory piece of $u$ in terms of $DD^{\ast}u$. In order to control the $\Sph ^{1}$--invariant diagonal piece $\Pi _{0}u_{D}$ we introduce appropriate weighted spaces. An issue similar to the one encountered in Section \ref{sec:Singularity} arises here when trying to define weighted spaces for which good multiplication properties and the surjectivity of $DD^{\ast}$ hold at the same time.

\subsection{Function spaces for gauge theory}
Models for our analysis are the paper \cite{Biquard:Jardim}, where Biquard and Jardim study doubly periodic instantons with quadratic curvature decay, and analytic results of Amrouche, Girault and Giroire \cite{Amrouche:Girault:Giroire:1,Amrouche:Girault:Giroire}.

Fix $R>0$ and work on the exterior domain $U=U_{R} \subset \R ^{2}$. Define weight functions
\begin{equation}\label{eqn:Weight:Infinity}
\omega (z)=\sqrt{1+r^{2}}
\end{equation}
Notice that
\begin{alignat}{2}\label{eqn:Property:Weight:Exterior}
|\nabla\omega |\leq 1, \qquad & -\omega \triangle\omega +|\nabla\omega |^{2}=2
\end{alignat}
An important consequence of introducing the weight function $\omega$ is the following Poincar\'e inequality.

\begin{lemma}\label{lem:Poincare:Inequality:Laplacian:Exterior}
There exists a constant $C=C(c,R)$ such that
\[
\| \omega ^{-(\delta + 1)} \,u \| _{L^{2}} \leq \frac{C}{|\delta|}\left(  \| \omega ^{-\delta} \,\nabla u \| _{L^{2}} + \| \omega ^{-\delta} \,[\Phi,u] \| _{L^{2}} \right)
\]
for all $\delta \neq 0$ and all $u \in C ^{\infty}_{0}\left( \overline{U}\right)$ subject to the additional restriction $\Pi _{0}u_{D}|_{\partial U} =0$ when $\delta > 0$.
\proof
Decompose $u = \Pi _{0}u_{D}+\Pi_{\perp}u_{D}+u_{T}$. \eqref{eqn:Control:OffDiagonal:Infinity} and \eqref{eqn:Control:Oscillatory:Infinity} imply that if $\Pi _{0}u_{D}=0$ we have
\[
\| \omega ^{-(\delta + 1)} \,u \| _{L^{2}} \leq \frac{C}{\sqrt{1+R^{2}}} \left(  \| \omega ^{-\delta} \,\nabla u \| _{L^{2}} + \| \omega ^{-\delta} \,[\Phi,u] \| _{L^{2}} \right)
\]
Therefore suppose that $u=\Pi _{0}u_{D}$. The estimate is analogous to Lemma \ref{lem:Control:Radial:Weighted:Spaces} and is proved by integration by parts, \cf \cite[Theorem 1.2]{Biquard:Paraboliques}:
\begin{eqnarray*}
\int { \omega ^{-2(\delta +1)} u^{2} } & = & -\frac{1}{2\delta}\int{ d\left(\frac{1}{\omega ^{\delta}}
 \right) \wedge \frac{ u^{2} \ast dr}{r}} \leq \int {\frac{u(\partial _{r}u)}{r \omega ^{\delta} }}\\
 & \leq & C \left( \int { \omega ^{-2(\delta +1)}u^{2}}\right) ^{1/2}\left( \int { \omega ^{-2\delta}|\nabla u|^{2}}\right) ^{1/2}.
\end{eqnarray*}
The first inequality follows because under the hypothesis on $u$ the boundary term is always non-positive and the second one follows from H\"older's inequality with $C=\frac{\sqrt{2+R^{2}}}{R}$.
\qed
\end{lemma}

\begin{definition}\label{def:Weighted:Sobolev:Spaces:Infinity}
For a smooth $V$--valued form $u \in \Omega \left( U \times \Sph ^{1};V \right)$ and $\delta \in \R$ we define norms:
\begin{itemize}
\item[(i)] $\| u \| _{L^{2}_{\omega , \delta}}=\| \omega ^{-(\delta +1)}u \| _{L^{2}}$
\item[(ii)] $\| u \| ^{2}_{W^{1,2}_{\omega ,\delta}} = \int{ \omega ^{-2\delta -2}|u|^{2} } + \int{ \omega ^{-2\delta } \left( |\nabla _{A}u|^{2} + |[\Phi, u]|^{2} \right) }$
\item[(iii)] $\| u \| ^{2}_{W^{2,2}_{\omega , \delta}} = \| u \| ^{2}_{L^{2}_{\omega , \delta}} + \| \nabla _{A}u \| ^{2}_{L^{2}_{\omega, \delta -1}} + \| [\Phi, u] \| ^{2}_{L^{2}_{\omega, \delta -1}} + \| \nabla _{A}(D^{\ast}u) \|^{2}_{L^{2}_{\omega, \delta -2}} + \| [\Phi, D^{\ast}u] \| ^{2}_{L^{2}_{\omega, \delta -2}}$
\item[(iv)] $\| u \| ^{2}_{\widetilde{W}^{2,2}_{\omega , \delta}} = \| u \| ^{2}_{L^{2}_{\omega , -\delta}} + \| \nabla _{A}u \| ^{2}_{L^{2}_{\omega, \delta -1}} + \| [\Phi, u] \| ^{2}_{L^{2}_{\omega, \delta -1}} + \| \nabla _{A}(D^{\ast}u) \| ^{2}_{L^{2}_{\omega, \delta -2}} +  \| [\Phi, D^{\ast}u] \| ^{2}_{L^{2}_{\omega, \delta -2}}$
\end{itemize}
The corresponding weighted Sobolev spaces are defined as the closure of the space of smooth compactly supported forms with respect to these norms.
\end{definition}

\begin{remark}
\begin{itemize}
\item[(i)] Since $(A,\Phi)$ is a solution to the Bogomolny equation, the $W^{1,2}_{\omega , \delta}$--norm of a compactly supported form $u \in C^{\infty}_{0}(U\times \Sph ^{1})$ is equivalent to $\| u \| _{L^{2}_{\omega, \delta}} + \| D^{\ast}u \| _{L^{2}_{\omega, \delta -1}}$ by the Weitzenb\"ock formula Lemma \ref{lem:Weitzenbock} for $DD^{\ast}$.
\item[(ii)] In view of \eqref{eqn:Control:OffDiagonal:Infinity} and \eqref{eqn:Control:Oscillatory:Infinity}, if $u \in W^{1,2}_{\omega,\delta}$ then $\Pi _{\perp}u_{D}, u_{T} \in L^{2}_{\omega, \delta -1}$. In particular, the only difference between the spaces $W^{2,2}_{\omega, \delta}$ and $\widetilde{W}^{2,2}_{\omega, \delta}$ consists in the chosen weighted $L^{2}$--norm of $\Pi _{0}u_{D}$. It follows from the proof of Lemma \ref{lem:Multiplication:Embedding:Infinity} below that when $\delta <0$ we have an extension
\[
0\ra W^{2,2}_{\omega, \delta} \ra \widetilde{W}^{2,2}_{\omega, \delta} \ra \R \, \hat{\sigma} \ra 0.
\]
\end{itemize}
\end{remark}

\begin{definition}\label{def:Weighted:Sobolev:Spaces:Gauge:Theory:Infinity}
Fix $\delta <0$.
\begin{itemize}
\item[(i)] $\mathcal{G}^{\infty}_{\delta}$ is the space of sections $g$ of $\text{Aut}(V)$ over $U \times \Sph ^{1}$ such that $(d_{1}g)g^{-1} \in W^{1,2}_{\omega,\delta-1}$.
\item[(ii)] $\mathcal{C}^{\infty}_{\delta}$ is the space of pairs $(A,\Phi)$ on $V$ of the form $(A_{\infty},\Phi _{\infty}) + (a,\psi)$, where $\xi = (a,\psi)$ is a section of $(\Lambda ^{1} \oplus \Lambda ^{0}) \otimes V$ of class $W^{1,2}_{\omega,\delta-1}$.
\item[(iii)] Infinitesimal gauge transformations are elements of $\widetilde{W}^{2,2}_{\omega, \delta -2}(U \times \Sph ^{1};V)$.
\end{itemize}
\end{definition}

\begin{lemma}\label{lem:Multiplication:Embedding:Infinity}
Fix $\delta \in (-1,0)$.
\begin{itemize}
\item[(i)] If $\xi = \Pi _{0}\xi_{D} + \Pi _{\perp}\xi _{D} + \xi _{T} \in W^{1,2}_{\omega, \delta -1}$ is a $V$--valued differential form then
\[
\omega ^{-\delta}\Pi _{0}\xi_{D}, \omega ^{-\delta +1}\Pi _{\perp}\xi _{D}, \omega ^{-\delta +1}\xi _{T} \in L^{p}
\]
for all $2 \leq p \leq 6$ and the inclusions are continuous.
\item[(ii)] $\widetilde{W}^{2,2}_{\omega, \delta} \hookrightarrow C^{0}$ is a continuous embedding.
\end{itemize}
The following products are continuous:
\begin{itemize}
\item[(iii)] $\widetilde{W}^{2,2}_{\omega,\delta} \times \widetilde{W}^{2,2}_{\omega,\delta} \ra \widetilde{W}^{2,2}_{\omega,\delta}$
\item[(iv)] $\widetilde{W}^{2,2}_{\omega,\delta} \times W^{m,2}_{\omega,\delta -2 +m} \ra W^{m,2}_{\omega,\delta -2 +m}$ for $m=0,1$
\item[(v)] $W^{1,2}_{\omega,\delta -1} \times W^{1,2}_{\omega,\delta -1} \ra L^{2}_{\omega,\delta -2}$
\end{itemize}
Moreover, the maps $\widetilde{W}^{2,2}_{\delta} \ra W^{m,2}_{\delta -2+m}$ and $W^{1,2}_{\omega,\delta -1} \ra L^{2}_{\omega,\delta -2}$ induced by (iv) and (v) by fixing the second argument are compact. Here the products are those induced by the Lie bracket on $\Lie{su}(2)$ under the identification $V \simeq \ad\, P$.
\proof
\begin{itemize}
\item[(i)] It is a consequence of the Sobolev embedding theorem $W^{1,2} \hookrightarrow L^{6}$ in $3$ dimensions and the fact that if $\xi \in W^{1,2}_{\omega,\delta -1}$ then $\omega ^{-\delta +1}\Pi _{\perp}\xi _{D}, \omega ^{-\delta +1}\xi _{T} \in L^{2}$.
\item[(ii)] For the oscillatory and off-diagonal part this is a consequence of the standard Sobolev embedding $W^{2,2} \hookrightarrow C^{0}$. In fact we have more: if $\Pi _{0}u_{D}=0$ then $\omega ^{-(\delta -1)}u \in W^{2,2}$ and therefore $u \in \omega ^{\delta -1}C^{0}$.

Suppose instead that $u = \Pi _{0}u_{D}$, so that we can work on $U \subset \R ^{2}$. First of all we can replace $\omega$ with $r$ because the two weights are equivalent (with a constant depending on $R$) on $U$. If $\nabla u \in W^{1,2}_{\omega, \delta -1}$, $r^{-\delta +1}\nabla u \in W^{1,2}_{cyl}$, where the latter is the standard Sobolev space with respect to the cylindrical metric $r^{-2}g_{\R ^{2}}$. Thus $r^{-\delta + 1}\nabla u \in L^{p}_{cyl}$ for all $p \in [2,\infty)$ by the standard Sobolev embedding. By an inversion $r=\frac{1}{\rho}$ we consider the function $\tilde{u}(\rho e^{i\theta})= u(\rho^{-1} e^{i\theta})$ defined on a punctured ball $B_{1/R} \subset \R ^{2}$. It is integrable because $u \in L^{2}_{\omega, -\delta}$ and $\delta >-1$ ($\delta >-2$ would be enough). Moreover, $\tilde{u}$ has gradient in $L^{p}$ for all $p < \frac{2}{1+\delta}$. Since $\delta <0$ we can choose $p>2$ and apply Morrey's estimate \cite[Theorem 4, \S 5.6.2]{Evans} to show that $\tilde{u}$, and therefore $u$, is continuous. In particular there exists a well-defined limit of $u_{\infty}=\lim_{r\ra\infty}{ u(re^{i\theta}) }$ and, by Lemma \ref{lem:Poincare:Inequality:Laplacian:Exterior}, $u-u_{\infty} \in W^{2,2}_{\omega, \delta}$.
\end{itemize}
The rest of the Lemma now follows easily in a way similar to Lemma \ref{lem:Embeddings:Products:Singularity}. It is crucial to observe that terms of the form $u_{D} \cdot v_{D}$ do not appear in the products.
\qed
\end{lemma}

\begin{prop}\label{prop:Weighted:Spaces:Infinity}
For all $\delta \in (-1,0)$, $\mathcal{G}^{\infty}_{\delta}$ is a Banach Lie group acting smoothly on $\mathcal{C}^{\infty}_{\delta}$.

Moreover, the map $\Psi\co \mathcal{C}^{\infty}_{\delta} \ra L^{2}_{\omega ,\delta -2}(U \times \Sph ^{1};\Lambda ^{1} \otimes V); (A,\Phi) \mapsto \ast F_{A}-d_{A}\Phi$ is smooth. 
\end{prop}

\subsection{Elliptic theory}

We now study the equation $DD^{\ast}u=\nabla _{A}^{\ast}\nabla _{A}u-\ad {(\Phi )}^{2}u =f$ for $f \in L^{2}_{\omega, \delta -2}$ and $u \in \widetilde{W}^{2,2}_{\omega,\delta}$ with $\delta <0$ sufficiently close to $0$.

\begin{prop}\label{prop:Laplacian:Infinity}
There exists $-1 \leq \delta_{0} <0$ and $R_{0}>0$ such that if either
\begin{itemize}
\item[(i)] $\delta \in (\delta_{0},0)$ and $R>0$ is arbitrary, or
\item[(ii)] $\delta \in (-1,0)$ is arbitrary and $R\geq R_{0}$,
\end{itemize}
then the following holds. For all $f\in L^{2}_{\omega,\delta -2}$ and $\varphi \in \partial \widetilde{W}^{2,2}_{\omega, \delta}$ there exists a unique solution $u \in \widetilde{W}^{2,2}_{\omega,\delta}$ to the Dirichlet problem
\[
\begin{dcases*}
DD^{\ast}u=f & in $U_{R} \times \Sph ^{1}$\\
u= \varphi & on $\partial U_{R} \times \Sph ^{1}$
\end{dcases*}
\]
Moreover there exists a constant $C=C(\delta)>0$ independent of $u$ and $f$ such that
\[
\| u \| _{\widetilde{W}^{2,2}_{\omega, \delta}} \leq C \left( \| f\| _{L^{2}_{\omega,\delta -2}} + \| \varphi \| _{\partial \widetilde{W}^{2,2}_{\omega, \delta}} \right).
\]
\proof
First suppose that $f=\Pi _{0}f_{D}$ so that we work on the exterior domain $U_{R} \subset \R ^{2}$. In this case, one can take $\delta \in (-1,0)$ and $R>0$ arbitrary. The proof is by separation of variables as for Proposition \ref{prop:Dirichlet:Problem:Weighted:Spaces:Singularity}. It is necessary to consider the extension $\widetilde{W}^{2,2}_{\omega, \delta}$ of $W^{2,2}_{\omega, \delta}$ by the constant functions to be able to solve the Dirichlet problem for arbitrary boundary data.

Assume instead that $\Pi_{0}f_{D}=0$. Then \eqref{eqn:Control:OffDiagonal:Infinity} and \eqref{eqn:Control:Oscillatory:Infinity} imply
\begin{equation}\label{eqn:Poincare:Proof:Dirichlet:Problem:Infinity}
c \int_{\Sph ^{1}}{ |u|^{2} } \leq \int_{\Sph ^{1}}{ |\nabla _{A}u|^{2} + |[\Phi,u]|^{2} }.
\end{equation} 
Since $L^{2}_{\omega,\delta -2} \subset L^{2}$ we obtain a solution $u \in L^{2} \subset L^{2}_{\omega , -\delta}$ to the Dirichlet problem by direct minimisation of the functional $\frac{1}{2}\int {|\nabla _{A}u|^{2}+|[\Phi,u]|^{2} }-\int{ \langle u,f \rangle }$. We have to show that $u \in \widetilde{W}^{2,2}_{\omega,\delta}$.

\begin{itemize}
\item[Step 1.] We can always reduce to the case $\varphi = 0$ by extending $\varphi$ to $\widetilde{\varphi} \in \widetilde{W}^{2,2}_{\omega, \delta}$ such that $\| \widetilde{\varphi} \| _{\widetilde{W}^{2,2}_{\omega, \delta}} \leq \| \varphi \| _{\partial \widetilde{W}^{2,2}_{\omega, \delta}} $ and replacing $u$ with $u-\widetilde{\varphi}$ and $f$ with $f-DD^{\ast}\widetilde{\varphi}$.

\item[Step 2.] Since $u$ vanishes on the boundary, an integration by parts yields (all integrals are taken over $U_{R} \times \Sph ^{1}$):
\begin{align*}
\int{ \omega^{-2\delta} \langle \nabla _{A}^{\ast}\nabla _{A}u - \ad^{2}(\Phi)u,u \rangle } = &\int{ \omega^{-2\delta}\left( |\nabla _{A}u|^{2} + |[\Phi,u]|^{2} \right) } \\
&-2\delta \int{ \omega^{-2\delta -1}\langle \nabla_{A}u,u\otimes d\omega\rangle }
\end{align*}
To control the last term, use H\"older's inequality, \eqref{eqn:Property:Weight:Exterior} and \eqref{eqn:Poincare:Proof:Dirichlet:Problem:Infinity}:
\[
\left| 2\delta \int{ \omega^{-2\delta -1} \langle \nabla_{A}u,u\otimes d\omega\rangle } \right| \leq C|\delta| \| \omega ^{-1} \| _{L^{\infty}} \int{ \omega^{-2\delta}\left( |\nabla _{A}u|^{2} + |[\Phi,u]|^{2} \right) }
\]
Thus if $|\delta|$ is sufficiently small or if $R$ is sufficiently large we deduce
\[
\| \nabla_{A}u \| _{L^{2}_{\omega, \delta -1}} + \| [\Phi,u]  \| _{L^{2}_{\omega, \delta -1}} \leq C \| DD^{\ast}u \| _{L^{2}_{\omega , \delta -2}}.
\]

In other words, in view of \eqref{eqn:Poincare:Proof:Dirichlet:Problem:Infinity} and the definition of $D^{\ast}$, we proved
\[
\| u \| _{L^{2}} + \| \omega ^{-\delta} D^{\ast}u \| _{L^{2}} \leq C\| DD^{\ast}u \| _{L^{2}_{\omega , \delta -2}}.
\] 

\item[Step 3.] Notice that if $\chi$ is a smooth function supported in a compact set $K \subset U_{R} \times \Sph ^{1}$, then
\[
\| DD^{\ast}(\chi u) \| _{L^{2}_{\omega, \delta -2}} \leq C \left( \| \nabla ^{2}\chi \| _{L^{2}} \| u \| _{L^{2}(K)} + \| \nabla \chi \|_{L^{2}} \| \nabla _{A}u \| _{L^{2}(K)} + \| DD^{\ast}u \| _{L^{2}_{\omega, \delta -2}} \right)
\]
and similarly $\| \omega ^{-\delta }D^{\ast}(\chi u) \| _{L^{2}} \leq C \left( \| \nabla \chi \| _{L^{2}}\| u \| _{L^{2}(K)} + \| \omega ^{-\delta} D^{\ast}u \| _{L^{2}} \right)$.

Choose $\chi \in C^{\infty}$ with $\chi \equiv 1$ on $\{ r \leq R+1 \}$ and $\chi \equiv 0$ if $r \geq R+2$. Write $u = \chi u + (1-\chi)u$. By Step 2 and standard elliptic regularity close to the boundary (\cf for example \cite[Theorem 8.12]{Gilbarg:Trudinger}),
\[
\| \chi u \| _{\widetilde{W}^{2,2}_{\omega, \delta}} \leq C \| \chi u \|_{W^{2,2}} \leq C \left( \| DD^{\ast}(\chi u) \| _{L^{2}} + \| \chi u \| _{L^{2}} \right) \leq C \| f \| _{L^{2}_{\omega, \delta -2}}.
\]

Hence we reduced to prove
\begin{equation}\label{eqn:Bochner:Integration:By:Parts}
\| \nabla _{A}\xi \| _{L^{2}_{\omega, \delta -2}} + \| [\Phi, \xi] \| _{L^{2}_{\omega, \delta -2}} \leq C \left( \| D\xi \| _{L^{2}_{\omega, \delta -2}} + \| \xi \| _{L^{2}_{\omega, \delta -1}} \right)
\end{equation} 
for $\xi = D^{\ast}\big( (1-\chi)u \big)$, \ie with $\xi$ vanishing in a neighbourhood of $\partial U _{R} \times \Sph ^{1}$.

\item[Step 4.] The Weitzenb\"ock formula for $D^{\ast}D$ in Lemma \ref{lem:Weitzenbock} implies
\begin{equation}\label{eqn:Bochner}
\frac{1}{2}d^{\ast}d\left( |\xi |^{2}\right) = -|\nabla _{A}\xi |^{2}-|\Phi\xi |^{2}+\langle D^{\ast}D\xi ,\xi\rangle -2\langle d_{A}\Phi\cdot\xi ,\xi\rangle .
\end{equation}
Integrate this Bochner-type identity against $\omega^{-2\delta +2}$ and integrate by parts:
\begin{align}\label{eqn:Bochner:Integration:By:Parts:1}
\int{ \omega^{-2\delta +2} \left( |\nabla _{A}\xi |^{2}+| [\Phi,\xi ] |^{2} \right) } &\leq  
\int{ \omega^{-2\delta +2}| D\xi |^{2} } -2 \int{ \omega^{-2\delta +2}\langle d_{A}\Phi \cdot \xi , \xi \rangle }\\ \nonumber
&{} + 2(1-\delta) \int{ \omega^{-2\delta +1} \langle D\xi , d\omega \cdot \xi \rangle }\\ \nonumber
&{} +  2(1-\delta) \int{ \omega^{-2\delta +1} \langle \nabla _{A}\xi , d\omega \otimes \xi \rangle }\\ \nonumber
&{} + 2(1-\delta)^{2} \int{ \omega ^{-2\delta} |\xi| ^{2} }
\end{align}

Consider the term $\int{ \omega^{-2\delta +2}\langle d_{A}\Phi \cdot \xi , \xi \rangle }$. Since $(A,\Phi)$ is reducible this term only involve $\xi _{T}$. Moreover, by \eqref{eqn:Growth:Higgs:Field:Infinity} $\omega |d_{A}\Phi| \leq c$. Then H\"older's and Young's inequality with $\varepsilon > 0$ imply
\[
\left| \int{ \omega^{-2\delta +2}\langle d_{A}\Phi \cdot \xi_{T} , \xi_{T} \rangle } \right| \leq \frac{c}{\varepsilon _{1}} \int{ \omega^{-2\delta} |\xi_{T}|^{2} } + c\varepsilon _{1} \int{ \omega^{-2\delta +2} |\xi_{T}|^{2} }
\]
for any $\varepsilon _{1} >0$. Moreover, by \eqref{eqn:Growth:Higgs:Field:Infinity}
\[
c\varepsilon _{1}\int{ \omega^{-2\delta +2} |\xi_{T}|^{2} } \leq \varepsilon _{1}\int{ \omega^{-2\delta +2} |[\Phi, \xi]|^{2} }.
\]

Secondly, by H\"older's inequality
\[
\left| \int{\omega ^{-2\delta +1} \langle D\xi, d\omega \cdot \xi \rangle  } \right| \leq \| \omega ^{-\delta +1}D\xi \| _{L^{2}} \|\omega ^{-\delta}\xi \| _{L^{2}} \leq \frac{1}{2}\| \omega ^{-\delta +1}D\xi \| ^{2}_{L^{2}} + \frac{1}{2}\|\omega ^{-\delta}\xi \| _{L^{2}}^{2}
\]
because $|d\omega| \leq 1$ by \eqref{eqn:Property:Weight:Exterior}. Similarly, for any $\varepsilon _{2}>0$:
\begin{align*}
\left| \int{\omega ^{-2\delta +1} \langle \nabla _{A}\xi, d\omega \cdot \xi \rangle  } \right| & \leq \| \omega ^{-\delta +1}\nabla _{A}\xi \| _{L^{2}} \|\omega ^{-\delta}\xi \| _{L^{2}}\\
& {} \leq \varepsilon _{2} \| \omega ^{-\delta +1}\nabla _{A}\xi \| ^{2}_{L^{2}} + \frac{1}{\varepsilon _{2}} \|\omega ^{-\delta}\xi \| _{L^{2}}^{2}
\end{align*}

Now choose $\varepsilon_{1}, \varepsilon _{2} <1$ so that the appropriate terms can be absorbed in the LHS of \eqref{eqn:Bochner:Integration:By:Parts:1} to obtain \eqref{eqn:Bochner:Integration:By:Parts}. \qedhere\end{itemize}
\end{prop}

\begin{remark}\label{rmk:Bochner:Integration:By:Parts}
For later use, notice that the a priori estimate \eqref{eqn:Bochner:Integration:By:Parts} holds for any $\delta \in \R$.
\end{remark}

\section{Construction of the moduli spaces}\label{sec:Construction:Moduli:Space}

In this section the local analysis of Sections \ref{sec:Singularity} and \ref{sec:Big:End} is used to prove that moduli spaces of $SO(3)$ periodic monopoles (with singularities) are, when non-empty, smooth hyperk\"ahler manifolds for generic choices of the parameters defining the boundary conditions of Definition \ref{def:Boundary:Conditions:SO(3)}. Before proceeding with the proof, we make precise definitions of the spaces of connections, Higgs fields and gauge transformations combining Definitions \ref{def:Weighted:Spaces:Singularities} and \ref{def:Weighted:Sobolev:Spaces:Gauge:Theory:Infinity}.

Fix a collection $S$ of $n$ distinct points $p_{1}, \ldots ,p_{n} \in \R ^{2} \times \Sph ^{1}$. Let $V \ra (\RS) \setminus S$ be an $SO(3)$--bundle such that $w_{2}(V) \cdot [\Sph ^{2}_{p_{i}}] = 1$. Denote by $P$ the associated principal $SO(3)$--bundle. Choose parameters $k_{\infty} \in \Z _{\geq 0}$ with $k_{\infty} \equiv n$ (mod $2$) and $v,b \in \R \times \R /\Z$, $q \in \RS$, with $v>0$ if $k_{\infty}=0$. Let $\mathcal{C}=\mathcal{C}(p_{1},\ldots ,p_{n},k_{\infty},v,b,q)$ be the space of smooth pairs of a connection and a Higgs field on $V$ as in Definition \ref{def:Boundary:Conditions:SO(3)}.

Fix a smooth pair $c=(A,\Phi) \in \mathcal{C}$, which we will refer to as the background pair. One such pair will be constructed in Section \ref{sec:Dimension}. We can always assume that there exist preferred gauges over $B_{\sigma}(p_{i}) \setminus \{ p_{i} \}$ and $U_{R} \times \Sph ^{1}$, for small $\sigma >0$ and large $R>0$, such that $c$ coincides with the asymptotic models over these regions. Given $c$, we use it as a background to define spaces $W^{m,2}_{\rho, \delta _{1}}$ and $W^{m,2}_{\omega , \delta _{2}}$ of forms with values in $V|_{B_{\sigma}(p_{i}) \setminus \{ p_{i} \} }$ and $V|_{U_{R} \times \Sph ^{1} }$ as in Definitions \ref{def:Weighted:Spaces:Singularities:1} and \ref{def:Weighted:Sobolev:Spaces:Infinity}.

\begin{definition}\label{def:Weighted:Spaces:Global}
Given $\sigma, R>0$, set $K_{\sigma ,R}= \left( \overline{B}_{R} \times \Sph ^{1} \right)\setminus \bigcup_{i=1}^{n}{ B_{\sigma}(p_{i}) }$.
\begin{itemize}
\item[(i)] A $V$--valued form $u \in L^{2}_{loc}$ on $(\RS) \setminus S$ belongs to the global weighted Sobolev space $L^{2}_{(\delta_{1}, \delta _{2})}$ if, in the preferred gauges around each singularity and at infinity, $u|_{B_{\sigma}(p_{i})\setminus\{ p_{i}\} } \in L^{2}_{\rho, \delta_{1} }$ and $u|_{U_{R} \times \Sph ^{1} } \in L^{2}_{\omega, \delta _{2} }$. We define a norm on $L^{2}_{ (\delta_{1}, \delta _{2}) }$ by taking the maximum of the semi-norms $\| u|_{ B_{\sigma}(p_{i}) \setminus \{ p_{i}\}} \| _{ L^{2}_{\rho, \delta_{1} } }$, $\| u |_{U_{R}\times\Sph ^{1}} \| _{ L^{2}_{\omega,\delta_{2} } }$ and $\| u |_{K_{\frac{\sigma}{2},2R}} \| _{L^{2}}$.

The spaces $\widetilde{W}^{2,2}_{ (\delta_{1}, \delta _{2}) }, W^{2,2}_{ (\delta_{1}, \delta _{2}) }$ and $W^{1,2}_{ (\delta_{1}, \delta _{2}) }$ are defined in a similar way.

\item[(ii)] Given $\delta >0$ with $\delta < \min \{ \frac{1}{2}, \delta _{0} \}$, where $\delta _{0}$ is given by Proposition \ref{prop:Laplacian:Infinity}, denote with $\underline{\delta}$ the pair $(\delta , -\delta)$ and set $\underline{\delta}-m=(\delta -m, -\delta -m)$ for any integer $m$. Define $\mathcal{C}_{ \underline{\delta} }$ as the space of pairs of a connection and a Higgs field on $V$ of the form $c + \xi$ with $\xi \in W^{1,2}_{ \underline{\delta}-1}$. 

\item[(iii)] The group $\mathcal{G}_{ \underline{\delta} }$ of gauge transformations is defined as the space of sections $g$ of $P \times_{\Ad} SO(3)$ such that $c+(d_{1}g)g^{-1} \in \mathcal{C}_{ \underline{\delta} }$.
\end{itemize}
\end{definition}

The fact that $\mathcal{G}_{ \underline{\delta }}$ is a group of continuous gauge transformations acting smoothly on $\mathcal{C}_{ \underline{\delta} }$ follows from Propositions \ref{prop:Weighted:Spaces:Singularities} and \ref{prop:Weighted:Spaces:Infinity}. Infinitesimal gauge transformations are sections of $V$ of class $\widetilde{W}^{2,2}_{ \underline{\delta} }$. Finally, by Propositions \ref{prop:Weighted:Spaces:Singularities} and \ref{prop:Weighted:Spaces:Infinity}, $(A,\Phi) \mapsto \ast F_{A} - d_{A}\Phi$ defines a smooth map $\Psi\co \mathcal{C}_{ \underline{\delta} } \ra L^{2}_{ \underline{\delta}-2}$.

The moduli space $\mathcal{M}_{n,k}$ is $\Psi ^{-1}(0)/\mathcal{G}_{\underline{\delta}}$. We will see below that the only singularities of $\mathcal{M}_{n,k}$ arise from reducible monopoles in $\mathcal{C}_{\underline{\delta}}$. Here a pair $(A,\Phi)$ is said to be \emph{reducible} if $V \simeq \underline{\R} \oplus M$ for an $SO(2)$--bundle $M \ra (\RS) \setminus S$ and $(A,\Phi)$ is induced by an abelian monopole on $M$. It is therefore important to understand when reducible monopoles exist. Denote by $c_{v,b}$ the abelian flat monopole $\left( d+ ib\, dt,v \right)$ and with $c_{p}$ the periodic Dirac monopole of charge $1$ with singularity at $p \in \RS$, \cf Definition \ref{def:Periodic:Dirac:Monopole}. Recall that we defined $k=\frac{k_{\infty}+n}{2} \in \Z_{\geq 0}$ as the non-abelian charge of the $SO(3)$--pair $(A,\Phi) \in \mathcal{C}_{ \underline{ \delta } }$.

\begin{lemma}\label{lem:Reducible:Pairs}
If $n<k$ every monopole in $\mathcal{C}_{ \underline{ \delta } }$ is irreducible.

If $n\geq k$, reducible monopoles in $\mathcal{C}_{ \underline{ \delta } }$ are in one to one correspondence with subsets $\{ p_{i_{1}}, \ldots, p_{i_{k}} \}$ of $S=\{ p_{1},\ldots ,p_{n}\}$ of cardinality $k$ and such that $p_{i_{1}}+ \ldots + p_{i_{k}} = \frac{1}{2}\left( \sum _{i=1}^{n}{p_{i}} + k_{\infty}q \right)$ in $\R^{2} \times \R / 2\pi\Z$.

After reordering the $p_{i}$'s if necessary, assume that $\{ p_{1},\ldots ,p_{k}\}$ satisfies this condition. Then the unique reducible monopole corresponding to this choice is
\[
c_{v,b} + \sum _{i=1}^{k}{c_{p_{i}}}-\sum _{i=k+1}^{n}{c_{p_{i}}}.
\]
\proof
If $(A,\Phi )\in\mathcal{C}_{ \underline{ \delta } }$ is a reducible monopole then $\Phi =\varphi\, \hat{\sigma}$ for a harmonic function $\varphi$ on $(\RS) \setminus S$ with prescribed behaviour at the punctures and at infinity. Here $\hat{\sigma}$ is the trivialising unit-norm section of the first factor in the decomposition $V \simeq \underline{\R} \oplus M$. After possibly reordering the $p_{i}$'s, $\varphi$ is of the form $\varphi =v+\sum _{i=1}^{n'}{G_{p_{i}}}-\sum _{i=n'+1}^{n}{G_{p_{i}}}$ for some $0 \leq n' \leq n$.

To conclude, use Lemmas \ref{lem:Asymptotics:Periodic:Dirac:Higgs:Field}, \ref{lem:Asymptotics:Periodic:Dirac:Connection} and \ref{lem:Asymptotics:Periodic:Dirac:Translations} to compare the asymptotics of the sum of Dirac monopoles $c_{v,b}+\sum _{i=1}^{n'}{c_{p_{i}}}-\sum _{i=n'+1}^{n}{c_{p_{i}}}$ with the boundary conditions of Definition \ref{def:Boundary:Conditions:SO(3)}: $n'=k$ because the charge at infinity has to be $2k-n = k_{\infty}$ and $p_{1} + \ldots + p_{k} = \frac{1}{2}\left( \sum _{i=1}^{n}{p_{i}} + k_{\infty}q \right)$ for the terms of order $\frac{1}{r}$ to coincide.
\qed
\end{lemma}

\subsection{The deformation complex}

Let $(A,\Phi) = c + \xi \in \Psi ^{-1}(0) \subset \mathcal{C}_{ \underline{ \delta } }$ be a solution to the Bogomolny equation \eqref{eqn:Bogomolny}. In order to prove that $\mathcal{M}_{n,k}$ is a smooth manifolds in a neighbourhood of $(A,\Phi)$ we have to show that:
\begin{itemize}
\item[(i)] The deformation complex \eqref{eqn:Deformation:Complex} defines a Fredholm complex $\widetilde{W}^{2,2}_{ \underline{ \delta } } \ra W^{1,2}_{ \underline{ \delta }-1 } \ra L^{2}_{ \underline{ \delta }-2 }$.
\item[(ii)] If $(A,\Phi)$ is irreducible, \ie $d_{1}$ is injective, then $d_{2}$ is surjective.
\end{itemize}

We will need the following elliptic regularity result for the Laplacians of the deformation complex.

\begin{lemma}\label{lem:Elliptic:Regularity}
Let $(A,\Phi) = c + \xi \in \mathcal{C}_{ \underline{\delta} }$. Then there exists $\sigma, R$ and $C$ depending on $\xi$ such that
\[
\| u \| _{\widetilde{W}^{2,2}_{ \underline{\delta} }} \leq C \left( \| DD^{\ast}u \| _{L^{2}_{ \underline{\delta}-2 }} + \| u \| _{L^{2}(K_{\sigma, R})} \right)
\]
for all $u \in \Omega \left( (\RS) \setminus S;V \right)$.
\proof
Denote with $D_{0}$ the Dirac operator \eqref{eqn:Dirac:Operator} twisted by the background pair $c$. By Lemmas \ref{lem:Embeddings:Products:Singularity} and \ref{lem:Multiplication:Embedding:Infinity}
\[
\| DD^{\ast}u - D_{0}D_{0}^{\ast}u \| _{L^{2}_{\underline{\delta}-2}} \leq C \| \xi \| _{W^{1,2}_{\underline{\delta}-1}}\| u \| _{\widetilde{W}^{2,2}_{\underline{\delta}}}.
\]
Choose $\sigma$ and $R$ so that $\| \xi|_{B_{\sigma}(p_{i})\setminus \{ p_{i}\} } \| _{W^{1,2}_{\underline{\delta}-1}}$ and $\| \xi|_{U_{R} \times \Sph ^{1}} \| _{W^{1,2}_{\underline{\delta}-1}}$ are sufficiently small. From the estimates in Propositions \ref{prop:Dirichlet:Problem:Weighted:Spaces:Singularity}
 and \ref{prop:Laplacian:Infinity} we deduce
\[
\| u \| _{\widetilde{W}^{2,2}_{ \underline{\delta} }} \leq C \left( \| DD^{\ast}u \| _{L^{2}_{ \underline{\delta}-2 }} + \| u|_{K_{\sigma, R}} \| _{\widetilde{W}^{2,2}_{\underline{\delta}}} \right).
\]
Therefore to prove the Lemma it is enough to show that for all compact sets $K' \subset K \subset X^{\ast}$, there exists $C=C(K,K',\xi)$ such that
\begin{equation}\label{eqn:Elliptic:Regularity}
\| u \| _{W^{2,2}(K')} \leq C \left( \| DD^{\ast}u \| _{L^{2}(K)} + \| u \| _{L^{2}(K)} \right).
\end{equation}
Here $W^{2,2}$ is the unweighted covariant Sobolev norm
\[
\| u \| ^{2}_{W^{2,2}} = \| u \| ^{2}_{L^{2}} + \| \nabla _{A}u \| ^{2}_{L^{2}} + \| [\Phi,u] \| ^{2}_{L^{2}} + \| \nabla _{A}(D^{\ast}u) \| ^{2}_{L^{2}} + \| [\Phi,D^{\ast}u] \| ^{2}_{L^{2}}.
\]
Choose a cut-off function $\chi$ supported on $K$ and such that $\chi \equiv 1$ on $K'$. Using the Weitzenb\"ock formula for $DD^{\ast}$, we have
\begin{align*}
\int{ \chi ^{2} \left( |\nabla_{A}u|^{2} + [\Phi,u]^{2} \right) } &+2 \int{\chi \langle \nabla _{A}u , \nabla \chi \otimes u \rangle} + \int{ \langle \Psi \cdot u , \chi ^{2}u \rangle }\\ &{}= \int{ \langle DD^{\ast}u, \chi ^{2}u \rangle } \leq \| DD^{\ast}u \| _{L^{2}{(K)}} \| u \| _{L^{2}(K)}
\end{align*}
where $\Psi = \ast F_{A}-d_{A}\Phi$. Now use Young's inequality with $\varepsilon >0$ to estimate
\[
\int{\chi \langle \nabla _{A}u , \nabla \chi \otimes u \rangle} \leq \varepsilon ^{2} \int{ \chi ^{2}|\nabla _{A}u|^{2} } + \frac{1}{\varepsilon ^{2}} \| u \| ^{2}_{L^{2}(K)}.
\]
and, together with H\"older's inequality,
\[
\int{|\Psi|\, |\chi u|^{2} } \leq \| \Psi \| _{L^{2}(K)} \| \chi u \| _{L^{2}}^{\frac{1}{2}}\| \chi u\| _{L^{6}}^{\frac{3}{2}} \leq \varepsilon ^{2} \| \chi u \| ^{2}_{L^{6}} + C_{\varepsilon} \| \Psi \| _{L^{2}(K)}^{4}\| \chi u \| _{L^{2}}^{2}.
\]
The Sobolev embedding $W^{1,2} \hookrightarrow L^{6}$ now implies
\[
\int{|\Psi|\, |\chi u|^{2} } \leq \varepsilon ^{2}\| \chi \nabla _{A}u \| ^{2}_{L^{2}} + C_{\varepsilon}(1+\|\Psi \| ^{4}_{L^{2}(K)}) \| u \| ^{2}_{L^{2}(K)}
\]
Choosing $\varepsilon$ small enough we obtain
\[
\| \nabla _{A}u \| ^{2}_{L^{2}(K')} + \| [\Phi,u] \| ^{2}_{L^{2}(K')} \leq C \| DD^{\ast}u \| ^{2}_{L^{2}(K)} + C(1+\|\Psi \| ^{4}_{L^{2}(K)})\| u \| ^{2}_{L^{2}(K)}.
\]
The second order estimate is obtained in a similar way, restricting to an even smaller compact set $K'' \subset K'$ and using the Weitzenb\"ock formula for the operator $D^{\ast}D$.

Thus we obtained \eqref{eqn:Elliptic:Regularity} for a constant $C$ depending on $\| d_{A}\Phi \| _{L^{2}(K)}$ and $\| \Psi \| _{L^{2}(K)}$. To conclude observe that, since $(A,\Phi)=c+\xi$, with $c$ smooth and $\xi \in W^{1,2}_{ \underline{\delta}-1 }$ (in particular $\xi \in W^{1,2}_{loc}$), $\| d_{A}\Phi \| _{L^{2}(K)}$ and $\| \Psi \| _{L^{2}(K)}$ are bounded in terms of $K$, the background $c$ and $\| \xi \| _{W^{1,2}_{ \underline{ \delta }-1 }}$.
\qed
\end{lemma}

\subsection{Slice to the action of the gauge group}

\begin{prop}\label{prop:Laplacian:Isomorphism}
The operator $DD^{\ast}\co \widetilde{W}^{2,2}_{ \underline{\delta} } \ra L^{2}_{ \underline{\delta}-2 }$ is Fredholm. If $(A,\Phi)$ is irreducible then $DD^{\ast}$ is an isomorphism.
\proof
For all $\varepsilon >0$ we can find $\sigma , R >0$ such that $\| \xi |_{B_{\sigma }(p_{i})} \| _{W^{1,2}_{\rho, \delta -1}} < \varepsilon$ and $\| \xi |_{U_{R} \times \Sph ^{1}} \| _{W^{1,2}_{\omega, -\delta -1}} < \varepsilon$. By choosing $\varepsilon >0$ sufficiently small, Propositions \ref{prop:Dirichlet:Problem:Weighted:Spaces:Singularity} and \ref{prop:Laplacian:Infinity} and the continuity of the products in Lemmas \ref{lem:Embeddings:Products:Singularity} and \ref{lem:Multiplication:Embedding:Infinity} imply that the Dirichlet problem for the operator $DD^{\ast}$ on $B_{\sigma }(p_{i}) \setminus \{ p_{i} \}$ and $U_{R} \times \Sph ^{1}$ is an isomorphism. Thus we obtain inverses of $DD^{\ast}$ in a neighbourhood of the singularities and at infinity by solving Dirichlet problems with vanishing boundary conditions. The fact that $DD^{\ast}$ is a Fredholm operator now follows by gluing these inverses with a parametrix on the compact set $K_{\sigma ,R}$, \cf for example R\aa de's \cite[Lemma 3.2]{Rade:Global:Theory}.

To show that $DD^{\ast}$ is an isomorphism if $(A, \Phi)$ is irreducible, we proceed in three steps.

\begin{enumerate}
\item By the Weitzenb\"ock formula Lemma \ref{lem:Weitzenbock}, if $(A,\Phi)$ is irreducible than $DD^{\ast}$ is injective. Indeed,
\[
0=\int{\langle \nabla _{A}^{\ast}\nabla _{A}u-\ad^{2}(\Phi)u,u \rangle } = \int{|\nabla _{A}u|^{2} + |[\Phi,u]|^{2}}
\]
If $\delta$ is in the range specified the integration by parts can be justified using a sequence of cut-off functions converging to $1$. Observe also that, since $DD^{\ast}$ is injective, Lemma \ref{lem:Elliptic:Regularity} and a standard argument by contradiction using Rellich's compactness imply that there exists a constant $C >0$ such that $\| u \| _{\widetilde{W}^{2,2}_{\underline{\delta}}} \leq C \| DD^{\ast}u \| _{L^{2}_{ \underline{\delta}-2 }}$. The Proposition will follow from the fact that the index of $DD^{\ast}$ vanishes.

\item Choosing $\sigma$ sufficiently small and $R$ sufficiently large we can deform $(A,\Phi)$ into a new pair $c'=(A', \Phi')=(A,\Phi) + \chi \xi$ which coincides with the background pair $c$ outside of $K_{\sigma, R}$. Here $\chi$ is a cut-off function with support in $K_{\sigma, R}$. By the compactness of the products in Lemmas \ref{lem:Embeddings:Products:Singularity} and \ref{lem:Multiplication:Embedding:Infinity} the index of $D_{c'}D_{c'}^{\ast}$ and $DD^{\ast}$ coincide. Moreover, since $D_{c'}D_{c'}^{\ast}$ is of the form $DD^{\ast} + T$, where the operator norm of $T$ is controlled by $\| (1-\chi)\xi \| _{W^{1,2}_{\underline{\delta}-1}}$, it follows from Step 1 that $D_{c'}D_{c'}^{\ast}$ remains injective provided $\sigma$ and $R$ are chosen so that $\| (1-\chi)\xi \| _{W^{1,2}_{\underline{\delta}-1}}$ is sufficiently small.

\item For notational convenience we drop the subscript $\,_{c'}$ in the rest of the proof. It remains to be shown that $DD^{\ast}$ is surjective. Start considering the map $DD^{\ast}\co W^{2,2}_{ \underline{ \delta } } \ra L^{2}_{ \underline{ \delta } -2 }$. Since $\delta$ and $-\delta$ are non-exceptional weights for $DD^{\ast}$ close to the singularities and at infinity, standard theory of weighted Sobolev spaces implies that the cokernel of this map is identified with the kernel of $DD^{\ast}$ in $L^{2}_{ \underline{\delta}^{\ast} }$, where $\underline{\delta}^{\ast} =(-\delta-1, \delta )$ (\cf for example \cite[Theorem 10.2.1]{Pacard:Connected:Sums}). Denote this finite dimensional space by $\ker{(DD^{\ast})_{ \underline{\delta}^{\ast} } }$. We claim that there is an injective map $\ker{(DD^{\ast})_{ \underline{\delta}^{\ast} }} \ra \R ^{4(n+1)}$. This can be shown by solving the Dirichlet problem on balls $B_{\sigma}(p_{i}) \setminus \{ p_{i} \}$ and on $U_{R} \times \Sph ^{1}$ (for some small $\sigma$ and large $R$) to write any element $u=u_{0}+u_{1}dx+u_{2}dy+u_{3}dt \in \ker{(DD^{\ast})_{ \underline{\delta}^{\ast} } }$ as
\begin{alignat*}{2}
u_{\alpha}|_{B_{\sigma}(p_{i})}= \frac{\lambda_{\alpha, i}}{\rho}\,\hat{\sigma} + u'_{\alpha, i}\qquad & u|_{U_{R} \times \Sph ^{1}}=\lambda_{\alpha, \infty}(\log{r})\, \hat{\sigma} + u'_{\alpha, \infty}
\end{alignat*}
with $u'_{\alpha, i} \in \widetilde{W}^{2,2}_{\rho, \delta }$ and $u'_{\alpha, \infty} \in \widetilde{W}^{2,2}_{\omega,-\delta }$, $\alpha = 0,1,2,3$. Here $\hat{\sigma}$ stands for the trivialising section of the diagonal factor in the decomposition $V \simeq \underline{\R} \oplus M$, with $M=H_{p_{i}}$ over $B_{\sigma}(p_{i}) \setminus \{ p_{i} \}$ and $M=L_{v,b} \otimes L_{q}^{k_{\infty}}$ on $U_{R} \times \Sph ^{1}$. Since $\widetilde{W}^{2,2}_{ \underline{ \delta} }$ is an extension of $W^{2,2}_{ \underline{ \delta } }$ by a $4(n+1)$--dimensional space and $DD^{\ast}\co \widetilde{W}^{2,2}_{ \underline{\delta} } \ra L^{2}_{ \underline{\delta}-2 }$ remains injective by Step 2, we conclude that $DD^{\ast}$ is an isomorphism.
\qedhere
\end{enumerate}
\end{prop}

\begin{remark}\label{rmk:Cokernel:Laplacian}
When $(A,\Phi)$ is reducible $DD^{\ast}$ has a $4$--dimensional cokernel. This is a consequence of the parabolicity of $\RS$ (\ie the fact that every Green's function changes sign): a necessary condition to solve $\triangle u =f$ on $\RS$ with $\nabla u \in L^{2}$ is that $f$ has mean value zero.
\end{remark}

Observe that if $(A,\Phi)$ is a monopole $DD^{\ast}$ acting on $\Omega ^{0} \oplus \{ 0 \} \subset \Omega$ coincides with $d_{1}^{\ast}d_{1}$. Hence standard theory \cite[Chapter 4]{Donaldson:Kronheimer} now implies that
\[
S_{(A,\Phi),\epsilon} = \left\{ (A,\Phi)+(a,\psi) \; | \; d_{1}^{\ast}(a,\psi) = 0 ,\, \| (a,\psi) \| _{W^{1,2}_{ \underline{\delta} -1 }} < \epsilon \right\}
\]
is a local slice for the action of $\mathcal{G}_{ \underline{ \delta } }$ on $\mathcal{C}_{ \underline{\delta} }$.

\subsection{Fredholm property of the Dirac operator $D$}

\begin{prop}\label{prop:Deformation:Operator:Fredholm}
Let $(A,\Phi) \in \mathcal{C}_{ \underline{\delta} }$ be a solution to the Bogomolny equation. Then $D\co W^{1,2}_{ \underline{\delta}-1 } \ra L^{2}_{ \underline{\delta}-2 }$ is a Fredholm operator, surjective when $(A,\Phi)$ is irreducible.
\proof
The cokernel of $D \co W^{1,2}_{ \underline{\delta}-1 } \ra L^{2}_{ \underline{\delta}-2 }$ is identified with $\ker{ D^{\ast} } \cap L^{2}_{ \underline{\delta}^{\ast} } \subset  \ker{ DD^{\ast} } \cap L^{2}_{ \underline{\delta}^{\ast} }$, where $\underline{\delta}^{\ast} =(-\delta-1, \delta )$. By Proposition \ref{prop:Laplacian:Isomorphism} and Remark \ref{rmk:Cokernel:Laplacian} we already know that this vanishes when $(A,\Phi)$ is irreducible and is $4$--dimensional otherwise.

It remains to show that the image of $D \co W^{1,2}_{ \underline{\delta}-1 } \ra L^{2}_{ \underline{\delta}-2 }$ is closed and the kernel finite dimensional.  Both statements follow by standard arguments from the estimate ($K$ a compact subset of $X$)
\begin{equation}\label{eqn:Fredholm:Estimate}
\| \xi \| _{W^{1,2}_{ \underline{\delta}-1 }} \leq C \left( \| D\xi \| _{L^{2}_{ \underline{\delta}-2 }} + \| \xi \| _{L^{2}(K)} \right) .
\end{equation}

From the estimates in Propositions \ref{prop:Dirichlet:Problem:Weighted:Spaces:Singularity} and \ref{prop:Laplacian:Infinity} and Lemma \ref{lem:Elliptic:Regularity} we deduce
\begin{equation}\label{eqn:Elliptic:Regularity:D}
\| \xi \| _{W^{1,2}_{ \underline{\delta}-1 }} \leq C \left( \| D\xi \| _{L^{2}_{ \underline{ \delta}-2 }} + \| \xi \| _{L^{2}_{ (\delta-1,-\delta-1) }} \right) .
\end{equation}
We can also fix $\sigma, R>0$ as small, large as needed and deform $(A,\Phi)$ to $(A',\Phi')$ so that it coincides with the model Dirac monopoles on $B_{2\sigma}(p_{i})$ and $U_{R} \times \Sph ^{1}$. By Lemmas \ref{lem:Embeddings:Products:Singularity} and \ref{lem:Multiplication:Embedding:Infinity} such a modification changes $D$ by a compact operator. Moreover, \eqref{eqn:Elliptic:Regularity:D} continues to hold.

We proceed with the proof of \eqref{eqn:Fredholm:Estimate}. Using a cut-off function we write $\xi = \xi _{1}+\xi _{2}$ with $\xi _{1}$ supported on $B_{2\sigma}(p_{i})$ and $U_{R} \times \Sph ^{1}$ and $\xi _{2}$ supported on $K_{\sigma, 2R}$. Notice that if $\chi$ is a compactly supported function, then
\[
\| D(\chi \xi) \| _{L^{2}_{\underline{\delta}-2}} \leq C \left( \| D\xi \| _{L^{2}_{ \underline{ \delta}-2 }} + \| \xi \| _{L^{2}(\text{spt}\,\chi)} \right).
\]
With $\xi = \xi _{2}$, \eqref{eqn:Elliptic:Regularity:D} is in fact equivalent to \eqref{eqn:Fredholm:Estimate}. Thus we reduced the problem to prove \eqref{eqn:Fredholm:Estimate} assuming that $\xi$ is supported on $B_{2\sigma}(p_{i})$ and $U_{R} \times \Sph ^{1}$. Since $(A',\Phi')$ is reducible on the support of $\xi$, we decompose $\xi = \xi _{D}+\xi _{T}$ and study separately the two terms.

\begin{enumerate}
\item On the diagonal part we can appeal to standard theory for the Laplacian in weighted Sobolev spaces. First, by Propositions \ref{prop:Dirichlet:Problem:Weighted:Spaces:Singularity} and \ref{prop:Laplacian:Infinity} there exists a unique solution $u$ of
\[
\begin{dcases}
\triangle u = D\xi _{D} & \text{on } \bigcup_{i=1}^{n}{B_{2\sigma}(p_{i})} \cup \left( U_{R} \times \Sph ^{1} \right) \\
u =0 & \text{on } \partial K_{\sigma, 2R}
\end{dcases}
\]
with $\| D^{\ast}u \| _{W^{1,2}_{\underline{\delta}-1}} \leq C \| D\xi \| _{L^{2}_{\underline{\delta}-2}}$. Thus $\xi _{D}=D^{\ast}u + \eta$ with $D\eta =0$. Fix a cut-off function which vanishes in a neighbourhood of $\partial K_{\sigma, 2R}$. Since $\delta -1$ and $-\delta-1$ are non-exceptional weights for the Laplacian and there are no harmonic functions in $W^{1,2}_{\underline{\delta}-1}$ vanishing on $\partial K_{\sigma, 2R}$
\[
\| \chi \eta \| _{L^{2}_{\underline{\delta}-1}} \leq C \| \triangle (\chi\eta ) \| _{L^{2}_{\underline{\delta}-2}} \leq C \| \eta \| _{W^{1,2}(\text{spt}\,\chi)}
\]
(\cf \cite[Proposition 6.2.2]{Pacard:Connected:Sums}). Therefore by standard elliptic estimates
\[
\| \eta \| _{L^{2}_{\underline{\delta}-1}} \leq C \| \eta \| _{W^{1,2}(\text{spt}\,\chi)} \leq C \| \eta \| _{L^{2}(K)}
\]
with $K = \bigcup_{i=1}^{n}{ B_{2\sigma}(p_{i}) \setminus B_{\sigma}(p_{i}) } \cup \left( B_{R+1}\setminus B_{R} \right) \times \Sph ^{1}$.

\item In order to prove the estimate for the off-diagonal component on the exterior domain $U_{R} \times \Sph ^{1}$, we exploit the Bochner formula \eqref{eqn:Bochner}. We showed in Step 4 of the proof of Proposition \ref{prop:Laplacian:Infinity} that \eqref{eqn:Bochner} implies
\[
\int{\omega ^{2\delta +2}\left( |\nabla _{A}\xi _{T}|^{2} + |[\Phi, \xi _{T}]|^{2} \right) } \leq C \left( \int{ \omega ^{2\delta +2}|D\xi _{T}|^{2}  } + \int{ \omega ^{2\delta}|\xi _{T}|^{2}  } \right) .
\]
The integrations by parts are justified because $\xi \in L^{2}_{\omega, -\delta -1}$. Since $|[\Phi, \xi _{T}]| \geq c|\xi _{T}|$ by \eqref{eqn:Growth:Higgs:Field:Infinity}, we can choose $R$ large enough so that $cR^{2}>C$ and therefore
\[
(cR^{2} - C) \int{ \omega^{2\delta} | \xi _{T} |^{2} } \leq C \int{ \omega^{2\delta +2}| D\xi_{T} |^{2} }.
\]

\item Since $(A',\Phi')$ coincides with an Euclidean Dirac monopole of mass $0$ on the ball $B_{2\sigma}(p_{i})$, $d_{A}(\rho \Phi)=0$ in this region. In particular, the Weitzenb\"ock formulas of Lemma \ref{lem:Weitzenbock} imply that $D(\rho D^{\ast}\xi _{T} )=D^{\ast}(\rho D\xi _{T})$. An integration by parts (justified because $\xi _{T} \in W^{1,2}_{\rho,\delta-1}$) yields
\[
\int{ \rho ^{-2\delta +1} |D^{\ast}\xi _{T} |^{2} } -  \int{ \rho ^{-2\delta +1} |D \xi _{T} |^{2} } = 2\delta  \int{ \rho ^{-2\delta} \langle D^{\ast}\xi _{T}-D\xi _{T}, d\rho \cdot \xi _{T} \rangle }.
\]
Now use the algebraic identity $2[\Phi,\xi _{T}]=D\xi _{T}-D^{\ast}\xi _{T}$:
\[
4|[\Phi, \xi _{T}]|^{2}=|D^{\ast}\xi _{T}|^{2}-|D\xi _{T}|^{2}+4\langle D\xi _{T}, [\Phi, \xi _{T}] \rangle
\]
and therefore
\[
\int{ \rho ^{-2\delta +1}|[\Phi, \xi _{T}]|^{2}  } = -\delta \int{ \rho ^{-2\delta} \langle [\Phi,\xi _{T}], d\rho \cdot \xi _{T} \rangle } + \int{ \rho ^{-2\delta+1}\langle D\xi _{T}, [\Phi, \xi _{T}] \rangle }.
\]
Finally, by the Cauchy--Schwarz inequality
\[
\int{ \rho ^{-2\delta +1}|[\Phi, \xi _{T}]|^{2}  } \leq \delta^{2} \int{ \rho ^{-2\delta-1} |\xi _{T}|^{2} } + \int{ \rho ^{-2\delta+1} |D\xi _{T}|^{2} }.
\]
Conclude using $\delta <\frac{1}{2}$ and $|[\Phi,\xi _{T}]|=\frac{1}{2}\rho ^{-1}|\xi _{T}|$.\qedhere
\end{enumerate}
\end{prop}

In view of Propositions \ref{prop:Laplacian:Isomorphism} and \ref{prop:Deformation:Operator:Fredholm} and the discussion of irreducibility in Lemma \ref{lem:Reducible:Pairs}, standard theory \cite[Chapter 4]{Donaldson:Kronheimer} implies that the moduli space $\mathcal{M}_{ \underline{\delta} }= \Psi ^{-1}(0)/\mathcal{G}_{\underline{\delta} }$ is a smooth manifold for generic choices of $p_{1}, \ldots, p_{n}, q \in X$ whenever it is non-empty.

\subsection{The $L^{2}$--metric}

The final task is to show that the $L^{2}$--metric is well-defined on $\mathcal{M}_{ \underline{\delta} }$. We will need the following lemma on the decay at infinity of monopoles in $\mathcal{C}_{\underline{\delta} }$.

\begin{lemma}\label{lem:Decay:Construction:Moduli:Space}
Let $(A,\Phi) = c+ \xi \in \mathcal{C}_{ \underline{\delta} }$ be an irreducible solution to the Bogomolny equation. Then there exist $R>0$ and $g \in \mathcal{G}_{ \underline{\delta} }$ such that on the exterior region $U_{R} \times \Sph ^{1}$ we have $g(A,\Phi)=c+\xi'$ with $\xi' \in W^{1,2}_{\omega, -\delta -1}$ and $\xi '_{D} = O(r^{-\delta -1})$, $\xi'_{T}=O(r^{\mu})$ for all $\mu \in \R$.
\proof
The line of proof follows \cite[Lemma 5.3]{Biquard:Boalch}.
\begin{itemize}
\item[Step 1.] First we put $(A,\Phi)$ in ``Coulomb gauge'' with respect to the background pair $c$ near infinity. Fix $R_{0}>0$ and a cut-off function $\chi_{R_{0}} \equiv 1$ on $B_{R_{0}} \times \Sph ^{1}$ and $\chi _{R_{0}} \equiv 0$ on $U_{2R_{0}} \times \Sph ^{1}$. Define a new pair $c'=(A',\Phi')=c+\chi_{R_{0}}\xi$. Then $c' \equiv c$ on $U_{2R_{0}} \times \Sph ^{1}$. As in Proposition \ref{prop:Laplacian:Isomorphism}, we can choose $R_{0}$ sufficiently large so that $d_{1}^{\ast}d_{1}\co \widetilde{W}^{2,2}_{\delta _{1}, \delta _{2}} \ra L^{2}_{\delta _{1}-2, \delta _{2}-2}$ remains invertible.

For all $R > R_{0}$ consider the pair $c' + \xi_{R}$ defined by $\xi _{R}=(1-\chi _{R})\xi$. Here $\chi_{R}$ is a cut-off function with the same properties of $\chi_{R_{0}}$ but with $R$ in place of $R_{0}$. The Implicit Function Theorem implies that, choosing $R$ large so that $\| \xi _{R} \| _{W^{1,2}_{ \underline{\delta}-1 } }$ is sufficiently small, there exists $g \in \mathcal{G}_{ \underline{\delta} }$ such that $g(c'+\xi _{R})=c'+\xi '$ with $\xi ' \in W^{1,2}_{ \underline{\delta}-1 }$ and $d_{1}^{\ast}\xi'=0$.

Since $c' + \xi_{R} = c+\xi$ on $U_{2R} \times \Sph ^{1}$, restricting to this exterior region $\xi'$ is a solution to $D\xi' + \xi' \cdot \xi' =0$. Here $D$ is the Dirac operator \eqref{eqn:Dirac:Operator} twisted by the background pair $c$.

\item[Step 2.] Renaming $\xi  =\xi'$, we reduced the problem to study the decay of solutions $\xi \in W^{1,2}_{\omega , -\delta -1}$ to $D\xi = -\xi \cdot \xi$. We start by proving an initial decay $\xi = O(r^{-\delta})$ and then improve to the required rate.

Apply $D^{\ast}$ to the equation and use the Weitzenb\"ock formula Lemma \ref{lem:Weitzenbock} to derive the differential inequality
\[
d^{\ast}d(|\xi|) \lesssim |d_{A}\Phi|\, |\xi| + \left( |\nabla _{A}\xi| + |[\Phi, \xi]| \right)\, |\xi|.
\]
Hence $|\xi| \in W^{1,2}$ is a subsolution to $dd^{\ast}u \leq (A_{1}+A_{2}) u$, where $A_{1}=|d_{A}\Phi| \in L^{\infty}$ and $A_{2}=|\nabla _{A}\xi| + |[\Phi, \xi]| \in L^{2}$. Then Moser iteration on a $3$--ball $B_{1}(p)$ centred at any point $p \in U_{3R} \times \Sph ^{1}$ as in \cite[Theorem 8.17]{Gilbarg:Trudinger} yields
\[
\sup _{B_{\frac{1}{2}}(p)}{|\xi|} \leq C\| \xi \| _{L^{2}(B_{1}(p))} \leq Cr^{-\delta}\| \xi \| _{L^{2}_{\omega, -\delta-1}}
\]
for a constant $C$ depending on the $L^{\infty}$--norm of $A_{1}$ and $\| A_{2} \| _{L^{2}}$. Here we used that $\omega \sim \omega (p)\sim r$ in $B_{1}(p)$. Hence $|\xi| \leq Cr^{-\delta}$ on $U_{3R} \times \Sph ^{1}$ for a constant $C$ depending on the background $c$, $R$ and $\| \xi \| _{W^{1,2}_{\underline{\delta} -1}}$.

\item[Step 3.] Recall that the background pair $c$ is abelian on $U_{R} \times \Sph ^{1}$. We  decompose $\xi = \xi _{D} + \xi _{T}$ into diagonal and off-diagonal part and exploit the fact that $\xi \in W^{1,2}_{\omega, -\delta -1} \Rightarrow \omega ^{-\delta +1}\xi _{T} \in W^{1,2}$ to improve the decay of $\xi_{T}$ first in an integral sense, then as a pointwise statement.

In order to justify the integrations by parts it is necessary to introduce a sequence of cut-off functions $\chi _{i}$ vanishing in a neighbourhood of infinity, such that $|d\chi _{i}| \leq \frac{2}{r}$ and converging to $1$ as $i \ra \infty$. Set $\xi _{i}= \chi _{i}\, \xi_{T}$; then $D\xi _{i} = d\chi _{i}\cdot \xi_{T} - \xi \cdot \xi _{i}$.

If $\xi _{T} \in W^{1,2}_{\omega , \mu -1}$ then $D\xi _{i} \in L^{2}_{\omega, \mu -2 +\delta}$ since $\omega ^{-\mu +1}\xi _{i}, \omega ^{-\mu +1}\xi _{T}, \omega ^{-\delta}\xi \in W^{1,2}$ and $\delta > -1$. Moreover, $\xi _{i} \in L^{2}_{\omega, \mu-1+\delta}$ because $\delta >-1$. The a priori estimate of Proposition \ref{prop:Laplacian:Infinity} now implies $\xi _{i} \in W^{1,2}_{\omega , \mu -1-\delta}$---an improvement. By iterating and letting $i \ra \infty$, we conclude that $\xi _{T} \in W^{1,2}_{\omega , \mu -1}$ for all $\mu \in (-\infty, -\delta ].$
\item[Step 4.] We repeat the argument of Step 2 with the equation $D\xi _{T} = -\xi \cdot \xi _{T}$. We have a differential inequality
\[
d^{\ast}d(|\xi_{T}|) \lesssim  |d_{A}\Phi|\,|\xi_{T}| + (|\nabla _{A}\xi| + |[\Phi, \xi]|)\, |\xi_{T}| + (|\nabla _{A}\xi_{T}| + |[\Phi, \xi_{T}]|)\, |\xi|
\]
of the form $d^{\ast}d u \lesssim A_{1}u + A_{2}u +  f$, where $u=|\xi_{T}| \in W^{1,2}$, $A_{1}=|d_{A}\Phi| \in L^{\infty}$, $A_{2}=|\nabla _{A}\xi| + |[\Phi, \xi]| \in L^{2}$ and $f=(|\nabla _{A}\xi_{T}| + |[\Phi, \xi_{T}]|)\, |\xi| \lesssim |\nabla _{A}\xi_{T}| + |[\Phi, \xi_{T}]| \in L^{2}$ by Step 2. Moser iteration and Step 3 yield $|\xi _{T}| =O(r^{\mu})$ on $U_{4R} \times \Sph ^{1}$ for all $\mu \in \R$. 
\item[Step 5.] The diagonal part $\xi _{D} \in W^{1,2}_{\omega , -\delta -1}$ is a solution to the equation
\[
\triangle \xi _{D} = D^{\ast}(\xi _{T} \cdot \xi _{T}) \in L^{2}_{\omega , \mu -2}
\]
for all $\mu \in \R$. By elliptic regularity $\xi _{D} \in W^{2,2}_{\omega, -\delta -1}$ and an argument analogous to the proof of Lemma \ref{lem:Multiplication:Embedding:Infinity}.(ii) yields the weighted Sobolev embedding $W^{2,2}_{\omega, -\delta -1} \hookrightarrow \omega ^{-\delta -1}C^{0}$. \qed
\end{itemize}
\end{lemma}
\begin{remark}\label{rmk:Decay:Construction:Moduli:Space}

\begin{itemize}
\item[(i)] In fact we could say a bit more: $|\xi _{D}|=O(r^{-2})$, the rate of decay of $L^{2}_{\omega, -\delta -1}$--harmonic functions on $\RS$.
\item[(ii)] An analogous argument yields the same decay for solutions to $D\xi = 0$.
\end{itemize}
\end{remark}

We summarise what we have proved so far in the following theorem.

\begin{thm}\label{thm:Construction:Moduli:Space}
Choose data $v,b,k_{\infty}, p_{1}, \ldots, p_{n},q$ defining the boundary conditions of Definition \ref{def:Boundary:Conditions:SO(3)}. Fix $\delta >0$ sufficiently small and suppose that the parameters $k_{\infty}, p_{1}, \ldots, p_{n},q$ are chosen so that every monopole $(A,\Phi) \in \Psi ^{-1}(0) \subset \mathcal{C}_{ \underline{\delta} }$ is irreducible. Then the moduli space $\mathcal{M}_{n,k}$ of $SO(3)$ periodic monopoles with non-abelian charge $k=\frac{k_{\infty}+n}{2}$, centre $q$ and singularities at $p_{1}, \ldots, p_{n}$ is a smooth manifold, provided it is non-empty. Moreover, the tangent space of $\mathcal{M}_{n,k}$ at a point $[(A,\Phi)]$ is identified with the $L^{2}$--kernel of $D$ and the $L^{2}$--metric is a hyperk\"ahler metric on $\mathcal{M}_{n,k}$.
\proof
In view of Proposition \ref{prop:Deformation:Operator:Fredholm}, only the last two statements need justification.

For the first, by \eqref{eqn:Elliptic:Regularity:D} it is enough to prove that if $\xi \in L^{2}$ satisfies $D\xi =0$ then $\xi \in L^{2}_{ \underline{\delta}-1 }$.
\begin{itemize}
\item[(i)] On a small ball $B_{\sigma}(p_{i})$, let $\hat{\xi}$ be the lift of $\xi$ to a $4$--ball as in Definition \eqref{eqn:Lift:Forms}. Then $\hat{\xi}$ is a solution to $\hat{D}\hat{\xi}=0$, where $\hat{D}$ is the Dirac operator twisted by the smooth connection $\hat{A}$ obtained from $(A,\Phi)$ as in \eqref{eqn:Lift:Connection}. By elliptic regularity $|\hat{\xi}|=\sqrt{\rho} |\xi|$ is bounded.
\item[(ii)] Near infinity we use Lemma \ref{lem:Decay:Construction:Moduli:Space} to write $(A,\Phi)=c+\eta$ with $\eta = O(r^{-\delta -1})$. Then $\xi$ is a solution to $D\xi +\eta \cdot \xi =0$, where $D$ is the Dirac operator \eqref{eqn:Dirac:Operator} twisted by the background pair $c$. It follows that $D\xi \in L^{2}_{\omega, -\delta -2}$ on $U_{R} \times \Sph ^{1}$ for some $R$ large enough. By Proposition \ref{prop:Laplacian:Infinity} we can write $\xi = \xi ' + D^{\ast}u$, where $u \in \widetilde{W}^{2,2}_{\omega, -\delta}$ and $\xi ' \in L^{2}$ with $D\xi'=0$. Since $c$ coincides with the model periodic Dirac monopole on $U_{R} \times \Sph ^{1}$, the diagonal component of $\xi'$ is an $L^{2}$ harmonic function and therefore $\xi'_{D} = O(r^{-2})$. On the other hand, by Remark \ref{rmk:Decay:Construction:Moduli:Space}.(ii)
 $\xi'_{T} = O(r^{\mu})$ for all $\mu \in \R$. 
\end{itemize}

Finally, the fact that the $L^{2}$--metric is hyperk\"ahler is an instance of a hyperk\"ahler quotient in infinite dimension, \cf \cite{Atiyah:Hitchin}. The only analytic point to be checked is that the equality
\[
\langle \xi , d_{1}u \rangle_{L^{2}} = \langle d_{1}^{\ast}\xi , u \rangle _{L^{2}}
\]
holds for $\xi \in W^{1,2}_{ \underline{\delta}-1 }$ and $u \in \widetilde{W}^{2,2}_{ \underline{\delta} }$. This can be verified by using a sequence of cut-off functions on $\RS$ converging to $1$.
\qed
\end{thm}

\section{The dimension of the moduli spaces}\label{sec:Dimension}

In order to conclude the proof of Theorem \ref{thm:Main:Theorem} it remains to calculate the dimension of the moduli space $\mathcal{M}_{n,k}$. In this section we prove the following index theorem.

\begin{thm}\label{thm:Index:Theorem}
Let $(A,\Phi)$ be a pair in $\mathcal{C}_{\underline{\delta}}$. Then the index of $D \co W^{1,2}_{\underline{\delta}-1} \ra L^{2}_{\underline{\delta}-2}$ is $4k-4$.
\end{thm}

The operator $D=\tau \slashed{D}_{A} + [\Phi, \cdot ]$ is of Callias-type, \ie it is a Dirac operator plus a potential. Index theorems for such operators on complete odd-dimensional manifolds have been obtained by  Callias \cite{Callias}, Anghel \cite{Anghel} and R\aa de \cite{Raade:Index}. The common requirement of all these results is that the potential term is non-degenerate at infinity. For example, if we assume that $(A,\Phi)$ is a periodic charge $k$ $SU(2)$--monopole without singularities and we let $D$ act on sections of the associated rank $2$ complex vector bundle $E$, R\aa de's result yields $\ind(D,E) = 2k$.

When we couple $D$ with the adjoint bundle, however, such non-degeneracy condition doesn't hold because $[\Phi, \cdot\,]$ has a $1$--dimensional kernel. One approach to go round this difficulty is given by Kottke \cite{Kottke} in the case of (smooth) monopoles on asymptotically conical complete $3$--manifolds. In our situation an additional complication arises from the presence of singularities.

We will give a direct computation of the index of $D$ using the excision principle, very much in the spirit of the calculation of the dimension of the moduli space of instantons on a $4$--manifold, \cf \cite[\S 7.1]{Donaldson:Kronheimer}. By the compactness properties of Lemmas \ref{lem:Embeddings:Products:Singularity} and \ref{lem:Multiplication:Embedding:Infinity} the index of $D$ is independent of the pair $(A,\Phi) \in \mathcal{C}_{\underline{\delta}}$. Thus we will carry out the computation of the index for an explicit smooth pair $(A,\Phi)$. This is constructed patching together a sum of periodic Dirac monopoles with an Euclidean charge $k$ monopole. Comparing the corresponding Dirac operators on $\RS$ and $\R ^{3}$, the excision principle allows to compute the index of $D$ as a sum of contributions from the different pieces: on one side, the index of the Dirac operator twisted by a (smooth) Euclidean monopole has been calculated by Taubes in \cite{Taubes:Stability:Yang:Mills:Theories}; on the other, making the mass of the monopole very large, one can understand the contribution of the sum of Dirac monopoles.

\subsection{Construction of a background pair}

As the first step in the proof of Theorem \ref{thm:Index:Theorem}, we give the explicit construction of a smooth pair $(A, \Phi)$. This can be taken to be the background pair in the definition of the moduli space $\mathcal{M}_{n,k}$ at the beginning of Section \ref{sec:Construction:Moduli:Space}.

Fix a collection of $n$ distinct points $S=\{ p_{1}, \ldots, p_{n} \}$ in $\R ^{2} \times \Sph ^{1}$ and set $X = \left( \R ^{2} \times \Sph ^{1} \right) \setminus S$. Choose an additional point $q \in X$ such that $2B=B_{2}(q)$ is contained with its closure in $X$ and set $X^{\ast}=X \setminus \{ q\}$. We write $X = B \cup U_{\ext}$, where $U_{\ext} = X \setminus \frac{1}{2}B$, and $X^{\ast}=B^{\ast} \cup U_{\ext}$, where $B^{\ast}=B \setminus \{ q \}$. Similarly, set $Y=\R ^{3} = B \cap U'_{\ext}$ and $Y^{\ast} = \R ^{3} \setminus \{ 0 \}= B^{\ast} \cup U'_{\ext}$, where $B=B_{1}(0) \subset \R ^{3}$ and $U'_{\ext} = \R ^{3} \setminus \frac{1}{2}B$. The notation suggests that we fix an identification of a neighbourhood of $q$ in $\RS$ with a neighbourhood of the origin in $\R ^{3}$.

Next we define $SO(3)$--bundles with connection and Higgs field on $X, X^{\ast}, Y$ and $Y^{\ast}$. Over $B^{\ast}$ and $U'_{\ext}$ fix the reducible $SO(3)$--bundle $\underline{\R} \oplus H^{2k}$, where $H$ is the radial extension of the Hopf line bundle. On this bundle we consider the reducible monopole induced by an Euclidean Dirac monopole of charge $2k$, singularity at the origin and mass $\lambda >0$. We denote by $c_{B^{\ast}}$ and $c_{U'_{\ext}}$ such pair regarded as a configuration on $B^{\ast}$ and $U'_{\ext}$, respectively.

Over $B$ consider the trivial bundle $B \times \Lie{su}_{2}$ and a pair $c_{B}$ defined as follows. Start with an Euclidean $SU(2)$ monopole of charge $k$, mass $\lambda >0$ and centre at the origin and the induced $SO(3)$--monopole on the adjoint bundle. Using a cut-off function $\chi$ with $\chi \equiv 1$ on $\frac{1}{4}B$ and $\chi \equiv 0$ outside of $\frac{1}{2}B$, we modify this initial configuration to define $c_{B}$ so that it coincides with $c_{B^{\ast}}$ on the annulus $B \setminus \frac{1}{2}B$. In order to carry out this step, it is necessary to fix an isomorphism $\eta \co B^{\ast} \times \Lie{su}_{2} \ra \underline{\R} \oplus H^{2k}$.

Finally, over $U_{\ext}$ consider the reducible bundle $\underline{\R} \oplus M$, where $M=L_{v,b} \otimes L^{2k}_{q} \otimes \bigotimes _{i=1}^{n}{ L_{p_{i}}^{-1} }$ is endowed with the corresponding sum of periodic Dirac monopoles. Here we choose $v$ so that $v+ka_{0}-\sum _{i=1}^{n}{ G_{p_{i}}(q) } = \lambda$. Furthermore, using a cut-off function one can modify this initial configuration to define a pair $c_{U_{\ext}}$ that agrees with $c_{B^{\ast}}$ on $B \setminus \frac{1}{2}B$.

Using the isomorphism $\eta$, we can now define pairs $c_{X}, c_{Y}, c_{X^{\ast}}, c_{Y^{\ast}}$. The former two pairs are smooth configurations on $X$ and $Y$, respectively; up to the modification appearing in the definition of $c_{U_{\ext}}$, the latter two pairs are (sums of) Dirac monopoles on $X^{\ast}$ and $Y^{\ast}$, respectively.

\subsection{Weighted spaces and the Fredholm property}

To each of the pairs $c_{X}, c_{Y}, c_{X^{\ast}}, c_{Y^{\ast}}$ we associate the corresponding Dirac operator $D_{X}, D_{Y}, D_{X^{\ast}}, D_{Y^{\ast}}$ acting on section of the adjoint bundle. We now introduce weighted Sobolev spaces and prove that these operators extend to Fredholm operators between these spaces. We begin with $\R ^{3}$.

\begin{definition}\label{def:Weighted:Spaces:Y}
Let $\rho$ be the distance from the origin in $\R ^{3}$.
\begin{itemize}
\item[(i)] Define $W^{1,2}(Y)$ to be the closure of the space of smooth compactly supported sections $\Omega (\R ^{3};\Lie{su}_{2})$ with respect to the norm
\[
\| \xi \| ^{2}_{W^{1,2}(Y)} = \| \xi \| ^{2}_{L^{2}} + \int_{Y}{ (1+\rho ^{2})\left( |\nabla _{A}\xi|^{2} + |[\Phi,\xi]|^{2} \right) }.
\]
\item[(ii)] We say that a $1$--form $f$ with values in the trivial $SO(3)$--bundle over $Y=\R ^{3}$ is in the space $L^{2}(Y)$ if $\sqrt{1+\rho ^{2}}f \in L^{2}$.
\item[(iii)] $W^{1,2}(Y^{\ast})$ is the closure of the space of compactly supported smooth forms with values in $\underline{\R} \oplus H^{2k}$ over $Y^{\ast}$ with respect to the norm
\[
\| \xi \| ^{2}_{W^{1,2}(Y^{\ast})} = \| \xi \| ^{2}_{L^{2}} + \int_{Y}{ \rho ^{2}\left( |\nabla _{A}\xi|^{2} + |[\Phi,\xi]|^{2} \right) }.
\]
\item[(iv)] A $1$--form $f$ with values in $\underline{\R} \oplus H^{2k}$ is in $L^{2}(Y^{\ast})$ if and only if $\rho f \in L^{2}$.
\end{itemize}
In (i) and (iii) $(A,\Phi) = c_{Y}$ and $c_{Y^{\ast}}$, respectively.
\end{definition}

\begin{prop}\label{prop:Fredholm:Y}
The operator $D_{Y} \co W^{1,2}(Y) \ra L^{2}(Y)$ is Fredholm.
\proof
This is essentially Proposition 7.2 in \cite{Taubes:Stability:Yang:Mills:Theories}, but we give an overview of the proof since, to be better suited to the presence of singularities, our spaces are slightly different from the one used by Taubes.

First we show that $D_{Y}$ has finite dimensional kernel and closed range. Both statement follow once we show that there exists $C>0$ and a compact set $K \subset Y$ such that for all $\xi \in W^{1,2}(Y)$
\begin{equation}\label{eqn:Estimate:Fredholm:Y}
\| \xi \| _{W^{1,2}(Y)} \leq C \left( \| D\xi \| _{L^{2}(Y)} + \| \xi |_{K} \| _{L^{2}} \right).
\end{equation}
As a preliminary, we claim that there exists $C>0$ such that
\begin{equation}\label{eqn:Apriori:Estimate:Y}
\| \xi \| _{W^{1,2}(Y)} \leq C \left( \| D\xi \| _{L^{2}(Y)} + \| \xi \| _{L^{2}} \right)
\end{equation}
for all $\xi \in C^{\infty}_{0}$. Indeed, use the Weitzenb\"ock formula for $D^{\ast}D$ to derive the Bochner-type identity
\[
\frac{1}{2}d^{\ast}d(|\xi|^{2}) = \langle D^{\ast}D\xi, \xi \rangle - |\nabla_{A}\xi|^{2} - |[\Phi, \xi]|^{2} - \langle (\ast F_{A} + d_{A}\Phi) \cdot \xi, \xi \rangle .
\]
Integrating by parts against $1+\rho ^{2}$ yields
\begin{align*}
\int{ (1+\rho ^{2}) \left( |\nabla_{A}\xi|^{2} + |[\Phi, \xi]|^{2} \right) } &= 2\int{ \rho \langle D\xi , d\rho \cdot \xi \rangle } + \int{ (1+\rho ^{2}) |D\xi|^{2} } \\
&{}+ \int{ (1+\rho ^{2})  \langle (\ast F_{A} + d_{A}\Phi) \cdot \xi, \xi \rangle } + 3\int{|\xi|^{2}}\\
& {}\leq 2 \int{ (1+\rho ^{2}) |D\xi|^{2} } + C_{1} \int{ |\xi|^{2} }.
\end{align*}
The constant $C_{1}=4+\| (1+\rho ^{2}) \left( \ast F_{A} + d_{A}\Phi \right) \| _{L^{\infty}}$ is bounded because $|d_{A}\Phi|=|F_{A}|=O(\rho ^{-2})$.

Thus \eqref{eqn:Apriori:Estimate:Y} is proved and we proceed with the proof of \eqref{eqn:Estimate:Fredholm:Y}. Let $R>0$ be sufficiently large and fix a cut-off function $\chi$ with $\chi \equiv 1$ on $Y \setminus B_{R}$ and $\chi \equiv 0$ on $B_{R-1}$. Write $\xi = \xi _{1}+\xi_{2} = \chi \xi + (1-\chi)\xi$. Applying \eqref{eqn:Apriori:Estimate:Y} to $\xi_{2}$ and observing that $\| D\xi_{i} \| _{L^{2}(Y)} \leq C\left( \| \xi |_{B_{R}} \| _{L^{2}} + \| D\xi \| _{L^{2}(Y)} \right)$, we reduce to prove
\[
\| \xi_{1} \| _{W^{1,2}(Y)} \leq C \left( \| D\xi_{1} \| _{L^{2}(Y)} + \| \xi_{1} |_{K} \| _{L^{2}} \right)
\]
for some $C>0$ and compact set $K \subset Y$ independent of $\xi _{1}$.

Since the pair $(A,\Phi)$ is reducible on the support of $\xi _{1}$ we decompose into diagonal and off-diagonal part $\xi _{1}= \xi _{1,D}+\xi _{1,T}$. On the off-diagonal part the estimate follows from \eqref{eqn:Apriori:Estimate:Y} provided $R$ is sufficiently large. Indeed, since $\lim_{\rho \ra \infty}{|\Phi|} = \lambda$, we can choose $R$ so that $\| \Phi \| \geq \frac{\lambda}{2}$ when $\rho \geq R$. Then
\[
\int{ (1+\rho ^{2})|[\Phi, \xi _{1,T}]|^{2} } \geq (1+R^{2})\frac{\lambda}{2} \int{ |\xi _{1,T}|^{2} }.
\]
Choosing $R$ even larger if necessary, the term $\| \xi \| _{L^{2}}$ can be absorbed in the left-hand-side of \eqref{eqn:Apriori:Estimate:Y} to obtain \eqref{eqn:Estimate:Fredholm:Y}. On the diagonal part we can appeal to standard theory of weighted Sobolev spaces for the scalar Laplacian on $\R ^{3}$ as in (1) in the proof of Proposition \ref{prop:Deformation:Operator:Fredholm}.

Finally, we claim that the cokernel of $D=D_{Y}$ is finite dimensional. By duality in weighted Sobolev spaces, this cokernel is identified with the kernel of $D^{\ast}$ in $\sqrt{1+\rho^{2}} L^{2}$ and its finite dimensionality follows from the a priori estimate
\begin{equation}\label{eqn:Estimate:Cokernel:Y}
\| (1+\rho ^{2})^{-\frac{1}{2}}u \| ^{2}_{L^{2}} + \| \nabla _{A}u \| ^{2}_{L^{2}} + \| [\Phi, u] \| ^{2}_{L^{2}} \leq C \left( \| D^{\ast}u \| ^{2}_{L^{2}} + \| u|_{2B \setminus B} \| ^{2}_{L^{2}} \right).
\end{equation}
by standard arguments. As \eqref{eqn:Apriori:Estimate:Y}, this estimate is obtained integrating by parts the Weitzenb\"ock formula for $DD^{\ast}$ in Lemma \ref{lem:Weitzenbock}. The constant $C$ depends on $\| \ast F_{A} -d_{A}\Phi \| _{L^{\infty}}$ which is supported on the annulus $2B \setminus B$.
\endproof
\end{prop}

\begin{prop}\label{prop:Fredholm:Y:ast}
The operator $D_{Y^{\ast}} \co W^{1,2}(Y^{\ast}) \ra L^{2}(Y^{\ast})$ is Fredholm.
\proof
It is immediate to check that $D=D_{Y^{\ast}}$ is surjective. Indeed, it follows from the Weitzenb\"ock formula for $DD^{\ast}$ and the fact that $c_{Y^{\ast}}$ is a solution to the Bogomolny equation that elements in the kernel of $D^{\ast}$ are constant diagonal sections in the decomposition $\underline{\R} \oplus H^{2k}$ and these are not in the space $\rho ^{-1}L^{2}$ dual to $L^{2}(Y^{\ast})$.

As before, it remains to prove the existence of a constant $C>0$ and a compact set $K \subset \R ^{3} \setminus \{ 0 \}$ such that for all $\xi \in W^{1,2}(Y^{\ast})$
\[
\| \xi \| _{W^{1,2}(Y^{\ast})} \leq C \left( \| D\xi \| _{L^{2}(Y^{\ast})} + \| \xi |_{K} \| _{L^{2}} \right).
\]
This is obtained as in Proposition \ref{prop:Fredholm:Y}. First, observe that an estimate analogous to \eqref{eqn:Apriori:Estimate:Y} holds because $\rho ^{2}|F_{A}|=\rho ^{2}|d_{A}\Phi|=k$ everywhere on $\R ^{3} \setminus \{ 0 \}$. In view of the proof of Proposition \ref{prop:Fredholm:Y}, we only have to explain why there exists $\sigma >0$ sufficiently small and $C>0$ such that
\begin{equation}\label{eqn:Estimate:Fredholm:Y:ast}
\| \xi \| _{W^{1,2}(Y^{\ast})} \leq C \left( \| D\xi \| _{L^{2}(Y^{\ast})} + \| \xi |_{B \setminus B_{\sigma }} \| _{L^{2}} \right)
\end{equation}
for all $\xi \in C^{\infty}_{0}(B^{\ast})$.

Since $c_{Y^{\ast}}$ is reducible, we decompose into diagonal and off-diagonal part. On the diagonal part one can argue as in (1) in the proof of Proposition \ref{prop:Deformation:Operator:Fredholm} to deduce \eqref{eqn:Estimate:Fredholm:Y:ast} from the theory of weighted Sobolev spaces for the Laplacian on $\R ^{3}$. On the off-diagonal part, if $\lambda$ vanished we could deduce \eqref{eqn:Estimate:Fredholm:Y:ast} from the arguments of (3) in the proof of Proposition \ref {prop:Deformation:Operator:Fredholm}. Since $\lambda >0$ yields a lower order term in the equation, one can then show that there exists $\sigma = \sigma (\lambda)$ such that
\[
\| \xi \| ^{2}_{L^{2}(B^{\ast}_{\sigma})} \leq C \left( \| D\xi \| ^{2}_{L^{2}} + \| \xi \| ^{2}_{L^{2}(B \setminus B_{\sigma})} \right). \qedhere
\]
\end{prop}

\begin{definition}\label{def:Weighted:Spaces:X}
Fix $\delta >0$ with $\delta < \min{\{ \tfrac{1}{2}, \delta _{0} \} }$ where $\delta _{0}$ is given by Proposition \ref{prop:Laplacian:Infinity}.
\begin{itemize}
\item[(i)] Set $W^{1,2}(X)=W^{1,2}_{\underline{\delta}-1}$ and $L^{2}(X)=L^{2}_{\underline{\delta}-2}$, with $W^{1,2}_{\underline{\delta}-1}$ and $L^{2}_{\underline{\delta}-2}$ defined in Definition \ref{def:Weighted:Spaces:Global}.
\item[(ii)] On $X^{\ast}$ we define $W^{1,2}(X^{\ast})$ to be the closure of the smooth compactly supported forms with values in $\underline{\R} \oplus M$ with respect to the norm defined by the maximum of the semi-norms:
\[
\| \xi|_{U_{\ext}} \|_{W^{1,2}(X)} \qquad \| \xi |_{B^{\ast}} \| _{W^{1,2}(Y^{\ast})}
\]
$L^{2}(X^{\ast})$ is defined similarly using the $L^{2}(X)$--norm on $U_{\ext}$ and the $L^{2}(Y^{\ast})$--norm on $B^{\ast}$.
\end{itemize}
\end{definition}

By Proposition \ref{prop:Deformation:Operator:Fredholm} $D_{X} \co W^{1,2}(X) \ra L^{2}(X)$ is a Fredholm operator. Combining the proof of Propositions \ref{prop:Deformation:Operator:Fredholm} and Proposition \ref{prop:Fredholm:Y:ast} we deduce the analogous statement for $D_{X^{\ast\ast}}$.

\subsection{Application of the excision principle}

The next step in the proof of the index formula Theorem \ref{thm:Index:Theorem} is an application of the excision principle for the index of Fredholm operators. A proof of the excision principle in the non-compact setting whic immediately applies to the operators $D_{X},D_{Y},D_{X^{\ast}},D_{Y^{\ast}}$ is given by Charbonneau in \cite[Appendix B]{Charbonneau:Thesis}.

\begin{lemma}\label{lem:Excision:principle}
$\ind{(D_{X^{\ast}})} + \ind{(D_{Y})}=\ind{(D_{X})} + \ind{ (D_{Y^{\ast}}) }$.
\end{lemma}

We apply the lemma to prove Theorem \ref{thm:Index:Theorem} by computing $\ind{(D_{Y})}$ and $\ind{(D_{X^{\ast}})} - \ind{ (D_{Y^{\ast}}) }$.

By \cite[Proposition 9.1]{Taubes:Stability:Yang:Mills:Theories} $\ind{(D_{Y})} = 4k$. Indeed, \eqref{eqn:Estimate:Cokernel:Y} shows that if $(1+\rho ^{2})^{-\frac{1}{2}}\xi \in L^{2}$ and $D^{\ast}\xi =0$, then $\nabla _{A}\xi, [\Phi, \xi] \in L^{2}$, \ie $\xi \in H_{c}$ in Taubes's notation. Then, by duality in weighted Sobolev spaces, $\ind{(D_{Y})} = -\overline{i}(\mathcal{D}_{c}^{\ast})=4k$ in the notation of \cite{Taubes:Stability:Yang:Mills:Theories}.

As for the indices $\ind{(D_{X^{\ast}})}, \ind{(D_{Y^{\ast}})}$, since $c_{X^{\ast}}$ and $c_{Y^{\ast}}$ are both reducible, we decompose the problem into diagonal and off-diagonal part. By Definitions \ref{def:Weighted:Spaces:Y} and \ref{def:Weighted:Spaces:X}, it is easy to see that $D_{X^{\ast}}$ acting on the diagonal component is injective but has a $4$--dimensional cokernel (dual to the subspace spanned by constant $0$ and $1$--forms) and that $D_{Y^{\ast}}$ is an isomorphism when acting on the diagonal components.

It does not seem immediate to calculate the index of $D_{X^{\ast}}$ and $D_{Y^{\ast}}$ acting on off-diagonal components individually, even if one guesses that they both vanish. However, we are only interested in the difference between the two indices and we are going to prove that this is zero. More precisely, we will show that taking $\lambda$ sufficiently large, we can make sure that
\begin{itemize}
\item[(i)] the two operators are surjective;
\item[(ii)]  elements in their kernels are concentrated in a small neighbourhood of $q$ and $0$, respectively. It follows that the kernels are isomorphic.
\end{itemize}

Before embarking in the proof of (i) and (ii) we explain why the index of $D_{X^{\ast}}$ and $D_{Y^{\ast}}$ does not depend on the mass $\lambda$. On $\R ^{3}$ this is clear by scaling. On $\R ^{2} \times \Sph ^{1}$ we consider the continuous family of Fredholm operators $D_{v'}\co W^{1,2}(X^{\ast}) \ra L^{2}(X^{\ast})$ defined as follows. Fix positive constants $C$ and $R$ such that $|\Phi| \geq C$ when $r \geq R$. If $n <2k$ we can take $C$ arbitrarily large; if $n=2k$ we have to assume that $C<v$. Define $D_{v'} = D_{X^{\ast}} + [\psi_{v'}, \cdot]$, where $\psi _{v'} = \chi (v'-v) \hat{\sigma}$. Here $\hat{\sigma} = \frac{\Phi}{|\Phi|}$ and $\chi$ is a smooth function $\chi \equiv 0$ when $r \leq R$ and $\chi \equiv 1$ if $r \geq 2R$. $D_{v'} \co W^{1,2}(X^{\ast}) \ra L^{2}(X^{\ast})$ is a bounded Fredholm operator for all $v' \in (v-C, v+C)$. In particular, the index of $D_{X^{\ast}}$ is independent of the mass $v \in \R$, subject to the only constraint $v >0$ when $n=2k$.

\subsection{The index of the Dirac operator $D$ twisted by a Dirac monopole}

In the rest of the proof we will assume that $\lambda$ is as large as needed. We want to prove that $D_{X^{\ast}} \co W^{1,2}(X^{\ast}) \ra L^{2}(X^{\ast})$ and $D_{Y^{\ast}} \co W^{1,2}(Y^{\ast}) \ra L^{2}(Y^{\ast})$ acting on the off-diagonal component (i) are surjective and (ii) have isomorphic kernel whenever $\lambda$ is sufficiently large. The proof is modelled on \cite[\S 7.1.2]{Donaldson:Kronheimer} and the main technical ingredient is to exhibit right inverses $Q_{X^{\ast}} \co L^{2}(X^{\ast}) \ra W^{1,2}(X^{\ast})$ and $Q_{Y^{\ast}} \co L^{2}(Y^{\ast}) \ra W^{1,2}(Y^{\ast})$ of $D_{X^{\ast}}$ and $D_{Y}^{\ast}$, respectively, which are bounded independently of $\lambda$.

We collect some important properties of $c_{Y^{\ast}}$ and $c_{X^{\ast}}$.
\begin{itemize}
\item[(a)] $c_{Y^{\ast}}$ is an exact solution to the Bogomolny equation. On the other hand, if $(A,\Phi) = c_{X^{\ast}}$, $\Psi=\ast F_{A} -d_{A}\Phi$ is supported in the region $2B \setminus B$ and $|\ast F_{A} -d_{A}\Phi| \leq C$ for a constant independent of $\lambda$.
\item[(b)] In both cases, there exists $\lambda_{0}$ such that if $\lambda >\lambda _{0}$ then $|\Phi| \geq \frac{\lambda}{2}$ outside of $\frac{1}{2}B$. On $Y^{\ast}$ this is clear because $\Phi = \left( \lambda -\frac{k}{\rho} \right) \hat{\sigma}$. On $X^{\ast}$ the statement is true for the sum of periodic Dirac monopoles provided $\lambda$ is sufficiently large (this follows from the maximum principle). In particular, $|\Phi| \geq \frac{\lambda}{2}$ on $B$ and outside of $2B$. In the annulus $2B \setminus B$, $\Phi = \left( \lambda -\frac{k}{\rho} \right) \hat{\sigma} + O(\rho)$ and therefore, taking $\lambda$ even larger if necessary, $|\Phi| \geq \frac{\lambda}{2}$.
\end{itemize}

\begin{lemma}\label{lem:Right:Inverse:X:ast}
There exists $\lambda _{0} \geq 2$ and $C>0$ such that if $\lambda > \lambda _{0}$ then the following holds. For all $f=f_{T} \in L^{2}(X^{\ast})$ there exists $\xi \in W^{1,2}(X^{\ast})$ such that $D_{X^{\ast}}\xi = f$ and $\| \xi \| _{W^{1,2}(X^{\ast})} \leq C \| f \| _{L^{2}(X^{\ast})}$.
\proof
We proceed in two steps. First we solve the equation $D\xi=f$ for $\xi$ of the form $\xi = D^{\ast}u$ by variational methods. Then we obtain the estimate integrating the Weitzenb\"ock formula for $D^{\ast}D$.

By the Weitzenb\"ock formula $DD^{\ast}= \nabla _{A}^{\ast}\nabla _{A} - \ad{(\Phi)}^{2} + \Psi$ and the fact that $u=u_{T}$ we deduce that $\| D^{\ast}u \| ^{2}_{L^{2}}$ is uniformly equivalent to $\| \nabla _{A}u \| ^{2}_{L^{2}} + \| [\Phi, u] \| _{L^{2}}^{2}$ provided $\lambda$ is sufficiently large. Indeed, we integrate by parts the Weitzenb\"ock formula and use the inequality
\[
\left| \int{ \langle \Psi \cdot u , u \rangle } \right| \leq \| \Psi \| _{L^{\infty}} \| u|_{\text{supp}(\Psi)} \| _{L^{2}}^{2} \leq \frac{C}{\lambda} \| [\Phi, u] \| ^{2}_{L^{2}}
\]
which follows from (a) and (b) above.

To show that $\| \nabla _{A}u \| ^{2}_{L^{2}} + \| [\Phi, u] \| _{L^{2}}^{2}$ is a norm, fix a cut-off function $\chi$ with $\chi \equiv 1$ on $\frac{1}{2}B$ and vanishing outside of $B$. Then by the Poincar\'e inequality on $B$ and the fact that $|\Phi| \geq 1$ outside of $\frac{1}{2}B$, we deduce
\[
\| u \|^{2}_{L^{2}} \leq \| \chi u \| _{L^{2}}^{2} + \| (1-\chi)u \| _{L^{2}}^{2} \leq C \left( \| \nabla _{A}u \| ^{2}_{L^{2}} + \| [\Phi, u] \| _{L^{2}}^{2} \right) .
\]
Define $H$ to be closure of smooth compactly supported forms $u$ with values in the line bundle $M$ with respect to the norm $\| \nabla _{A}u \| ^{2}_{L^{2}} + \| [\Phi, u] \| _{L^{2}}^{2}$. Hardy's inequality $\| \rho ^{-1}u \| _{L^{2}} \leq 2\| \nabla _{A} u \| _{L^{2}}$ on $B$ and the fact that $|\Phi| \geq c_{1}\rho _{i}^{-2}$ in a neighbourhood of the singularity $p_{i}$ imply that $\langle f , u \rangle _{L^{2}}$ defines a continuous functional on $H$ for all $f \in L^{2}(X^{\ast})$. Then a weak solution $\xi$ to $D\xi = f$ with $\| \xi \| _{L^{2}} \leq C \| f \| _{L^{2}(X^{\ast})}$ is found by minimising the functional $\frac{1}{2}\int{|D^{\ast}u|^{2}} - \langle f,u \rangle _{L^{2}}$ on $H$.

It remains to show the existence of a uniform constant $C>0$ such that
\begin{equation}\label{eqn:A:Priori:Estimate:Right:Inverse}
\| \xi \| _{W^{1,2}(X^{\ast})} \leq C \left( \| D\xi \| _{L^{2}(X^{\ast})} + \| \xi \| _{L^{2}} \right).
\end{equation}
Fix $\sigma >0$ such that $2B$ and $B_{2\sigma}(p_{i})$ are all disjoint. Observe that given a smooth compactly supported function $\chi$ on $X^{\ast}$ then $\| D(\chi\xi) \| _{L^{2}(X^{\ast})} \leq C_{1}(\chi) \| \xi \| _{L^{2}} + C_{2} \| D\xi \| _{L^{2}(X^{\ast})}$. Hence to prove \eqref{eqn:A:Priori:Estimate:Right:Inverse} we can suppose that $\xi$ is supported on either $B$, $B_{2\sigma}(p_{i})$ or the complement $U_{\sigma}$ of $\tfrac{1}{2}B \cup \bigcup _{i=1}^{n}{ B_{\sigma _{0}}(p_{i})} $.

\begin{enumerate}
\item If $\xi \in C^{\infty}_{0}(B)$, as in \eqref{eqn:Apriori:Estimate:Y} we can find a constant $C$ depending only on $\| \rho ^{2} \left( \ast F_{A}+d_{A}\Phi \right) \| _{L^{\infty}}$ such that
\[
\int{ \rho ^{2} \left( |\nabla _{A}\xi|^{2} + |[\Phi, \xi]|^{2} \right) } \leq C \left( \| \rho D\xi \| _{L^{2}} + \| \xi \| _{L^{2}} \right).
\]
Since $(A,\Phi)$ is an Euclidean Dirac monopole up to terms of order $O(\rho)$, $\rho ^{2} \left( \ast F_{A}+d_{A}\Phi \right)$ is bounded independently of $\lambda$.

\item If $\xi \in C^{\infty}_{0}(U_{\sigma})$ the estimate was proved in Steps 3 and 4 in the proof of Proposition \ref{prop:Laplacian:Infinity}. The constant is uniform because $\omega d_{A}\Phi$ is bounded independently of $\lambda$.

\item Finally, when $\xi$ is compactly supported on $B_{2\sigma}(p_{i}) \setminus \{ p_{i} \}$ we have to slightly modify the arguments of (3) in the proof of Proposition \ref{prop:Deformation:Operator:Fredholm} to show that \eqref{eqn:A:Priori:Estimate:Right:Inverse} holds with a uniform constant $C$.

First, an integration by parts of the Weitzenb\"ock formula yields
\[
\int{ \rho _{i}^{-2\delta +1} \left( |\nabla _{A}\xi|^{2} + |[\Phi, \xi]| ^{2} \right) } \leq C \left( \| \rho _{i}^{-\delta +\frac{1}{2}}D\xi \| ^{2}_{L^{2}} + \| \rho _{i}^{-\delta +\frac{1}{2}}\xi \| ^{2}_{L^{2}} \right)
\]
for a constant $C$ depending only on $\| \rho_{i}^{2}d_{A}\Phi \| _{L^{\infty}}$ and therefore independent of $\lambda$. It remains to control the weighted norm $\| \rho _{i}^{-\delta +\frac{1}{2}}\xi \| _{L^{2}}$.

Suppose first that $(A,\Phi)$ coincides with an Euclidean Dirac monopole with $|\Phi| = \lambda +\frac{1}{2\rho_{i}}$. Set $h = |\Phi|$ and observe that $D(h^{-1}D^{\ast}\xi) = D^{\ast}(h^{-1}D\xi)$ because $d_{A}(h^{-1}\Phi) = 0$. Then the identities $D^{\ast}\xi -D\xi = -2[\Phi, \xi]$ and $4|[\Phi,\xi]|^{2} = |D^{\ast}\xi|^{2}-|D\xi|^{2} + 4\langle D\xi, [\Phi, \xi] \rangle$ imply
\[
\| \rho _{i}^{-\delta}h^{-\frac{1}{2}}[\Phi,\xi]\|_{L^{2}} \leq \delta \| \rho_{i} ^{-\delta -1}h^{-\frac{1}{2}}\xi \| _{L^{2}} + \| \rho_{i} ^{-\delta}h^{-\frac{1}{2}}D\xi \| _{L^{2}}.
\]
Since $\rho_{i} ^{-2\delta}h^{-1}|[\Phi,\xi]|^{2} \geq \frac{1}{4}\rho_{i} ^{-2\delta-2}h^{-1}|\xi|^{2}$ and $\delta <\frac{1}{2}$ we deduce
\[
\| \rho_{i} ^{-\delta -\frac{1}{2}}\xi \| _{L^{2}} \leq \frac{2}{1-2\delta} \| \rho _{i}^{-\delta +\frac{1}{2}}D\xi \| _{L^{2}}.
\]

Now, the pair $(A,\Phi)$ agrees with an Euclidean Dirac monopole with $|\Phi| = \lambda +\frac{1}{2\rho}$ up to bounded terms which only depend on $p_{1}, \ldots, p_{n}$ and $q$. Therefore we can find $\sigma >\sigma '>0$ and $C>0$ independent of $\lambda$ such that
\[
\int{ \rho _{i}^{-2\delta +1} \left( |\nabla _{A}\xi|^{2} + |[\Phi, \xi]| ^{2} \right) } \leq C \left( \| \rho _{i}^{-\delta +\frac{1}{2}}D\xi \| ^{2}_{L^{2}} + \| \xi \| ^{2}_{L^{2}(\sigma ' \leq \rho _{i} \leq \sigma )} \right).
\]
\end{enumerate}
Combining (1),(2) and (3) we obtain \eqref{eqn:A:Priori:Estimate:Right:Inverse} and the Lemma is proved.
\endproof
\end{lemma}

The analogous statement for $D_{Y^{\ast}}$ is actually easier to prove: the existence of a week solution follows immediately from the Hardy inequality, while the analogous of \eqref{eqn:A:Priori:Estimate:Right:Inverse} is \eqref{eqn:Apriori:Estimate:Y} with $\rho$ in place of $\sqrt{1+\rho ^{2}}$. The existence of uniformly bounded right inverses $Q_{X^{\ast}}$ and $Q_{Y^{\ast}}$ follows.

In order to conclude the proof of Theorem \ref{thm:Index:Theorem}, we have to show that the kernels of $D_{X^{\ast}}$ and $D_{Y^{\ast}}$ are isomorphic.

Fix a pair of cut-off functions $\gamma$ and $\beta$ with $\gamma \equiv 1$ on $\frac{1}{2}B$, $\gamma = 0$ outside of $B$, $\beta \equiv 1$ outside of $B$, $\beta \equiv 0$ in $\frac{1}{2}B$ and $\beta \equiv 1$ on the support of $d\gamma$. Notice that $\| (1-\gamma) \xi \| _{L^{2}} \leq \frac{c}{\lambda _{0}} \| \xi \| _{W^{1,2}(X^{\ast})}$ if $\lambda > \lambda _{0}$. Therefore, taking $\lambda _{0}$ larger if necessary, we can improve the estimate in Lemma \ref{lem:Right:Inverse:X:ast} to
\[
\| \xi \| _{W^{1,2}(X^{\ast})} \leq C \left( \| D\xi \| _{L^{2}(X^{\ast})} + \| \gamma\xi \| _{L^{2}} \right)
\]
(and similarly on $Y^{\ast}$).

Now, let $\xi \in W^{1,2}(X^{\ast})$ be such that $D_{X^{\ast}}\xi = 0$. We define an element in the kernel of $D_{Y^{\ast}}$ by $\xi' = \gamma \xi - Q_{Y^{\ast}}D_{Y^{\ast}}(\gamma\xi)$. We want to show that the map $\xi \mapsto \xi'$ is injective. By contradiction, suppose that there exists $\xi$ such that $\gamma \xi = Q_{Y^{\ast}}D_{Y^{\ast}}(\gamma\xi)$. Then
\[
\| \gamma \xi \| _{L^{2}} \leq \| \gamma \xi \| _{W^{1,2}(Y^{\ast})} \leq C \| D_{Y^{\ast}}(\gamma \xi) \| _{L^{2}(Y^{\ast})} \leq C \| \beta \xi \| _{L^{2}} \leq \frac{C}{\lambda} \| \xi \| _{W^{1,2}(X^{\ast})} \leq \frac{C}{\lambda} \| \gamma \xi \| _{L^{2}}.
\]
If $\lambda$ is sufficiently large we get a contradiction. Exchanging the role of $X^{\ast}$ and $Y^{\ast}$, we construct injective maps between the kernel of $D_{X^{\ast}}$ and $D_{Y^{\ast}}$, which are therefore isomorphic.

\bibliography{Monopoles}

\begin{thebibliography}{10}

\bibitem{Amrouche:Girault:Giroire:1}
C.~Amrouche, V.~Girault, and J.~Giroire.
\newblock Weighted {S}obolev spaces for {L}aplace's equation in {$\bold R^n$}.
\newblock {\em J. Math. Pures Appl. (9)}, 73(6):579--606, 1994.

\bibitem{Amrouche:Girault:Giroire}
C.~Amrouche, V.~Girault, and J.~Giroire.
\newblock Dirichlet and {N}eumann exterior problems for the {$n$}-dimensional
  {L}aplace operator: an approach in weighted {S}obolev spaces.
\newblock {\em J. Math. Pures Appl. (9)}, 76(1):55--81, 1997.

\bibitem{Anghel}
N.~Anghel.
\newblock On the index of {C}allias-type operators.
\newblock {\em Geom. Funct. Anal.}, 3(5):431--438, 1993.

\bibitem{Atiyah:Hyperbolic:Monopoles}
M.~F. Atiyah.
\newblock Magnetic monopoles in hyperbolic spaces.
\newblock In {\em Vector bundles on algebraic varieties ({B}ombay, 1984)},
  volume~11 of {\em Tata Inst. Fund. Res. Stud. Math.}, pages 1--33. Tata Inst.
  Fund. Res., Bombay, 1987.

\bibitem{Atiyah:Hitchin}
Michael Atiyah and Nigel Hitchin.
\newblock {\em The geometry and dynamics of magnetic monopoles}.
\newblock M. B. Porter Lectures. Princeton University Press, Princeton, NJ,
  1988.

\bibitem{Biquard:Paraboliques}
Olivier Biquard.
\newblock Fibr\'es paraboliques stables et connexions singuli\`eres plates.
\newblock {\em Bull. Soc. Math. France}, 119(2):231--257, 1991.

\bibitem{Biquard:Prolongement}
Olivier Biquard.
\newblock Prolongement d'un fibre holomorphe hermitien \`a\ courbure {$L^p$}
  sur une courbe ouverte.
\newblock {\em Internat. J. Math.}, 3(4):441--453, 1992.

\bibitem{Biquard:Boalch}
Olivier Biquard and Philip Boalch.
\newblock Wild non-abelian {H}odge theory on curves.
\newblock {\em Compos. Math.}, 140(1):179--204, 2004.

\bibitem{Biquard:Jardim}
Olivier Biquard and Marcos Jardim.
\newblock Asymptotic behaviour and the moduli space of doubly-periodic
  instantons.
\newblock {\em J. Eur. Math. Soc. (JEMS)}, 3(4):335--375, 2001.

\bibitem{Braam:Donaldson}
P.~J. Braam and S.~K. Donaldson.
\newblock Floer's work on instanton homology, knots and surgery.
\newblock In {\em The {F}loer memorial volume}, volume 133 of {\em Progr.
  Math.}, pages 195--256. Birkh\"auser, Basel, 1995.

\bibitem{Braam}
Peter~J. Braam.
\newblock Magnetic monopoles on three-manifolds.
\newblock {\em J. Differential Geom.}, 30(2):425--464, 1989.

\bibitem{Callias}
Constantine Callias.
\newblock Axial anomalies and index theorems on open spaces.
\newblock {\em Comm. Math. Phys.}, 62(3):213--234, 1978.

\bibitem{Charbonneau:Thesis}
Benoit Charbonneau.
\newblock {\em Analytic aspects of periodic instantons}.
\newblock 2004.
\newblock Ph.D. Thesis, MIT.

\bibitem{Charbonneau:Hurtubise}
Benoit Charbonneau and Jacques Hurtubise.
\newblock Singular {H}ermitian-{E}instein monopoles on the product of a circle
  and a {R}iemann surface.
\newblock {\em Int. Math. Res. Not. IMRN}, (1):175--216, 2011.

\bibitem{Cherkis:Talk}
Sergey~A. Cherkis.
\newblock Doubly {P}eriodic {M}onopoles and their {M}oduli {S}paces.
\newblock Talk given at the workshop Advances in hyperk\"ahler and holomorphic
  symplectic geometry, BIRS, Banff, March 2012. {A}vailable online at
  http://temple.birs.ca/\textasciitilde 12w5126/Cherkis.pdf.

\bibitem{Cherkis:Kapustin:1}
Sergey~A. Cherkis and Anton Kapustin.
\newblock Nahm transform for periodic monopoles and {$N=2$} super
  {Y}ang-{M}ills theory.
\newblock {\em Comm. Math. Phys.}, 218(2):333--371, 2001.

\bibitem{Cherkis:Kapustin:3}
Sergey~A. Cherkis and Anton Kapustin.
\newblock Hyper-{K}\"ahler metrics from periodic monopoles.
\newblock {\em Phys. Rev. D (3)}, 65(8):084015, 10, 2002.

\bibitem{Cherkis:Kapustin:2}
Sergey~A. Cherkis and Anton Kapustin.
\newblock Periodic monopoles with singularities and {$N=2$} super-{QCD}.
\newblock {\em Comm. Math. Phys.}, 234(1):1--35, 2003.

\bibitem{Donaldson:Floer:Homology}
S.~K. Donaldson.
\newblock {\em Floer homology groups in {Y}ang-{M}ills theory}, volume 147 of
  {\em Cambridge Tracts in Mathematics}.
\newblock Cambridge University Press, Cambridge, 2002.
\newblock With the assistance of M. Furuta and D. Kotschick.

\bibitem{Donaldson:Kronheimer}
S.~K. Donaldson and P.~B. Kronheimer.
\newblock {\em The geometry of four-manifolds}.
\newblock Oxford Mathematical Monographs. The Clarendon Press Oxford University
  Press, New York, 1990.
\newblock Oxford Science Publications.

\bibitem{Evans}
Lawrence~C. Evans.
\newblock {\em Partial differential equations}, volume~19 of {\em Graduate
  Studies in Mathematics}.
\newblock American Mathematical Society, Providence, RI, 1998.

\bibitem{Floer:Monopoles:1}
A.~Floer.
\newblock The configuration space of {Y}ang-{M}ills-{H}iggs theory of
  asymptotically flat manifolds.
\newblock In {\em The {F}loer memorial volume}, volume 133 of {\em Progr.
  Math.}, pages 43--75. Birkh\"auser, Basel, 1995.

\bibitem{Floer:Monopoles:2}
A.~Floer.
\newblock Monopoles on asymptotically flat manifolds.
\newblock In {\em The {F}loer memorial volume}, volume 133 of {\em Progr.
  Math.}, pages 3--41. Birkh\"auser, Basel, 1995.

\bibitem{Gluing}
Lorenzo Foscolo.
\newblock A gluing construction for periodic monopoles.
\newblock arXiv preprint, 2014.

\bibitem{Gibbons:Hawking}
G.~W. Gibbons and S.~W. Hawking.
\newblock Gravitational multi-instantons.
\newblock {\em Phys. Letters B}, 78(4):430--432, 1978.

\bibitem{Gilbarg:Trudinger}
David Gilbarg and Neil~S. Trudinger.
\newblock {\em Elliptic partial differential equations of second order}.
\newblock Classics in Mathematics. Springer-Verlag, Berlin, 2001.
\newblock Reprint of the 1998 edition.

\bibitem{Gross:Wilson}
Mark Gross and P.~M.~H. Wilson.
\newblock Large complex structure limits of {$K3$} surfaces.
\newblock {\em J. Differential Geom.}, 55(3):475--546, 2000.

\bibitem{Jaffe:Taubes}
Arthur Jaffe and Clifford Taubes.
\newblock {\em Vortices and monopoles}, volume~2 of {\em Progress in Physics}.
\newblock Birkh\"auser Boston, Mass., 1980.
\newblock Structure of static gauge theories.

\bibitem{Kottke}
Christopher~N. Kottke.
\newblock Callias' index theorem and monopole deformation.
\newblock arXiv:1210.3275, 2012.

\bibitem{Kronheimer:Mrowka}
P.~B. Kronheimer and T.~S. Mrowka.
\newblock Gauge theory for embedded surfaces. {I}.
\newblock {\em Topology}, 32(4):773--826, 1993.

\bibitem{Kronheimer:Thesis}
P.B. Kronheimer.
\newblock Monopoles and {T}aub-{NUT} metrics.
\newblock M.Sc. Dissertation, Oxford, 1985.

\bibitem{Kuwabara}
Ruishi Kuwabara.
\newblock On spectra of the {L}aplacian on vector bundles.
\newblock {\em J. Math. Tokushima Univ.}, 16:1--23, 1982.

\bibitem{Lockhart:McOwen}
Robert~B. Lockhart and Robert~C. McOwen.
\newblock Elliptic differential operators on noncompact manifolds.
\newblock {\em Ann. Scuola Norm. Sup. Pisa Cl. Sci. (4)}, 12(3):409--447, 1985.

\bibitem{Melrose}
Richard~B. Melrose.
\newblock {\em The {A}tiyah-{P}atodi-{S}inger index theorem}, volume~4 of {\em
  Research Notes in Mathematics}.
\newblock A K Peters Ltd., Wellesley, MA, 1993.

\bibitem{Pacard:Connected:Sums}
Frank Pacard.
\newblock Connected sum constructions in geometry and nonlinear analysis.
\newblock Lecture notes, 2008.

\bibitem{Pauly}
Marc Pauly.
\newblock Monopole moduli spaces for compact {$3$}-manifolds.
\newblock {\em Math. Ann.}, 311(1):125--146, 1998.

\bibitem{Raade:Index}
Johan R{\aa}de.
\newblock Callias' index theorem, elliptic boundary conditions, and cutting and
  gluing.
\newblock {\em Comm. Math. Phys.}, 161(1):51--61, 1994.

\bibitem{Rade:Global:Theory}
Johan R{\aa}de.
\newblock Singular {Y}ang-{M}ills fields---global theory.
\newblock {\em Internat. J. Math.}, 5(4):491--521, 1994.

\bibitem{Rade:Local:Theory:1}
Johan R{\aa}de.
\newblock Singular {Y}ang-{M}ills fields. {L}ocal theory. {I}.
\newblock {\em J. Reine Angew. Math.}, 452:111--151, 1994.

\bibitem{Rade:Local:Theory:2}
Johan R{\aa}de.
\newblock Singular {Y}ang-{M}ills fields. {L}ocal theory. {II}.
\newblock {\em J. Reine Angew. Math.}, 456:197--219, 1994.

\bibitem{Taubes:Stability:Yang:Mills:Theories}
Clifford~Henry Taubes.
\newblock Stability in {Y}ang-{M}ills theories.
\newblock {\em Comm. Math. Phys.}, 91(2):235--263, 1983.

\bibitem{Uhlenbeck:Removable:Singularities}
Karen~K. Uhlenbeck.
\newblock Removable singularities in {Y}ang-{M}ills fields.
\newblock {\em Comm. Math. Phys.}, 83(1):11--29, 1982.

\bibitem{Whitney}
Hassler Whitney.
\newblock Topological properties of differentiable manifolds.
\newblock {\em Bull. Amer. Math. Soc.}, 43(12):785--805, 1937.

\end{thebibliography}
\bibliographystyle{plain}

\end{document}